\documentclass[11pt]{amsart} \textwidth=14.5cm \oddsidemargin=1cm
\evensidemargin=1cm \usepackage{amsmath} \usepackage{amsxtra}
\usepackage{amscd} \usepackage{amsthm} \usepackage{amsfonts}
\usepackage{amssymb} \usepackage{eucal}

\newtheorem{cor}[subsubsection]{Corollary}
\newtheorem{lem}[subsubsection]{Lemma}
\newtheorem{prop}[subsubsection]{Proposition}
\newtheorem{propconstr}[subsubsection]{Proposition-Construction}

\newtheorem{conj}[subsubsection]{Conjecture}
\newtheorem{thm}[subsubsection]{Theorem}

\theoremstyle{definition}

\theoremstyle{remark}

\newcommand{\thmref}[1]{Theorem~\ref{#1}}
\newcommand{\secref}[1]{Sect.~\ref{#1}}
\newcommand{\lemref}[1]{Lemma~\ref{#1}}
\newcommand{\propref}[1]{Proposition~\ref{#1}}
\newcommand{\corref}[1]{Corollary~\ref{#1}}

\newcommand{\nc}{\newcommand}
\nc{\renc}{\renewcommand}
\nc{\ssec}{\subsection}
\nc{\sssec}{\subsubsection}
\nc{\on}{\operatorname}

\nc\ol{\overline}
\nc\wt{\widetilde}
\nc\tboxtimes{\wt{\boxtimes}}
\nc{\alp}{\alpha}

\nc{\ZZ}{{\mathbb Z}}
\nc{\NN}{{\mathbb N}}
\nc{\CC}{{\mathbb C}}
\nc{\OO}{{\mathbb O}}
\nc{\DD}{{\mathbb D}}
\nc{\GG}{{\mathbb G}}
\renewcommand{\AA}{{\mathbb A}}
\nc{\Fq}{{\mathbb F}_q}
\nc{\Fqb}{\ol{{\mathbb F}_q}}
\nc{\Ql}{\ol{{\mathbb Q}_\ell}}
\nc{\oY}{{\overline Y}}
\nc{\kap}{{\kappa}}
\nc{\id}{\text{id}}

\nc{\Eis}{\on{Eis}}
\nc{\Aut}{\on{Aut}}
\nc{\Rep}{\on{Rep}}
\nc{\Hom}{\on{Hom}}
\nc{\Loc}{\on{Loc}}
\nc{\Pic}{\on{Pic}}
\nc{\Bun}{\on{Bun}}
\nc{\IC}{\on{IC}}
\nc{\rk}{\on{rk}}
\nc{\Sh}{\on{Sh}}
\nc{\Perv}{\on{Perv}}
\nc{\Conv}{\on{Conv}}
\nc{\Sph}{\on{Sph}}
\nc{\BunBb}{\overline{\Bun}_B}
\nc{\Buno}{\overset{o}{\Bun}}
\nc{\BunPb}{\overline{\Bun}_P}
\nc{\BunBM}{\overline{\Bun}_{B(M)}}
\nc{\BunPbw}{\widetilde{\Bun}_P}
\nc{\BunBP}{\widetilde{\Bun}_{B,P}}
\nc{\GUb}{\overline{G/U}}
\nc{\GUPb}{\overline{G/U(P)}}
\nc{\GPPb}{\overline{G/[P,P]}}
\nc\Hl{H^{\leftarrow}}
\nc\hl{h^{\leftarrow}}
\nc\Hr{H^{\rightarrow}}
\nc\hr{h^{\rightarrow}}
\nc\h{\mathfrak h}

\nc{\Hhom}{\underline{\on{Hom}}}
\nc\syminfty{\on{Sym}^{\infty}}
\nc\lal{\ol{\lambda}}
\nc\xl{\ol{x}}
\nc\thl{\ol{\theta}}
\nc\nul{\ol{\nu}}
\nc\mul{\ol{\mu}}
\nc\Sum\Sigma

\nc{\M}{{\mathcal M}}
\nc{\N}{{\mathcal N}}
\nc{\F}{{\mathcal F}}
\nc{\D}{{\mathcal D}}
\nc{\Y}{{\mathcal Y}}
\nc{\G}{{\mathcal G}}
\nc{\E}{{\mathcal E}}
\nc{\CalC}{{\mathcal C}}
\renc{\Sb}{\overline{S}}
\nc\Dh{\widehat{\D}}
\renewcommand{\O}{{\mathcal O}}
\nc{\C}{{\mathcal C}}
\nc{\K}{{\mathcal K}}
\renewcommand{\H}{{\mathcal H}}
\renewcommand{\S}{{\mathcal S}}
\nc{\T}{{\mathcal T}}
\nc{\V}{{\mathcal V}}
\renc{\P}{{\mathcal P}}
\nc{\A}{{\mathcal A}}
\nc{\B}{{\mathcal B}}
\nc{\U}{{\mathcal U}}
\renewcommand{\L}{{\mathcal L}}
\nc{\Gr}{\on{Gr}}

\nc{\p}{\overline{\mathfrak p}}
\nc{\q}{\overline{\mathfrak q}}
\nc\f{{\mathfrak f}}
\renewcommand\k{{\mathfrak k}}
\nc{\qo}{{\mathfrak q}}
\nc{\po}{{\mathfrak p}}
\nc{\s}{{\mathfrak s}}
\nc\w{\text{w}}
\renewcommand{\r}{{\mathfrak r}}

\newcommand{\tf}{{\mathfrak t}}

\nc\mathi\iota
\nc\Spec{\on{Spec}}
\nc{\tw}{\widetilde{\mathfrak t}}
\nc{\pw}{\widetilde{\mathfrak p}}
\nc{\qw}{\widetilde{\mathfrak q}}
\nc{\jw}{\widetilde j}

\nc{\grb}{\overline{\Gr}}
\nc{\I}{\mathcal I}
\renewcommand{\i}{\mathfrak i}
\renewcommand{\j}{\mathfrak j}

\nc{\lambdach}{\check\lambda}
\nc{\Lambdach}{\check\Lambda}
\nc{\much}{\check\mu}
\nc{\omegach}{\check\omega}
\nc{\nuch}{\check\nu}
\nc{\etach}{\check\eta}
\nc{\alphach}{\check\alpha}
\nc{\betach}{\check\beta}
\nc{\rhoch}{\check\rho}

\nc{\Hb}{\overline{\H}}

\begin{document}

\title{Geometric Eisenstein series}

\author{A.~Braverman and D.~Gaitsgory}

\address{A.B.:Dept. of Math., MIT, Cambridge MA 02139; D.G:
Dept. of Math., Harvard University, Cambridge MA 02138}

\email{braval@math.mit.edu, gaitsgde@math.harvard.edu}

\maketitle

\tableofcontents

\section*{Introduction} \label{intr}

\ssec{Motivation}

\sssec{}

Let $X$ be a curve over $\Fq$ and let $G$ be a reductive group. The classical 
theory of automorphic forms is concerned with the space of functions on 
the quotient $G_{\AA}/G_\K$, where 
$\K$ (resp., $\AA$) is the field or rational functions on $X$ 
(resp., the ring of ad\`eles of $\K$).
In this paper, we will consider only the unramified situation, i.e. 
we will study functions
(and afterwards perverse sheaves) on the double quotient 
$G_{\OO}\backslash G_{\AA}/G_\K$.

Let $T$ be a Cartan subgroup of $G$. There is a well-known construction, 
called the {\it Eisenstein series operator} 
that attaches to a compactly supported 
function on $T_{\OO}\backslash T_{\AA}/T_\K$ a function on 
$G_{\OO}\backslash G_{\AA}/G_\K$:

Consider the diagram
$$
\CD
B_{\OO}\backslash B_{\AA}/B_\K   @>{\mathfrak q}>> T_{\OO}\backslash T_{\AA}/T_\K \\
@V{\mathfrak p}VV   \\
G_{\OO}\backslash G_{\AA}/G_\K,
\endCD
$$
where $B$ is a Borel subgroup of $G$. Up to a normalization factor, 
the Eisenstein series of a function $S$ on 
$T_{\OO}\backslash T_{\AA}/T_\K$ is ${\mathfrak p}_{!}({\mathfrak q}^*(S))$, where 
${\mathfrak q}^*$ denotes pull-back and ${\mathfrak p}_{!}$ is integration along the fiber.

Our goal is to study a geometric analog of this construction.

\sssec{}

Let $\Bun_G$ denote the stack of $G$--bundles on $X$. One may regard the derived category
of constructible sheaves on $\Bun_G$ (denoted $\Sh(\Bun_G)$, cf. \secref{conventions}) as a 
geometric analog of the space of functions on $G_{\OO}\backslash G_{\AA}/G_\K$. Then, by geometrizing
the Eisenstein series operator, we obtain an Eisenstein series functor
$\Eis':\Sh(\Bun_T)\to\Sh(\Bun_G)$ defined by a diagram
similar to the above one, where the intermediate stack is $\Bun_B$--the stack of $B$--bundles on $X$.

However, this construction has an immediate drawback--it is not 
sufficienly functorial (for example
it does not commute with Verdier duality), 
the reason being that the projection 
${\mathfrak p}:\Bun_B\to\Bun_G$ has non-compact fibers. 
Therefore, it is natural to look for a relative
compactification of $\Bun_B$ along the fibers of the 
projection ${\mathfrak p}$. 

Such a compactification, denoted $\BunBb$, was found by V.~Drinfeld. We will use it to define
the corrected functor $\Eis:\Sh(\Bun_T)\to\Sh(\Bun_G)$. This paper is devoted to the 
investigation of various properties of this functor. 

As we will see, all the technical results will essentially reduce to questions about the geometry 
of $\BunBb$ and the behaviour of the intersection cohomology sheaf on it.

\sssec{}

We should say right away that the pioneering work in this direction was done by G.~Laumon in \cite{La},
who considered the case of $G=GL(n)$ using his own compactification of the stack $\Bun_B$.
In the sequel we will explain how the two approaches are related.

\ssec{Survey of the main results}

\sssec{}

First, let us indicate how Drinfeld's compactifications are constructed. 

Let $G=GL(2)$. In this case,
the stack $\Bun_B$ classifies pairs $(\M,\L,\kappa:\L\hookrightarrow\M)$, where $\M$ is a rank $2$ vector bundle 
on $X$, $\L$ is a line bundle on $X$, and $\kap$ is a non-zero (=injective) bundle map. When $\M$ 
is fixed, the set of all possible $(\L,\kap)$ forms an algebraic variety. This variety is non-complete,
but there is a natural way to compactify it: 

The corresponding complete variety also consists of pairs $(\L,\kap:\L\to\M)$, but this time $\kap$ is
required to be just an injective map of coherent sheaves, and not necessarily a bundle map. In other words,
the quotient $\M/\L$ may have torsion. By definition, the stack $\BunBb$ for $GL(2)$ classifies triples
$(\M,\L,\kap)$, where $\kap$ is understood in the new sense.

For a general $G$ (we are assuming that the derived group of $G$ is simply-connected) the construction of $\BunBb$
is very similar.

Let $\F_G$ be a fixed $G$-bundle. Let us denote by $\V^{\lambdach}_{\F_G}$ the vector bundle on $X$
associated with $\F_G$ and a highest weight $G$-module $\V^{\lambdach}$. A point of $\Bun_B$ that projects
to our $\F_G\in\Bun_G$ defines in every $\V^{\lambdach}_{\F_G}$ as above a line subbundle $\L^{\lambdach}$,
such that the system of these subbundles satisfies the so-called {\it Pl\"ucker relations}, cf. Sect. 1.

By definition, $\BunBb$ classifies the data of $\F_G$ endowed with a system of maps
$\kap^{\lambdach}:\L^{\lambdach}\hookrightarrow \V^{\lambdach}_{\F_G}$ defined for every dominant integral weight
$\lambdach$, where $\kap^{\lambdach}$ are injective as maps of coherent sheaves, which satisfy the Pl\"ucker relations.

Along the same lines one can try to define a compactification of the stack $\Bun_P$, where $P$ is
a parabolic subgroup of $G$. However, in this case there will be two different natural compactifications
$\BunPb$ and $\BunPbw$, of which the second one will be more important for our purposes.

\medskip

Let us make a brief digression and discuss Laumon's compactification of $\Bun_B$ for $G=GL(n)$. In this case,
$\Bun_B$ is the stack classifying the data of a rank $n$ vector bundle $\M$ on $X$ endowed with a complete
flag $0=\M_0\subset \M_1\subset...\subset\M_n=\M$ of subbundles. Laumon's compactification $\BunBb^L$
classifies the data of ($\M$ + a flag of $\M_i$'s), but each $\M_i$ is required to be just a subsheaf
(rather than a subbundle) of the subsequent $\M_{i+1}$.

For $GL(2)$, Laumon's and Drinfeld's compactifications coincide, but this is no longer true for $n\geq 3$.
The advantage of Laumon's compactification is that $\BunBb^L$ is smooth, whereas $\BunBb$ is singular for
$G$ of semi-simple rank $>1$. The drawback is that $\BunBb^L$ makes sense only for $GL(n)$.

In fact, there is a natural (proper and representable) map of stacks $\BunBb^L\to\BunBb$ and it was shown
by A.~Kuznetsov in \cite{Ku} that $\BunBb^L$ is a small resolution of singularities of $\BunBb$.

\sssec{}

Once the stack $\BunBb$ is constructed, one can try to use it to define the "compactified" Eisenstein series
functor $\Eis:\Sh(\Bun_T)\to\Sh(\Bun_G)$. Let $\p$ and $\q$ denote the natural projections from $\BunBb$
to $\Bun_G$ and $\Bun_T$, respectively. The first idea would be to consider the functor
$\S\in\Sh(\Bun_T)\mapsto \p_{!}(\q^*(\S))\in\Bun_G$. However, this is too naive, since if we want 
our functor to commute with Verdier duality, we need to take into account the sigularities of $\BunBb$. 
Therefore, one introduces a {\it kernel} on $\BunBb$ given by its intersection cohomology sheaf. I.e., we define
the functor $\Eis$ by
$$\S\mapsto \p_{!}(\q^*(\S)\otimes\IC_{\BunBb}),$$
up to a cohomological shift and Tate's twist.
Similarly, one defines the functor $\Eis^G_M:\Sh(\Bun_M)\to\Sh(\Bun_G)$, where $M$ is the Levi quotient
of a parabolic $P$. 

The first test whether our definition of the functor $\Eis$ is "the right one" would be the assertion
that $\Eis$ (or more generally $\Eis^G_M$) indeed commutes with Verdier duality. Our $\Eis$ passes this test,
according to \corref{goodEisprince}.

\medskip

Let us again add a comment of how the functor $\Eis$ is conneceted to Laumon's work. One can define
functors $\Eis^L:\Sh(\Bun_T)\to\Sh(\Bun_G)$ using Laumon's compactification. (In the original work \cite{La},
Laumon did not consider $\Eis^L$ as a functor, but rather applied it to specific sheaves on $\Bun_T$.)
However, from the smallness result of \cite{Ku} it follows that the functors $\Eis^L$ and $\Eis$ are canonically
isomorphic.

The situation is more subtle when one considers non-principal Eisenstein series, i.e. when $B$
is replaced by a more general parablic. In this case, it is no longer true that the two definitions agree.
(In particular, for $P\neq B$, Laumon's Eisenstein series fails to commute with Hecke functors, cf. below.) 
At the end of Sect. 2 we suggest a conjecture of how they may be related.

\sssec{}

Once we defined the functors $\Eis=\Eis^G_T:\Sh(\Bun_T)\to\Sh(\Bun_G)$, $\Eis^G_M:\Sh(\Bun_M)\to\Sh(\Bun_G)$
and a similar functor for $M$, $\Eis^M_T:\Sh(\Bun_T)\to\Sh(\Bun_M)$, it is by all means natural to expect 
that these functors compose nicely, i.e. that $\Eis^G_T\simeq \Eis^G_M\circ \Eis^M_T$. 

For example, if instead of $\Eis$ we used the naive (uncompactified)
functor $\Eis'$, the analogous assertion would be a triviality, since
$\Bun_B\simeq \Bun_P\underset{\Bun_M}\times\Bun_{B(M)}$, where $B(M)$ is the Borel subgroup of $M$.

The problem with our definition of $\Eis^G_M$ is that there is no map between the relevant 
compactifications, i.e. from $\BunBb$ to $\BunPbw$. Neverthess, the assertion that 
$\Eis^G_T\simeq \Eis^G_M\circ\Eis^M_T$ does hold. This is a non-trivial theorem (\thmref{compos}), 
which is one of the main results of this paper. 

\sssec{}

Finally, we have to describe the following two results, which are geometric interpretations
of certain properties of the classical (i.e. function-theoretic) Eisenstein series operator.

\medskip

\noindent {\bf Behaviour with respect to the Hecke functors.}

\medskip

Classically, on the space of functions on the double quotient $G_{\OO}\backslash G_{\AA}/G_\K$
we have the action of $\underset{x\in X}\otimes \H_x(G)$, where $x$ runs over the set of places of $\K$,
and for each $x\in X$, $\H_x(G)$ denotes the corresponding spherical Hecke algebra of the group $G$.

Similarly, $\underset{x\in X}\otimes \H_x(T)$ acts on the space of functions on
$T_{\OO}\backslash T_{\AA}/T_\K$. In addition, for every $x$ as above, there is a canonical homomorphism
$\H_x(G)\to \H_x(T)$ described as follows:

Recall that there is a canonical isomorphism (due to Satake) between $\H_x(G)$ and the Grothendieck
ring of the category of finite-dimensional representations of the Langlands dual group $\check G$.
We have the natural restriction functor $\on{Rep}(\check G)\to \on{Rep}(\check T)$, and our homomorphism
$\H_x(G)\to \H_x(T)$ corresponds to the induced homomorphism $K(\on{Rep}(\check G))\to K(\on{Rep}(\check T))$
between Grothendieck rings.

The basic property of the Eisenstein series operators is that it intertwines the $\H_x(G)$-action on
$G_{\OO}\backslash G_{\AA}/G_\K$ and the $\H_x(T)$-action on $T_{\OO}\backslash T_{\AA}/T_\K$
via the above homomorphism.

\smallskip

Our \thmref{commutewithHecke} is a reflection of this phenomenon in the geometric setting. Now, instead
of the Hecke algebras, we have the action of the Hecke functors on $\Sh(\Bun_G)$. Namely,
for $x\in X$ and an object $V\in \on{Rep}(\check G)$, one defines the {\it Hecke functor} 
$$\S\mapsto {}_xH_G(V,\S)$$
from $\Sh(\Bun_G)$ to
itself (cf. \secref{prelaffgrass}) and we have for $\S\in\Sh(\Bun_T)$:
$$_xH_G(V,\Eis(\S))\simeq \Eis({}_xH_T(\on{Res}^G_T(V),\S)).$$

A similar statement holds for the non-principal Eisenstein series functor $\Eis^G_M$, cf. 
\thmref{commutewithHeckenonprince}.

As a coroallary, we obtain that if $E_{\check M}$ is an $\check M$-local system on $X$ and
$\Aut_{E_{\check M}}$ is a perverse sheaf (or a complex of sheaves) on $\Bun_M$, 
corresponding to it in the sense of the geometric Langlands correspondence, then the complex 
$\Eis^G_M(\Aut_{E_{\check M}})$ on $\Bun_G$ is a Hecke eigen-sheaf with respect to the induced
$\check G$-local system.

In particular, \thmref{commutewithHecke} establishes the geometric Langlands correspondence
for those homomorphisms $\pi_1(X)\to \check G$, whose image is contained in a maximal torus of $\check G$.

\medskip

\noindent {\bf The functional equation.}

\medskip

It is well-known that the classical Eisenstein series satisfy the functional equation. Namely, let $\chi$
be a character of the group $T_{\OO}\backslash T_{\AA}/T_\K$ and let $\w\in W$ be an element of the Weyl
group. We can translate $\chi$ by menas of $\w$ and obtain a new (grossen)-character $\chi^{\w}$.

The functional equation is the assertion that the Eisenstein series corresponding to $\chi$ and
$\chi^{\w}$ are equal, up to a ratio of the corresponding L-functions. 

Now let $\S$ be an arbitrary complex of sheaves on $\Bun_T$ and let $\w\cdot\S$ be its $\w$-translate.
One may wonder whether there is any relation between $\Eis(\S)$ and $\Eis(\w\cdot\S)$. 

We single out a subcategory in $\Sh(\Bun_T)$, corresponding to sheaves which we call "regular", for
which we answer the above question. We claim (cf. \thmref{functeq}) that for a regular $\S$ 
$$\Eis(\S)\simeq \Eis(\w\cdot\S).$$

(It is easy to see that one should not expect the functional equation to hold for
non-regular sheaves.)

A remarkable feature of this assertion is that the L-factors that enter the classical functional
equation have disappeared. An explanation of this fact is provided by \thmref{fullcompare},
which says that the definition of $\Eis$ via the intersection cohomology sheaf on $\BunBb$ already
incorporates the $L$--function. 

We remark that an assertion similar to the above functional equation should hold also for
non-principal Eisenstein series. Unfortunately, this is beyond access for our methods.

\medskip

By combining \thmref{commutewithHecke} and \thmref{functeq} we obtain a proof
of a special case of the Langlands conjecture.  Namely, we prove \thmref{langlands}, which says that if
we start with an unramified irreducible representation of $\pi_1(X)$ into $\check G$,
such that $\pi_1(X)^{geom}$ maps to $\check T\subset \check G$, then there 
exists an unramified automorphic form on $G_{\AA}$ which corresponds to this representation in the 
sense of Langlands.

This may be considered as an application of the machinery 
developped in this paper to the classical theory of automorphic forms.

\ssec{Contents}

\sssec{}

Now let us briefly discuss how this paper is organized.

In Section 1 we discuss the definitions of Drinfeld's compactifications.

In Section 2 we formulate all the main theorems concerning the Eisenstein series functor. In the rest
of the paper we deduce these theorems from various geometric properties of the stacks $\BunBb$,
$\BunPbw$ and their close relatives.

\sssec{}

In Sections 3 and 4 our goal is to prove that our Eisenstein series functor 
commutes with the Hecke functors (cf. above). In Section 3 we consider the principal case, 
i.e. $P=B$ and in Section 4 we make a generalization 
to the case of a general parabolic subgroup.
 
Our discussion in these two sections is based on the 
Lusztig-Drinfeld-Ginzburg-Mirkovic-Vilonen theory of spherical
perverse sheaves on the {\it affine Grassmannian} $\Gr_G$, which we 
review in \secref{prelaffgrass}.

\sssec{}

In Section 5 we formulate and prove \thmref{acycthm}, which is one of our
main technical tools. It says that the stack $\BunBb$ (or, more generally,
$\BunPbw$), is locally acyclic with respect to the natural projection 
$\q:\BunBb\to\Bun_T$ (resp., $\qw_P:\BunPbw\to\Bun_M$). In other words, 
this means that from the point of view of singularities,
all the fibers of the above projections look the same.

\thmref{acycthm} has two main applications. On the one hand, it implies 
\thmref{goodtopullbackprince}, which says that the functor $\Eis$
commutes with Verdier duality. On the other hand, it is a key ingredient in the
proof of \thmref{compos} about the composition of Esenstein (cf. above).

\sssec{}

In Sections 6 and 7 our goal is to prove of \thmref{compos} and the functional equation.

Section 6 is largely preparatory; in it we introduce various stratifications
of the stacks $\BunBb$ and $\BunPbw$. 

As we have mentioned before, the main problem in the proof of \thmref{compos} is
that the natural map $\Bun_B\to\Bun_P$ does not extend to a 
map from $\BunBb$ to $\BunPbw$. To cure this, we introduce the stack $\wt{\Bun}_{B,P}$, 
which will map to both $\BunBb$ and $\BunPbw$.

In Section 7 we prove \thmref{restriction}, which asserts that the map $\wt{\Bun}_{B,P}\to\BunPbw$ is 
stratified-small in the sense of Goresky-MacPherson. This will allow us to finally prove 
\thmref{compos}.

The paper is concluded by \secref{comprkone}, in which we complete the 
proof of the functional equation (\thmref{functeq}). As in the 
classical theory, after having established \thmref{compos}, it is sufficient 
to treat the case of reductive groups of semi--simple rank $1$, 
in which case the functional equatiuon isomorphism 
amounts to a Fourier--transform type argument, combined with the Serre duality.

\ssec{Conventions} \label{conventions}

We will work over the base field $\Fq$ and study Weil sheaves over 
$\Fq$--schemes and stacks. However, all the results of this paper 
carry over to the characteristic $0$ world, in which case one must replace
$\Ql$--adic sheaves by holonomic D--modules.

\sssec{The group $G$}

Throughout the paper, $G$ will be a connected reductive group over $\Fq$.
Moreover, we will assume that its derived group $G'=[G,G]$ is simply connected.

We will fix a Borel subgroup $B\subset G$ and let $T$ be the 
``abstract'' Cartan,
i.e., $T=B/U$, where $U$ is the unipotent radical of $B$. We will 
denote by $\Lambda$ the {\it coweight} lattice of
$T$ and by $\check \Lambda$--its dual, i.e. the weight lattice; $\langle\cdot,\cdot\rangle$ is the canonical
pairing between the two.

The semi--group of dominant coweights
(resp., weights) will be denoted by $\Lambda_G^+$ (resp., by $\Lambdach_G^+$). The set of vertices
of the Dynkin diagram of $G$ will be denoted by $\I$; for each $\i\in\I$ there corresponds a simple
coroot $\alpha_i$ and a simple root $\alphach_i$. The set of positive coroots will be denoted by $\Delta$
and their positive span inside $\Lambda$, by $\Lambda_G^{\on{pos}}$. (Note that, since $G$ is not semi--simple,
$\Lambda_G^+$ is not necessarily contained in $\Lambda_G^{\on{pos}}$.) By $\rhoch\in\Lambdach$ we will denote 
the half sum of positive roots of $G$ and by $\w_0$ the longest element in the Weyl group of $G$.
For $\lambda_1,\lambda_2\in\Lambda$ we will write that 
$\lambda_1\geq \lambda_2$ if $\lambda_1-\lambda_2\in\Lambda_G^{\on{pos}}$,
and similarly for $\Lambdach_G^+$.

It will be convenient to choose fundamental weights $\omegach_\i\in\Lambdach^+_G$, $\i\in\I$: each
$\omegach_\i$ is {\it non-uniquely} characterized by the property that 
$\langle \alpha_\j,\omegach_\i\rangle=\delta_{\i,\j}$.

Let $P$ be a standard parabolic of $G$, i.e. $P\supset B$; let $U(P)$ be its unipotent radical
and $M:=P/U(P)$--the Levi factor (note that the derived group of $M$ is also simply connected). To $M$
there corresponds a subdiagram $\I_M$ in $\I$. We will denote by $\Lambda_M^+\supset \Lambdach_G^+$, 
$\Lambda_M^{\on{pos}}\subset \Lambda_G^{\on{pos}}$, $\rhoch_M\in\Lambdach$, $\w_0^M\in W$, $\underset{M}\geq$ etc.
the corresponding objects for $M$.

To a dominant coweight $\lambdach\in\Lambdach$ one attaches the Weyl 
$G$-module, denoted $\V^{\lambdach}$, with a fixed highest weight
vector $v^{\lambdach}\in \V^{\lambdach}$. For a pair of weights 
$\lambdach_1,\lambdach_2$, there is a canonical map
$\V^{\lambdach_1+\lambdach_2}\to \V^{\lambdach_1}\otimes \V^{\lambdach_2}$ that sends 
$v^{\lambdach_1+\lambdach_2}$ to $v^{\lambdach_1}\otimes v^{\lambdach_2}$.

Similarly, for $\nuch\in \Lambdach_M^+$, we will denote by $\U^{\nuch}$ the corresponding $M$--module. 
There is a natural functor $\V\mapsto \V^{U(P)}$ from the category of $G$--modules to that of $M$--modules and 
we have a canonical morphism $\U^{\lambdach}\to (\V^{\lambdach})^{U(P)}$, which sends the corresponding
highest weight vectors to one another.

\smallskip

By $\F_G$ we will denote principal $G$-bundles on various schemes. The notation $\F^0_G$ is reserved
for {\it the trivial} $G$-bundle. For a $G$-bundle $\F_G$ 
(resp., a $T$-bundle $\F_T$) and a $G$-module $\V$ (resp., a character $\lambdach:T\to \GG_m$) let
$\V_{\F_G}$ (resp., $\L_{\F_T}^{\lambdach}$) denote the corresponding associated vector bundle
(resp., line bundle).

If $\F_T$ is a $T$-bundle on a scheme $\Y$, $\mu$ is a coweight of
$T$ and $D\subset\Y$ is a Cartier divisor, we will denote by $\F_T(\mu\cdot D)$ the new $T$-bundle,
defined by the condition that for a weight $\much$ of $T$, the associated line bundle 
$\L^{\much}_{\F_T(\mu\cdot D)}$ equals $\L^{\much}_{\F_T}(\langle \mu,\much\rangle\cdot D)$.

\medskip

Finally, we let $\check G$ denote the Langlands dual group of $G$. This is a reductive group over $\Ql$,
endowed with a canonical Cartan subgroup $\check T\subset \check G$,
whose lattice of characters is identified with $\Lambda$. For $\lambda\in\Lambda^+_G$, we will
denote by $V^\lambda$ the corresponding $\check G$-module and for $\mu\in\Lambda$, $V^\lambda(\mu)$
will designate the $\mu$-th weight subspace of $V^\lambda$. For an $M$--dominant coweight $\nu$, we will
denote the corresponding $\check M$--module by $U^\nu$.

\sssec{Stacks}

The main objects of study in this paper are various stacks attached to our curve $X$ (assumed to be smooth
and complete) and the group $G$.

For instance, the stack of $G$--bundles on $X$, denoted $\Bun_G$, attaches to a scheme $S$ the groupoid
of principal $G$--bundles on $S\times X$ and for a map of schemes $S_1\to S_2$ the corresponding functor
is given by the ordinary pull-back. 

In the main body of the paper, when giving a definition of a stack $\Y$, we will usually say
that $\Y$ classifies {\it something}, by which we mean that it is clear what an $S$--family of this
{\it something} is for a scheme $S$ and what the corresponding pull--back functor on the category
of such families is, for a map $S_1\to S_2$. In addition, in all the examples that the reader will encounter 
in this paper, the fact that our $\Y$ {\it is indeed} an algebraic stack, follows easily from the fact that
$\Bun_G$ is. 

Here is an example:

\smallskip

\noindent The stack $\Y$ classifying the data of $(\M,s)$, where $\M$ is a coherent sheaf on $X$ and $s$ 
is its section (resp., an injective section) means that $\Hom(S,\Y)$ is the groupoid of pairs $(\M_S,s_S)$,
where $\M_S$ is a coherent sheaf on $S\times X$, flat over $S$, and $s_S$ is a map $\O_{S\times X}\to \M_S$
(resp., a map, which is injective over any geometric point of $S$).

A special care must be used when we define stacks or schemes using the formal (resp., formal punctured) disc.
We refer the reader to the Introduction of \cite{FGV} for the corresponding conventions.

\sssec{Derived categories}

If $\Y$ is a scheme, we will denote by $\Sh(\Y)$ the corresponding bounded derived category
of constructible $\Ql$--adic Weil sheaves on $\Y$ (cf. \cite{De}). 
By abuse of language, we will call its objects ``sheaves'', but we will
always mean complexes of such. The constant sheaf on $\Y$ will be denoted 
by ${\Ql}_\Y$. We will denote by $\Ql(\frac{1}{2})$ the $1$-dimensional sheaf over 
$\on{Spec}(\Fq)$ which corresponds to a chosen once and for all square 
root of $q$. 

For a scheme or a stack $\Y$, $\on{Perv}(\Y)$ will denote the abelian category of perverse sheaves on $\Y$.
$\IC_{\Y}\in \on{Perv}(\Y)$ will denote the {\it intersection cohomology sheaf} on $\Y$, 
normalized so that it is pure of weight $0$. E.g., when $\Y$ is smooth, $\IC_{\Y}=
(\Ql(\frac{1}{2})[1])_\Y^{\otimes \on{dim}(\Y)}$. Unless specified otherwise, we work with the
perverse t-structure, i.e. all the cohomological degrees are counted in the perverse normalization.

When $\Y$ is a stack, the definition of the derived category is problematic. However, this poses
no difficulty in our case, since all our stacks are unions of open substacks of finite type, each being
a quotient of a scheme by an affine group. If $\Y=\underset{i}\cup \Y_i$ is such a stack, then
$\Sh(\Y):=\underset{\leftarrow}{\lim} \Sh(\Y_i)$, where each $\Sh(\Y_i)$ can be defined as an
equivariant derived category in the sense of Bernstein--Lunts \cite{BeLu}.

If $f:\Y_1\to \Y_2$ is a map of schemes (resp., of stacks) we will denote by $f_{!}, f_{*}, f^*, f^{!}$
the corresponding functors between $\Sh(\Y_1)$ and $\Sh(\Y_2)$ always understood in the derived sense.
Note, that if $f$ is non-representable, the functors $f_{!}$ and $f_{*}$ do not in general preserve our
bounded derived category. Luckily, in the main body of the paper, we will have to consider the
above direct image functors only for representable morphisms (an exception being the definition of regular sheaves
and the discussion in \secref{comprkone}, in which case their usage will be but very mild).

Another warning on our abuse of language: let $f,f':\Y_1\to \Y_2$ be two morphisms of stacks. We will sometimes
say that they coincide, by which we mean, of course, that they are isomorphic. At any event, in this situation,
the pull-back functors $f^*$ and $f'{}^*$ from $\Sh(\Y_2)\to\Sh(\Y_1)$ are isomorphic too.

If $i:Z\hookrightarrow \Y$ is a locally closed embedding and $\S\in\Sh(\Y)$, the notation $\S|_{Z}$
will mean the same as $i^*(\S)$. The Verdier duality functor will be denoted by $\DD$. For $\S$ as above,
$h^i(\S)\in \on{Perv}(\Y)$ is the $i$-th perverse cohomology of $\S$.

\medskip

Let $\F_H\to \Y$ be an $H$--torsor, where $H$ is an algebraic group and let $Z$ be an $H$--scheme.
We will denote by $Z\overset{H}\times \Y$ the corresponding fibration over $\Y$ with the typical fiber $Z$.
(The word ``fibration'' in this paper will be understood only in this sense.) Let $\S$ (resp., $\T$) be an object
of $\Sh(Z)$ (resp., $\Sh(\Y)$) and assume, moreover, that $\S$ is $H$--equivariant. In this case we can 
form their twisted external product, denoted $\S\tboxtimes\T$, which will be an object of 
$\Sh(Z\overset{H}\times \Y)$.
 
\medskip

\noindent{\bf Acknowledgements.}

\medskip

This paper owes its existence to V.~Drinfeld, who has explained to (one of) 
us the definition of $\BunBb$ and recommended to study it. 
Discussions with him and also with A.~Beilinson
and J.~Bernstein helped us overcome many technical and ideological 
difficulties on the way.

Besides the work of Laumon, the present paper has been much influenced by the Beilinson--Drinfled 
work on the geometric Langlands correspondense and by the papers of
Finkelberg and Mirkovic on semi-infinite flags.

\section{Definition of Drinfeld's compactifications}

\ssec{Strongly quasi-affine varieties}

\sssec{} Let $Y$ be a scheme over an algebraically closed 
field $k$. We shall say that $Y$ is
{\it strongly quasi-affine} if the following conditions hold:

1) The algebra $A_Y=\Gamma(Y,\O_Y)$ of global functions on $Y$ is finitely
generated as a $k$-algebra.

2) Let $\oY=\text{Spec}(A_Y)$. Then the natural map 
$Y\to\oY$ is an open embedding.

\smallskip

If $Y$ is strongly quasi-affine then we will call $\oY=\Spec(\A)$
the {\it affine closure} of $Y$. In this section we will need the following result.

\begin{thm} \label{quasi-affine}
Let $G$ be a reductive group over $k$ and let $P$ be a parabolic
subgroup of $G$. Let $H$ be either the derived group
$[P,P]$ of $P$ or the unipotent radical $U(P)$ of $P$. 
 Let $Y=G/H$.
Then
\begin{enumerate}
\item $Y$ is a strongly quasi-affine variety.
\item For any scheme $S$ over $k$ an $S$-point of
$\oY$ is described by a collection of maps of $\O_S$-modules
$\kap^\V : \O_S\otimes \V^H\to \O_S\otimes \V$
for every $G$-module $\V$ satisfying the following conditions:

a) If $\V$ is the trivial representation of $G$ then $\kap^\V=\id$.

b) Assume that we are given a morphism $a:\V^1\to \V^2$ of
$G$-modules. Then the following square must be commutative:
\begin{equation*} \label{pluckergeneralone}
\CD
(\V^1)^H\otimes \O_S @>{\kappa^{\V^1}}>> \V^1\otimes \O_S \\
@VVV          @VVV       \\
(\V^2)^H\otimes\O_S @>{\kappa^{\V^2}}>> \V^2\otimes \O_S
\endCD
\end{equation*}

c) For a pair of $G$-modules, $\V^1$ and $\V^2$,
the square
\begin{equation*} \label{pluckergeneraltwo}
\CD
(\V^1)^H \otimes (\V^2)^H \otimes \O_S
@>\kappa^{\V^1}\otimes \kappa^{\V^2}>> 
\V^1\otimes \V^2 \otimes \O_S\\
@VVV          @V{\on{id}}VV       \\
(\V^1\otimes\V^2)^H\otimes \O_S @>{\kappa^{\V^1\otimes\V^2}}>> 
\V^1\otimes\V^2\otimes \O_S
\endCD
\end{equation*}
should commute too.
\end{enumerate}
\end{thm}

The second assertion of \thmref{quasi-affine} is straightforward for any
$H$ for which $G/H$ is strongly quasi-affine. The first assertion is
straighforward when $H=[P,P]$ and when $H=U(P)$ this is a theorem
of F.~D.~Grosshans (cf. \cite{Gr}).

\sssec{} Let $G$, $H$ and $Y$ be as above. It will be convenient to describe $S$-points 
of $Y$ slightly differently. For a $P/H$-module $\U$, let $\on{Ind}(\U)$ denote
the $G$-module induced from the corresponding $P$-module, i.e. for a $G$-module $\V$, we have functorially:
$$\Hom_{P/H}(\U,\V^H)\simeq \Hom_G(\on{Ind}(\U),\V).$$
From \thmref{quasi-affine} we obtain using the Frobenius reciprocity:

\begin{lem}  \label{descrq-a}
For a $k$-scheme $S$, an $S$-point of $\oY$ is the same as a collection of maps
of $\O_S$-modules $\kap^\U : \O_S\otimes \U\to \O_S\otimes \on{Ind}(\U)$
for every $P/H$-module $\U$ satisfying the following conditions:

a) If $\U$ is the trivial representation of $P/H$ then $\on{Ind}(\U)$ is canonically trivial as well and 
$\kap^\U$ must equal the identity map.

b) Assume that we are given a morphism $a:\U^1\to \U^2$ of
$G$-modules. Then the following square must be commutative:

\begin{equation*} 
\CD
\U^1\otimes \O_S @>{\kappa^{\U^1}}>> \on{Ind}(\U^1)\otimes \O_S \\
@VVV          @VVV       \\
\U^2\otimes\O_S @>{\kappa^{\U^2}}>> \on{Ind}(\U^2)\otimes \O_S
\endCD
\end{equation*}

c) For a pair of $P/H$-modules, $\U^1$ and $\U^2$,
the square
\begin{equation*} 
\CD
(\U^1 \otimes \U^2) \otimes \O_S
@>\kappa^{\U^1\otimes\U^2}>> 
\on{Ind}(\U^1\otimes \U^2) \otimes \O_S\\
@V\on{id}VV          @VVV       \\
(\U^1\otimes\U^2) \otimes \O_S @>{\kappa^{\U^1}\otimes\kap^{\U^2}}>> 
\on{Ind}(\U^1)\otimes\on{Ind}(\U^2) \otimes \O_S
\endCD
\end{equation*}
should commute too.

\end{lem}

\ssec{Drinfeld structures}
Consider the stacks $\Bun_G$, $\Bun_B$ and $\Bun_T$. 
The maps of groups $B\hookrightarrow G$ and $B\twoheadrightarrow T$
give rise to a diagram of stacks:

$$
\CD
\Bun_B @>{\qo}>> \Bun_T \\
@V{\po}VV  \\
\Bun_G.
\endCD
$$

All the three stacks are smooth. It is easy to see, that the projection 
$\qo$ is smooth, but in general 
non-representable and that the projection $\po$ is representable, but neither smooth nor proper.

\smallskip

Our first goal is to define (following Drinfeld) a compactification of $\Bun_B$ along the 
fibers of the projection $\po$; in other words
we will construct a stack $\BunBb$ which contains $\Bun_B$ as an open sub-stack and which 
is endowed with projections 
$$\p:\BunBb\to\Bun_G \text{ and } \q:\BunBb\to\Bun_T,$$ 
such that the map $\p$ is representable is proper.

\sssec{}  \label{intrcompact}

Note that the stack $\Bun_B$, by definition, classifies the following data:
$$(\F_G;\F_T;\kappa:\F_G\to G/U{\overset{T}\times}\F_T),$$
where $\F_G$ is a $G$-bundle, $\F_T$ is a $T$-bundle and $\kappa$ is a $G$-equiavariant map. 

\smallskip

Recall that $G/U$ is a strongly quasi-affine variety and let $\GUb$ denote its affine closure. 
The groups $G$ and $T$ act on $G/U$ and therefore also on $\GUb$.

\smallskip

We define the functor $\BunBb$ on the category of $k$-schemes as follows:
an $S$-point of $\BunBb$ is a triple $(\F_G,\F_T,\kappa)$, where
$\F_G$ (resp., $\F_T$) is an $S$-point of $\Bun_G$ (resp., of $\Bun_T$) and
$\kappa$ is a $G$-equivariant map
$$\F_G\to \GUb{\overset{T}\times}\F_T,$$
such that for every geometric point $s\in S$ there is a Zariski-open subset $X^0\subset X\times s$ such that the map 
$$\kap|_{X^0}:\F_G|_{X^0} \to \GUb{\overset{T}\times}\F_T|_{X^0}$$ factors through
$G/U{\overset{T}\times}\F_T|_{X^0}\subset \GUb{\overset{T}\times}\F_T|_{X^0}$.

In the sequel, we will simply say that $\BunBb$ classifies the data of
$(\F_G;\F_T;\kappa:\F_G\to \GUb{\overset{T}\times}\F_T)$.

\medskip

In order to show that $\BunBb$ is an algebraic stack in the smooth topology with the
desired properties, we will slightly rewrite the above definition. 
Let us apply \lemref{descrq-a} in our case when $H=U$. Since $P/H=B/U=T$ is commutative, it
is enough to consider only $1$-dimensional $P/H$-modules $\U$. If $\lambdach$ is a character of
$T$, the corresponding induced $G$-module vanishes unless $\lambdach\in \Lambdach_G^+$ and in the latter
case it equals $\V^{\lambdach}$.

From \lemref{descrq-a} we obtain that
an $S$-point of $\BunBb$ is a triple $(\F_G,\F_T,\kappa^{\lambdach},\,\,\forall \lambdach\in \Lambdach_G^+)$, where
$\F_G$ and $\F_T$ are as above, and $\kappa^{\lambdach}$ is a map of coherent sheaves 
$$\L_{\F_T}^{\lambdach}\hookrightarrow \V^{\lambdach}_{\F_G},$$
(cf. our conventions in \secref{conventions})
such that for evey geometric point $s\in S$ the restriction $\kappa^{\lambdach}|_{X\times s}$ is an
injection. The last condition is equivalent to saying that $\kappa^{\lambdach}$ is an injection such that the
quotient $\V^{\lambdach}_{\F_G}/\on{Im}(\kappa^{\lambdach})$ is $S$-flat.

The system of embeddings $\kappa^{\lambdach}$ must satisfy conditions (a)-(c) of \lemref{descrq-a}. Let us
write them down explicitly in our case. We obtain the so-called Pl\"ucker relations:

\smallskip

First, for $\lambdach=0$, $\kappa^0$ must be the identity map 
$\O\simeq \L^0_{\F_T}\to \V^0_{\F_G}\simeq \O$. Secondly, 
for two dominant integral weights $\lambdach$ and $\much$, the map
\begin{equation*}  \label{pluckerone}
\L^{\lambdach}_{\F_T}\otimes\L^{\much}_{\F_T}
\overset{\kappa^{\lambdach}\otimes\kappa^{\much}}\longrightarrow
\V^{\lambdach}_{\F_G}\otimes \V^{\much}_{\F_G}\simeq  
(\V^{\lambdach}\otimes \V^{\much})_{\F_G}
\end{equation*}
must coincide with the composition
\begin{equation*} \label{pluckertwo}
\L^{\lambdach}_{\F_T}\otimes\L^{\much}_{\F_T}\simeq \L^{\lambdach+\much}_{\F_T}
\overset{\kappa^{\lambdach+\much}}\longrightarrow
\V^{\lambdach+\much}_{\F_G}\to (\V^{\lambdach}\otimes \V^{\much})_{\F_G}.
\end{equation*}

\begin{prop} \label{reprprince}
The functor $\BunBb$ is an algebraic stack in the smooth topology. The natural map
$\BunBb\to\Bun_G$ given by $(\F_G,\F_T,\kappa)\mapsto\F_G$ is representable and proper.
\end{prop}

\begin{proof}

Let $\F_G$ be an $S$-point of $\Bun_G$. We have to prove that the Cartesian product
$S\underset{\Bun_G}\times\BunBb$ is representable by a scheme proper over $S$. 

For a weight $\lambdach$ consider the Hilbert scheme $\on{Hilb}^{\lambdach}$
whose $S'$ points are pairs:

\smallskip

\hskip1cm (A map $f:S'\to S$; a subsheaf $\L$ of $f^*(\V^\lambda_{\F_G})$ on $X\times S'$
such that $\L$ is locally free of rank $1$ and the quotient $f^*(\V^\lambda_{\F_G})/\L$ is $S'$-flat.)

\medskip

It is well-known that $\on{Hilb}^{\lambdach}$ exists and is proper over $S$. Now, let 
$(\F_G,\F_T,\kappa)$ be an $S'$-point of $S\underset{\Bun_G}\times\BunBb$. To it we can
attach an $S'$-point of $\on{Hilb}^{\lambdach}$ for every $\lambdach$, by setting $\L=\L^{\lambdach}_{\F_T}$.

Moreover, since every collection of $\{\kappa^{\lambdach}\}$'s satisfying the
Pl\"{u}cker relations is uniquely determined by its values on the fundamental weights, i.e. by the
$\kappa^{\check\omega_{\i}}$'s, $\i\in \I$, we obtain that the above morphisms of functors define
a closed embedding
$$S\underset{\Bun_G}\times\BunBb\hookrightarrow \underset{\i\in\I}\times \on{Hilb}^{\omegach_{\i}}.$$

Therefore, our Cartesian product is also representable by a scheme proper over $S$.

\end{proof}

Obviously, we have an open embedding of $\Bun_B$ into $\BunBb$. The image corresponds to those triples
$(\F_G,\F_T,\kappa)$, for which $\kappa$ factors through $G/U{\overset{T}\times}\F_T\subset
\GUb{\overset{T}\times}\F_T$. In terms of the embeddings $\kappa^{\lambdach}$, an $S$-point 
belongs to $\Bun_B$ if and only if each
$$\kappa^{\lambdach}:\L^{\lambdach}_{\F_T}\to \V^{\lambdach}_{\F_G}$$
is a bundle map, i.e. the quotient is not only $S$-flat, but is in fact $X\times S$-flat.
We will denote this open embedding by $j_B$, or simply by $j$.

We will denote the natural map $\BunBb\to\Bun_G$ by $\p$, since it extends the map $\po:\Bun_B\to\Bun_G$.
The map $(\F_G,\F_T,\kappa)\mapsto \F_T:\BunBb\to\Bun_T$ will be denoted by $\q$ for the same reason.

\begin{prop} \label{denseprince}
$\Bun_B$ is dense in $\BunBb$.
\end{prop}

The proof will be given in \secref{stra1}. We remark that the simply-connectedness 
assumption on $G'$ was necessary precisely in order to make \propref{denseprince} true. 

\sssec{} We will call $\Fqb$--valued points of $\BunBb$ ``Drinfeld's structures''. Let us
see explicitly what they look like.

\medskip

\noindent{\bf Note on the terminology:}
If $\M_1$ and $\M_2$ are vector bundles on a scheme $\Y$, $\kap:\M_1\to\M_2$ is an embedding
of coherent sheaves and $Z\subset \Y$ is a closed subscheme,
we will say that $\kap$ {\it has no zeroes on $Z$} or {\it is maximal at $Z$} if the following
equiavlent conditions hold:

\begin{enumerate}

\item The quotient $\M_2/\M_1$ has no torsion supported on $Z$.

\item There is no non-trivial $\N$ with $\M_1\subset \N\subset \M_2$, such that
$\N/\M_1$ is a torsion sheaf supported on $Z$.

\end{enumerate}

\medskip

\begin{prop} \label{pointsofBunBb}
Let $(\F_G,\F_T,\kap)$ be an $\Fqb$-point of $\BunBb$. Then 

\smallskip

\noindent{\em(1)} There is a Zariski-open subset $X^0\subset X$, such $(\F_G,\F_T,\kap)$
defines a reduction of $\F_G|_{X^0}$ to $B$.

\smallskip

\noindent{\em(2)} There exists a canonical principal $B$-bundle $\F'_B$ on $X$, such that
$\F'_B|_{X^0}$ is the $B$-bundle of point (1).
\end{prop}

\begin{proof}

Consider the embeddings $\kap^{\omegach_\i}:\L_{\F_T}^{\omegach_\i}\to \V_{\F_G}^{\omegach_\i}$, $\i\in\I$.
We take $X^0$ to be the maximal open subset of $X$, such that none of the above $\kap^{\omegach_\i}$'s has
zeroes on $x\in X^0$. 

From the Pl\"ucker relations, we obtain that on $X^0$, the other $\kap^{\lambdach}$'s have no zeroes either.
Hence, $\F_B:=(\F_G|_{X^0},\F_T|_{X^0},\kap|_{X^0})$ is the thought-for $B$-bundle on $X^0$.

\medskip

Let $X-X^0=\{x_1,...,x_n\}$. To each of the points $x_k$, we can assign an element $\nu_k\in\Lambda_G^{\on{pos}}$.
Namely, $\nu_k$ is such that $\langle \nu_k,\omegach_\i\rangle$ is the order of zero of $\kap^{\omegach_\i}$
at $x_k$.

If we put $\F'_T=\F_T(\underset{k}\Sigma\, \nu_k\cdot x_k)$, we obtain a system of maximal embeddings
$$\kap'{}^{\omegach_\i}:\L_{\F'_T}^{\omegach_\i}\to \V_{\F_G}^{\omegach_\i},$$
satisfying the Pl\"ucker relations. 

Hence $\F'_B:=(\F_G,\F'_T,\kap')$ is a $B$-bundle on $X$, which on $X^0$ coincides with $\F_B$.

\end{proof}

Let $(\F_G,\F_T,\kap)$, $\{x_1,...,x_n\}$ and $\{\nu_1,...,\nu_k\}$ be as above. We will say that our
Drinfeld's structure has {\it singularities} on $X-X^0$ and we will call $\nu_k$ its {\it defect} at $x_k$.
The $B$-bundle $\F_B$ will be called the saturation of $(\F_G,\F_T,\kap)$.

\ssec{Parabolic Drinfeld's structures}

Let us now generalize the above definitions to the case, when $B\subset G$ is 
replaced by a parabolic subgroup $P$. In this
case, there will be two different compactifications of the stack $\Bun_P$ 
along the fibers of the projection $\Bun_P\to\Bun_G$.

\sssec{}  \label{BunPb}

For a parabolic $P$, let $U(P)$, $M$ and $\I_M$ be as in \secref{conventions}. We will
denote by $\Lambda_{G,P}$ the quotient lattice
$$\Lambda/\on{Span}\{\alpha_\i,\,\i\in\I_M\}.$$
Let $\Lambdach_{G,P}$ denote the dual lattice. In other words, 
$\Lambdach_{G,P}$ is the lattice of characters of the torus $M/[M,M]=P/[P,P]$. We have:
\begin{equation*} \Lambdach_{G,P}=\{\lambdach\in\Lambdach\,|\, 
\langle \alpha_{\i},\lambdach\rangle=0 \text{ if } \i\in \I_M\}.
\end{equation*}

Let $\Lambda^{\on{pos}}_{G,P}$ denote the
positive part of $\Lambda_{G,P}$, i.e. the span of the projections of $\alpha_{\i}$, $\i\in \I-\I_M$.

\sssec{}

The stack $\Bun_P$ classifies the data of 
$$(\F_G,\F_{M/[M,M]},\kappa_P:\F_G\to G/[P,P]{\overset{M/[M,M]}\times}\F_{M/[M,M]}),$$
where $\F_G$ (resp., $\F_{M/[M,M]}$) is a $G$-bundle (resp., an $M/[M,M]$-bundle) and 
$\kappa$ is a $G$-equivariant map.

\smallskip

We let $\qo^\dagger_P$ (resp., $\po_P$) denote the natural projection 
$\Bun_P\to\Bun_{M/[M,M]}$ (resp., $\Bun_P\to\Bun_G$). As in the case of 
the Borel subgroup, the projection $\qo^\dagger_P$ is smooth (but non-representable) 
and the projection $\po_P$ is representable, but neither smooth nor proper.

\smallskip

The variety $G/[P,P]$ is strongly quasi-affine and let $\GPPb$ denote its affine closure. We let 
$\BunPb$ be the stack that classifies triples
$$(\F_G,\F_{M/[M,M]},\kappa_P:\F_G\to \GPPb{\overset{M/[M,M]}\times}\F_{M/[M,M]}),$$
where $\F_G$ and $\F_{M/[M,M]}$ are as above and $\kappa_P$ is a $G$-equivariant map,
which over the generic point of $X$ maps $\F_G$ to 
$$G/[P,P]{\overset{M/[M,M]}\times}\F_{M/[M,M]}\subset\GPPb{\overset{M/[M,M]}\times}\F_{M/[M,M]},$$
in the same sense as in the definition of $\BunBb$.

\medskip

Here is a Pl\"ucker-type description of $\BunPb$: 

An $S$-point of $\BunPb$ consists of 
$(\F_G,\F_{M/[M,M]},\kappa_P^{\lambdach},\,\,\forall\lambdach\in \Lambdach_{G,P}\cap \Lambda_G^+)$,
where $\F_G$ and $\F_{M/[M,M]}$ are as above, and $\kappa_P^{\lambdach}$ are maps of coherent sheaves
$$\L^{\lambdach}_{\F_{M/[M,M]}}\hookrightarrow \V^{\lambdach}_{F_G},$$
which satisfy the same conditions as in the case of $\BunBb$.

\smallskip

Similarly to the previous case, we obtain that $\BunPb$ is indeed an algebraic stack, representable 
and proper over $\Bun_G$. We will denote by $\p_P$ the natural map $\BunPb\to\Bun_G$ and by
$\q^\dagger_P$ the map $\BunPb\to\Bun_{M/[M,M]}$ given by $(\F_G,\F_{M/[M,M]},\kappa_P)\mapsto \F_{M/[M,M]}$.
By $j_P$ we will denote the open embedding $\Bun_P\hookrightarrow \BunPb$. 

\sssec{}  \label{pointsofBunPb}

We will call $\Fqb$-points of $\BunPb$ ``parabolic Drinfeld's structures''. An analog of \propref{pointsofBunBb}
holds in the parabolic case as well. Namely, to an $\Fqb$-point $(\F_G,\F_{M/[M,M]},\kappa_P)$ of $\BunPb$
one can attach a finite collection of points $x_1,...,x_n\in X$ and elements 
$\nu_1,...,\nu_n\in\Lambda^{\on{pos}}_{G,P}$, such that:

On $X^0:=X-\{x_1,...x_n\}$, $(\F_G,\F_{M/[M,M]},\kappa_P)$ defines a $P$-bundle. This $P$-bundle can be 
extended to a $P$-bundle $\F'_P=(\F_G,\F'_{M/[M,M]},\kappa'_P)$ on the entire $X$ with
$\F'_{M/[M,M]}=\F_{M/[M,M]}(\underset{k}\Sigma\, \theta_k\cdot x_k)$. We will call this $\F'_P$ the saturation
of $(\F_G,\F_{M/[M,M]},\kappa_P)$.

\sssec{}

Observe now that the inclusion of groups $B\subset P$ gives rise to a map of stacks $\Bun_B\to\Bun_P$.
One may wonder whether this map extends to a map between $\BunBb$ and $\BunPb$. In this case the answer
is positive:

Using the natural maps $T\to M/[M,M]$ and $\GUb\to \GPPb$, to every $S$-point $(\F_G,\F_T,\kappa)$ of
$\BunBb$ we can assign in a functorial way an $S$-point $(\F_G,\F_{M/[M,M]},\kappa_P)$ of $\BunPb$. We will
denote this map of stacks by $\s_P$. In terms of the Pl\"ucker picture, $\s_P$ "remembers"
only the subsheaves $\L^{\lambdach}_{\F_T}\hookrightarrow \V^{\lambdach}_{\F_G}$ for
$\lambdach\in\Lambdach_{G,P}\cap\Lambdach_G^{+}$.

We have $\p=\p_P\circ \s_p$. This implies, in particular, that $\s_P$ is representable and proper.

\medskip

In addition to the map $\qo^\dagger_P:\Bun_P\to \Bun_{M/[M,M]}$, there is a natural map 
$\qo_P:\Bun_P\to\Bun_M$; one can
show that this map does {\it not} extend to a map of stacks $$\BunPb\to\Bun_M.$$ 
We will construct now another relative compactification $\BunPbw$ of $\Bun_P$, such that 
$\qo_P$ will extend to a map $\qw_P:\BunPbw\to\Bun_M$.

The stack $\BunPbw$ will project naturally to $\BunPb$; however, the
map $\s_P:\BunBb\to\BunPb$ will {\it not} lift to a map $\BunBb\to\BunPbw$.

\sssec{}  \label{BunPbw}

The stack $\Bun_P$ can be also viewed as a stack that classifies triples
$$(\F_G,\F_M,\widetilde\kappa_P:\F_G\to G/U(P){\overset{M}\times}\F_M),$$
where $\F_G$ and $\F_M$ are $G$-- and $M$--bundles respectively and where 
$\widetilde\kappa_P$ is a $G$-equivariant map. 
The variety $G/U(P)$ is strongly quasi-affine, according to \thmref{quasi-affine},
and let $\GUPb$ denote its affine closure. 

We set $\BunPbw$ to be the stack that classifies triples 
$$(\F_G,\F_M,\widetilde\kappa_P:\F_G\to \GUPb{\overset{M}\times}\F_M),$$
where $\widetilde\kappa_P$ is a $G$-equivariant map, such that over 
the generic point of $X$, $\F_G$ is mapped to
$$G/U(P){\overset{M}\times}\F_M\subset \GUPb{\overset{M}\times}\F_M,$$
as in the definition of $\BunBb$.

To prove that $\BunPbw$ is indeed an algebraic stack, representable and proper over $\Bun_G$ we
will once again have to resort to the Pl\"ucker picture:

By \thmref{quasi-affine}, as $S$-point of $\BunPbw$ is a triple $(\F_G,\F_M,\wt\kappa_P^\V)$, where
$\wt\kappa_P^\V$ is a map of coherent sheaves, defined for every $G$-module $\V$:
$$(\V^{U(P)})_{\F_M}\to \V_{\F_G},$$
such that 

\noindent 1) For every geometric point $s\in S$, the restriction $\wt\kappa_P^\V|_{X\times s}$ is an injection.

\noindent 2) Pl\"ucker relations hold in the following sense:

a) For the trivial representation $\V$, 
$\widetilde\kappa_P^{\V}$ must be the identity map $\O\to\O$.

b) For a $G$-module map $\V^1\to\V^2$ the square
\begin{equation*} \label{pluckernonprinceone}
\CD
((\V^1)^{U(P)})_{\F_M} @>{\widetilde\kappa_P^{\V^1}}>> \V^1_{\F_G} \\
@VVV          @VVV       \\
((\V^2)^{U(P)})_{\F_M} @>{\widetilde\kappa_P^{\V^2}}>> \V^2_{\F_G}
\endCD
\end{equation*}
should commute. 

c) For a pair of $G$-modules, $\V^1$ and $\V^2$, the square
\begin{equation*} \label{pluckernonprincetwo}
\CD
((\V^1)^{U(P)} \otimes (\V^2)^{U(P)})_{\F_M} 
@>{\widetilde\kappa_P^{\V^1}\otimes \widetilde\kappa_P^{\V^2}}>> 
\V^1_{\F_G}\otimes \V^2_{\F_G} \\
@VVV          @V{\on{id}}VV       \\
((\V^1\otimes\V^2)^{U(P)})_{\F_M} @>{\widetilde\kappa_P^{\V^1\otimes\V^2}}>> 
\V^1_{\F_G}\otimes\V^2_{\F_G}
\endCD
\end{equation*}
should commute too.

\begin{prop}  \label{repnonprince}
The morphism of functors $\BunPbw\to\Bun_G$ 
is representable and proper.
\end{prop}

\begin{proof}

The proof will use the same idea as the proof of \propref{reprprince}, of which
the latter is a particular case, since for $P=B$, $\BunPb=\BunPbw=\BunBb$.

Let $\V$ be a rerpesentation of $G$; put $n=\dim(\V)$, $k=\dim(\V^{U(P)})$. Let $\Bun_n=\Bun_{GL(n)}$
denote the stack of rank $n$ vector bundles on $X$ and let $\Bun_{k,n}$ denote the stack classifying pairs
$\F_k\hookrightarrow \F_n$ where $\F_n$ is a rank $n$ vector bundle on $X$ and $\F_k$ is a rank $k$ subsheaf of
$\F_n$.

We have a natural map of stacks $\Bun_G\to\Bun_n$, which sends $\F_G$ to the vector bundle 
$\V_{\F_G}$. In addition, we have a map of functors $\BunPbw\to \Bun_{k,n}\underset{\Bun_n}\times \Bun_G$,
that sends a triple $(\F_G,\F_M,\wt\kap_P)$ to $(\V^{U(P)}_{\F_M}\hookrightarrow \V_{\F_G}, \F_G)$.

It is easy to see that when $\V$ is large enough so that the image of $\V^*\otimes \V^{U(P)}$ in
$A_{G/U(P)}:=\Gamma(G,\O_G)^{U(P)}$ generates it as an algebra, the above map
$$\BunPbw\to \Bun_{k,n}\underset{\Bun_n}\times \Bun_G$$
is a closed embedding of functors.

This proves the proposition, since the map $\Bun_{k,n}\to \Bun_n$ is known to be
representable and proper (cf. \cite{La}). 

\end{proof}

We will call $\Fqb$-valued points of $\BunPbw$ ``enhanced parabolic Drinfeld's structures''. A point
$(\F_G,\F_M,\wt\kap_P)\in\BunPbw(\Fqb)$ has
singularities at the same finite collection of points as the corresponding simple parabolic Drinfeld's structure,
i.e. $\r_p(\F_G,\F_M,\wt\kap_P)\in\BunPb(\Fqb)$. However, the structure of the defect in the enhanced case
is more complicated and we will study it in detail in
\secref{stra2}.
 
\begin{prop}  \label{densenonprince}
The stack $\Bun_P$ is dense in both $\BunPb$ and $\BunPbw$.
\end{prop}

The proof is postponed until \secref{stra2}.

\sssec{}

We will denote by $\jw_P$ the open embedding $\Bun_P\hookrightarrow \BunPbw$ and by $\pw_P$ (resp., $\qw_P$)
the natural map $\BunPbw\to \Bun_G$ (resp., $\BunPbw\to\Bun_M$) given by $(\F_G,\F_M,\wt\kap_P)\mapsto \F_G$
and $(\F_G,\F_M,\wt\kap_P)\mapsto \F_M$, respectively.

\smallskip

Let us observe now that there is a natural map of stacks $\r_P:\BunPbw\to\BunPb$, which extends the identity map
on $\Bun_P$: 

In terms of $(\F_G,\F_M,\wt\kap_P)$, it corresponds to the projection $M\to M/[M,M]$ and the map
$\GUPb\to\GPPb$. In terms of the Pl\"ucker picture, $\r_P$
``remembrs'' only the embeddings $\widetilde\kappa_P^{\V}$
for $\V=\V^{\lambdach}_{\F_G}$, 
$\lambdach\in\Lambdach_{G,P}\cap\Lambda_G^{+}$ (in this case the 
representation 
$\U^{\lambdach}$ of $M$ factors 
through $M\to M/[M,M]$). 

Since $\pw_P=\p_P\circ \r_P$, the map $\r_P$ is representable and proper. We will show in \secref{stra2} 
that its fibers can be identified with 
closed subvarieties of the affine Grassmannian of the group $M$.


\section{Main results on Eisenstein series}

\ssec{Principal Eisenstein series} 

First, we will discuss the principal Eisenstein series functor that maps 
the category $\Sh(\Bun_T)$ to $\Sh(\Bun_G)$.

\sssec{}

Let $\q^{!*}$ denote the functor $\Sh(\Bun_T)\to \Sh(\BunBb)$ given by
\begin{equation*}
\S\to \IC_{\BunBb}\otimes\q^*(\S)\otimes(\Ql(\frac{1}{2})[1])^{\otimes-\dim(\Bun_T)}.
\end{equation*}
Analogously, we define the functor $\qo^{!*}:=j^*\circ\q^{!*}:\, \Sh(\Bun_T)\to \Sh(\Bun_B)$. 
Note, that since the map $\qo$ is smooth, the latter functor is an ordinary pull-back up to a 
cohomological shift and Tate's twist.

\begin{thm} \label{goodtopullbackprince} 

\smallskip

(a) The functor $\q^{!*}$ is exact and commutes with Verdier duality.

\smallskip

(b) If $\S$ is perverse, we have:
$\q^{!*}(\S)\simeq j_{!*}\circ \qo^{!*}(\S)$.

\end{thm}

We define the functor $\Eis^G_T$ (or simply $\Eis$) $\Sh(\Bun_T)\to 
\Sh(\Bun_G)$ by:
\begin{equation*}
\Eis(\S):=\p_{!}\circ\q^{!*}(\S).
\end{equation*}

\begin{cor} \label{goodEisprince}
The functor $\Eis$ commutes with Verdier duality and maps pure 
complexes to pure ones.
\end{cor}

\sssec{}  \label{intrhecke}

Let $\H_G$ be the {\it Hecke} stack, i.e. an $S$-point of $\H_G$ is a quadruple: 
$(x,\F_G,\F'_G,\beta)$, where $x\in\Hom(S,X)$, $\F_G$ and $\F'_G$ are $S$-points of $\Bun_G$
(i.e. $G$-bundles on $X\times S$) and $\beta$ is an isomorphism
$$\beta:\F_G|_{X\times S-\Gamma_x}\to \F'_G|_{X\times S-\Gamma_x},$$
where $\Gamma_x\subset X\times S$ is the graph of the map $x:S\to X$. 

\smallskip

We have three projections: $\pi:\H_G\to X$, $\hr_G:\,\H_G\to \Bun_G$ and 
$\hl_G:\,\H_G\to \Bun_G$ 
$$\pi(x,\F_G,\F'_G,\beta)=x;\,\,\hr_G(x,\F_G,\F'_G,\beta)=\F'_G;\,\, 
\hl_G(x,\F_G,\F'_G,\beta)=\F_G.$$
i.e. we have a diagram:

$$\Bun_G\times X \overset{\hl_G\times\pi}\longleftarrow 
\H_G \overset{\hr_G}\longrightarrow \Bun_G.$$

For a dominant coweight $\lambda$ of $G$ we introduce a closed substack $\Hb_G^{\lambda}$ of $\H_G$ as follows:
By definition, a quadruple $(x,\F_G,\F'_G,\beta)$ belongs to $\Hb_G^{\lambda}$
if and only if for every $G$-module $\V$ such all its weights are 
$\leq\lambdach$ for some $\lambdach\in\check\Lambda_G^{+}$,
\begin{equation*}  \label{heckelambda}
\V_{\F_G}(-\langle\lambda,\lambdach\rangle\cdot \Gamma_x)\subset\V_{\F'_G}\subset
\V_{\F_G}(-\langle \w_0(\lambda),\lambdach\rangle\cdot \Gamma_x).
\end{equation*}
(Of course, the second inclusion is a corollary of the first one for
the dual representation $\V^*$.)
As we will see later, both maps $\hl_G\times\pi$ and $\hr_G\times\pi$ when restricted to $\Hb_G^{\lambda}$
are proper and, in addition, are locally-trivial fibrations in the smooth topology (cf. \secref{conventions}).

\smallskip

For instance, when $G=T$, the projection $\hl_G\times\pi$ defines an isomorphism 
$\Hb_T^{\mu}\to \Bun_T\times X$ and the composition
$\Bun_T\times X\simeq\Hb^\mu_T\overset{\hr_T}\longrightarrow \Bun_T$
sends 
$$(\F_T,x)\in \Bun_T\times X\longrightarrow 
\F'_T:=\F_T(-\mu\cdot x)\in\Bun_T.$$ 

\medskip

Each $\Hb_G^{\lambda}$ is an algebraic stack locally of finite type and let $\IC_{\Hb_G^{\lambda}}$
denote the intersection cohomology sheaf on it. We define the Hecke functor 
$H^{\lambda}_G:\,\Sh(\Bun_G)\to \Sh(\Bun_G\times X)$ by
\begin{equation*}
\S\to (\hl_G\times\pi)_{!}
(\hr_G{}^*(\S)\otimes \IC_{\Hb_G^{\lambda}})
\otimes(\Ql(\frac{1}{2})[1])^{\otimes-\dim(\Bun_G)}.
\end{equation*}

\smallskip

As the map $\hl_G\times \pi:\Hb_G^\lambda\to\Bun_G$ is proper and 
$\hr_G:\Hb_G^\lambda\to\Bun_G$ is a fibration, the functor $H^{\lambda}_G$ 
commutes with Verdier duality.

\medskip

Recall that for $\lambda\in\Lambda_G^{+}$, $V^{\lambda}$ denotes the corresponding irreducible
representation
of the Langlands dual group $\check G$ and for a coweight $\mu$, $V^{\lambda}(\mu)$ is
the corresponding weight subspace of $V^\lambda$.

\begin{thm} \label{commutewithHecke}
For every $\lambda\in\Lambda_G^{+}$ and $\S\in\Sh(\Bun_T)$ we have a functorial isomorphism:
$$H^{\lambda}_G\circ\Eis(\S)\simeq 
{\underset{\mu\in\Lambda}\oplus}(\Eis\boxtimes\on{id})\circ H^{\mu}_T(\S)
\otimes V^{\lambda}(\mu),$$
where we used the notation $\Eis\boxtimes\on{id}$ for the corresponding functor $\Sh(\Bun_T\times X)\to
\Sh(\Bun_G\times X)$, and for $\mu\in\Lambda$,
$H^{\mu}_T$ denotes the corresponding Hecke functor for the group $T$.
\end{thm}

\sssec{} 

The Hecke functors can be composed in the following way:
If $\lambda_1$ and $\lambda_2$ are two elements of $\Lambda_G^{+}$, we define the functor
$$H^{\lambda_1}_G\star H^{\lambda_2}_G:\Sh(\Bun_G)\to \Sh(\Bun_G\times X)$$ by the rule:
$$H^{\lambda_1}_G\star H^{\lambda_2}_G(\S)=(\on{id}\boxtimes\Delta_X^{*})
((H^{\lambda_1}_G\boxtimes\on{id})\circ H^{\lambda_2}_G(\S))\otimes
(\Ql(\frac{1}{2})[1])^{\otimes-1}$$
where $\Delta_X:X\to X\times X$ denotes the diagonal embedding.

It is well-known (\cite{BD}) that there exists a canonical isomorphism of functors:
$$H^{\lambda_1}_G\star H^{\lambda_2}_G(\S)\simeq 
\underset{\lambda\in\Lambda_G^{+}}\oplus H_G^\lambda(\S)\otimes \Hom_{\check G}
(V^{\lambda},V^{\lambda_1}\otimes V^{\lambda_2}).$$

An additional property of the isomorphism of \thmref{commutewithHecke} is that for 
$\S\in \Sh(\Bun_T)$ the two isomorphisms:
\begin{align*}
&(H^{\lambda_1}_G\star H^{\lambda_2}_G)\circ\Eis(\S)\simeq \underset{\mu_1,\mu_2}\oplus 
(\Eis\boxtimes\on{id})\circ 
(H^{\mu_1}_T\star H^{\mu_2}_T)(\S)\otimes
V^{\lambda_1}(\mu_1)\otimes V^{\lambda_2}(\mu_2) \simeq \\
&\simeq \underset{\lambda}\oplus\underset{\mu}\oplus
(\Eis\boxtimes\on{id})\circ H^{\mu}_T(\S)\otimes
\Hom_{\check G}(V^\lambda,V^{\lambda_1}\otimes V^{\lambda_2})\otimes V^{\lambda}(\mu)
\end{align*}
and
\begin{align*}
&H^{\lambda_1}_G\star H^{\lambda_2}_G\circ\Eis(\S)\simeq \underset{\lambda}
\oplus\, H_G^\lambda\circ \Eis(\S)\otimes \Hom_{\check G}
(V^\lambda,V^{\lambda_1}\otimes V^{\lambda_2})\simeq \\
&\simeq \underset{\lambda}\oplus\underset{\mu}\oplus
(\Eis\boxtimes\on{id})\circ H^{\nu'}_T(\S)
\otimes\Hom_{\check G}(V^\lambda,V^{\lambda_1}\otimes V^{\lambda_2})\otimes V^{\lambda}(\mu)
\end{align*}
coincide.

\sssec{}  \label{intreq}

For a coroot $\alpha\in\Delta^+$ of $G$, consider the corresponding projection 
$T\to T/{\mathbb G}_m$. This map of groups gives rise to a map of stacks
$$\f^{\alpha}:\Bun_T\to \Bun_{T/{\mathbb G}_m}.$$

\smallskip

We will call a sheaf $\S$ on $\Bun_T$ regular if 
$$\f^{\alpha}_{!}(\S)=0,\,\,\,\forall\alpha\in\Delta^+$$ (this definition
makes sense even though the map $\f^{\alpha}$ is non-representable). 
We let $\Sh(\Bun_T)^{reg}\subset\Sh(\Bun_T)$ denote the full
triangulated subcategory of regular sheaves.

\medskip

Observe now that the action of the Weyl group on $T$ gives rise to a $W$-action on the stack 
$\Bun_T: \F_T\to\F_T^{\w}$. We define a twisted $W$-action on $\Bun_T$ as follows:

For an element $\lambda\in\Lambda$, consider the $T$-bundle $\Omega^{\lambda}$ on $X$ induced from 
the canonical line bundle on $X$
by means of the homomorphism ${\mathbb G}_m\to T$ corresponding to $\lambda$.

\smallskip

We set
\begin{equation*}
\w\cdot \F_T=\F^{\w}_T\otimes \Omega^{\w(\rho)-\rho}.
\end{equation*}
Similarly, for $\w\in W$ we will denote by $\S\to \w\cdot \S$ the direct image functor 
$\Sh(\Bun_T)\to \Sh(\Bun_T)$ corresponding to the new action of $\w$ on $\Bun_T$. 

It is easy to see that for every $\w\in W$ the functor $\S\to \w\cdot \S$ preserves
the subcategory $\Sh(\Bun_T)^{reg}\subset \Sh(\Bun_T)$. Let $N(\check T)$ denote 
the normalizer of $\check T$ in $\check G$.

\smallskip

\begin{thm} \label{functeq}
For each $\on{w}\in W$ with a choice of a lift to an element 
$\wt{\on{w}}\in N(\check T)\subset \check G$, there exists a functorial isomorphism
$$\Eis(\on{w}\cdot\S)\overset{f.eq}\longrightarrow \Eis(\S),\,\S\in \Sh(\Bun_T)^{reg}.$$
\end{thm}

\sssec{}  \label{compfuncteq}
The isomorphisms of \thmref{functeq} and \thmref{commutewithHecke} are compatible
with the following way: the two isomorphisms 
\begin{equation*}    \label{feandhecke}
H_G^{\lambda}\circ \Eis(\w\cdot\S)\rightrightarrows
{\underset{\nu\in\Lambda}\oplus}(\Eis\boxtimes\on{id})\circ H^{\nu}_T(\S)
\otimes V^{\lambda}(\w(\nu))
\end{equation*}
(see below) coincide, where the first isomorphism is simply
\begin{align*}
&H_G^{\lambda}\circ\Eis(\w\cdot\S)\overset{f.eq}\longrightarrow H_G^{\lambda}\circ\Eis(\S)\simeq 
{\underset{\nu\in\Lambda}\oplus}(\Eis\boxtimes\on{id})\circ H^{\nu}_T(\S)
\otimes V^{\lambda}(\nu)\simeq   \\
&{\underset{\nu\in\Lambda}\oplus}(\Eis\boxtimes\on{id})\circ H^{\nu}_T(\S)
\otimes V^{\lambda}(\w(\nu)),
\end{align*} 
(the last arrow comes from the map $V^{\lambda}(\nu)\to V^{\lambda}(\w(\nu))$ given
by the action of $\wt{\w}\in\check G$), and the second isomorphism is the composition:
\begin{align*}
&H_G^{\lambda}\circ\Eis(\w\cdot\S)\simeq {\underset{\nu\in\Lambda}\oplus}(\Eis\boxtimes\on{id})\circ 
H^{\nu}_T(\w\cdot\S)\otimes V^{\lambda}(\nu)\simeq \\
&{\underset{\nu\in\Lambda}\oplus}(\Eis\boxtimes\on{id})\circ 
(\w\cdot\boxtimes\on{id})\circ H^{\w^{-1}(\nu)}_T(\S)\otimes V^{\lambda}(\nu)\overset{f.eq}\longrightarrow
{\underset{\nu\in\Lambda}\oplus}(\Eis\boxtimes\on{id})\circ 
H^{\w^{-1}(\nu)}_T(\S)\otimes V^{\lambda}(\nu) \\
&\simeq{\underset{\nu'\in\Lambda}\oplus}(\Eis\boxtimes\on{id})\circ H^{\nu'}_T(\S)
\otimes V^{\lambda}(\w(\nu')),
\end{align*}
where the second arrow used the isomorphism 
$H^{\nu}_T(\w\cdot\S)\simeq (\w\cdot\boxtimes\on{id})\circ H^{\w^{-1}(\nu)}_T(\S)$,
which holds for any $\S\in \Sh(\Bun_T)$.

\noindent{\it Remark.}
Our construction of the isomorphism of functors 
$$\Eis(\w\cdot\S)\overset{f.eq}\longrightarrow\Eis(\S):\Sh(\Bun_T)^{reg}\to\Sh(\Bun_G)$$
of \thmref{functeq} involves an additional choice of a representation of $\w$ as 
a product of simple reflections. At the moment, we can prove that this isomorphism is 
independent of such a representation only when $G$ is simply-laced.

\begin{conj}
For $\S\in \Sh(\Bun_T)^{reg}\cap \Perv(\Bun_T)$, 
the object $\Eis(\S)$ is a perverse sheaf. Moreover, if $\S$ is irreducible, then so is $\Eis(\S)$.
\end{conj}

\noindent{\it Remark.}
When $G$ has semi--simle rank $1$, or when $G=GL(n)$ and $\S$ is a local system,
this conjecture is in fact a theorem (cf. \cite{Ga}). 

\ssec{Relation with the classical theory}

Let now $E_{\check T}$ be a $\check T$-local system on $X$. 
For $\lambda\in\Lambda$ we will denote by 
$E^\lambda_{\check T}$ the $1$-dimensional local system on $X$ induced from $E_{\check T}$ by means 
of the homomorphism $\check T\overset{\lambda}\to {\mathbb G}_m$.

\sssec{}

The (geometric) abelian class field theory produces from $E_{\check T}$
a $1$-dimensional local system, which we will denote by $\Aut'_{E_{\check T}}$ on $\Bun_T$, 
whose fiber at $\F_T=\F^0_T(\underset{i}\Sigma \,\mu_i\cdot x_i)$ is 
$\underset{i}\otimes (E^{-\mu_i}_{\check T})_{x_i}$, where the subscript ``$x_i$''
means ``fiber at $x_i$''.

\smallskip

Let $\sqrt{(\Aut'_{E_{\check T}})_{\Omega_X^{2\rho}}}$ be a $1$-dimensional vector space such that
$(\sqrt{(\Aut'_{E_{\check T}})_{\Omega_X^{2\rho}}})^{\otimes 2}\simeq
(\Aut'_{E_{\check T}})_{\Omega_X^{2\rho}}$. We define a perverse sheaf 
$\Aut_{E_{\check T}}$ on $\Bun_T$ as a 
tensor product of $\Aut'_{E_{\check T}}$ with the $1$--dimensional vector 
space $\sqrt{(\Aut'_{E_{\check T}})_{\Omega_X^{2\rho}}}\otimes
(\Ql(\frac{1}{2}))^{\otimes (g-1)\dim(T)}$,
placed in the cohomological degree $(1-g)\cdot\dim(T)$.

\smallskip

By construction, the perverse sheaf $\Aut_{E_{\check T}}$ has the following property with respect 
to the twisted $W$-action on $\Bun_T$:
$$\w\cdot \Aut_{E_{\check T}}\simeq \Aut_{E^{\w}_{\check T}},$$
where $E^{\w}_{\check T}$ is a $\check T$-local system induced from $E_{\check T}$ by means of 
$\check T\overset{\w}\longrightarrow \check T$.

\smallskip

It is clear from the definitions that the perverse sheaf $\Aut_{E_{\check T}}$ is a 
Hecke eigen-sheaf on $\Bun_T$ with respect to 
$E_{\check T}$, i.e. for $\lambda\in\Lambda_T$, we have:
$$H_T^{\mu}(\Aut_{E_{\check T}})\simeq 
\Aut_{E_{\check T}}\boxtimes E^\mu_{\check T}\otimes \Ql(\frac{1}{2})[1].$$

\medskip

Let $E_{\check G}$ denote the $\check G$-local system on $X$ induced from $E_{\check T}$ by means 
of the canonical
embedding $\check T\hookrightarrow\check G$. \thmref{commutewithHecke} yields the following result:

\begin{thm} \label{exHecke}
The object $\Eis(\Aut_{E_{\check T}})\in \Sh(\Bun_G)$ is a Hecke eigen-sheaf on $\Bun_G$ with 
respect to $E_{\check G}$, i.e. for each $\lambda\in\Lambda^{+}_G$ we have:
$$H^{\lambda}_G\circ\Eis(\Aut_{E_{\check T}})\simeq 
\Eis(\Aut_{E_{\check T}})\boxtimes V^{\lambda}_{E_{\check G}}\otimes 
\Ql(\frac{1}{2})[1],$$ where $V^{\lambda}_{E_{\check G}}$ 
is a local system on $X$ associated to $E_{\check G}$ and the $\check G$-module $V^{\lambda}$.
\end{thm}

\sssec{}

Assume now that the $\check T$-local system $E_{\check T}$ is regular, i.e. 
that for every $\alpha\in\Delta^+$ the $1$-dimesnional
local system $E^\alpha_{\check T}$ on $X$ is non-trivial. This condition is equivalent to the fact that 
the perverse sheaf $\Aut_{E_{\check T}}$ belongs to $\Sh(\Bun_T)^{reg}$.

\smallskip

For $\w\in W$ let us choose a representative $\wt{\w}\in N(\check T)$. Note that this gives 
an isomorphism $E_{\check G}\simeq \on{Ind}_{\check T}^{\check G}(E^{\w}_{\check T})$.
By applying \thmref{functeq}, we obtain the following assertion:

\begin{thm} \label{exfuncteq}
For $\wt{\on{w}}\in N(\check T)$ and $E_{\check G}$ as above there is an isomorphism
(functional equation):
$$\Eis(\Aut_{E_{\check T}})\overset{\text{f.eq}}\simeq \Eis(\Aut_{E^{\on{w}}_{\check T}}).$$
Moreover, the diagram
$$
\CD
H_G^{\lambda}\circ\Eis(\Aut_{E_{\check T}}) @>>> 
\Eis(\Aut_{E_{\check T}})\boxtimes V^\lambda_{E_{\check G}}\otimes\Ql(\frac{1}{2})[1]  \\
@V{H_G^{\lambda}(\text{f.eq})}VV     @V{\text{f.eq}\boxtimes \on{id}}VV  \\
H_G^{\lambda}\circ\Eis(\Aut_{E^{\on{w}}_{\check T}}) @>>> \Eis(\Aut_{E^{\on{w}}_{\check T}})\boxtimes 
V^\lambda_{E_{\check G}}\otimes\Ql(\frac{1}{2})[1] 
\endCD
$$
commutes.
\end{thm}

\medskip

The stack $\Bun_T$ splits into connected components $\Bun_T^\mu$, 
numbered by the elements of $\Lambda$: by definition $\F^0_T(\Sum\,\mu_i\cdot x_i)\in\Bun_T^{-\mu}$
if $\Sum\,\mu_i=\mu$. Let $\Aut_{E_{\check T}}^\mu$ be the direct summand of $\Aut_{E_{\check T}}$ 
concentrated on $\Bun_T^\mu$.

\thmref{exfuncteq} implies the following assertion, which is hard to see directly:

\begin{cor}
Let $E_{\check T}$ be a regular $\check T$-local system and let $U\subset\Bun_G$ be an open substack
of finite type. Then the sheaves $\Eis(\Aut_{E_{\check T}}^\mu)|_U$ are zero except for
finitely many $\mu\in\Lambda$.
\end{cor}

\begin{proof}

The proof is a combination of \thmref{exfuncteq} and of the following (obvious) statement:

\begin{lem} \label{finite}
Let $U\subset\Bun_G$ be an open sub-stack of finite type. Then the set of $\mu\in\Lambda$ for which 
the intersection
$\p^{-1}(U)\cap \q^{-1}(\Bun_T^\mu)$ is non-empty has the form
$$\mu\geq\mu',$$
where $\mu'$ is some fixed element of $\Lambda$.
\end{lem}

\end{proof}

\sssec{}

\thmref{exfuncteq} implies the following special case of the Langlands' conjecture:

\smallskip

Let $\K$ be the field of rational functions on $X$ and let $G_{\AA}$, 
(resp., $G_{\OO}\subset G_{\AA}$, $G_{\K}\subset G_{\AA}$)
be the corresponding adele group (resp., the group of integral points of $G_{\AA}$, the group of 
rational points of $G_{\AA}$).
Recall that the double quotient $G_{\OO}\backslash G_{\AA}/G_{\K}$ identifies with the set of 
isomorphism classes of objects of
the category $\Bun_G({\mathbb F}_q)$, i.e. with the set of isomorphism classes of $G$-bundles over $X$.

\medskip

\begin{thm} \label{langlands}
Let $E_{\check G}$ be an irreducible $\check G$-local system on $X$, such that 
if we ``forget'' the Weil structure on $E_{\check G}$, it admits 
a reduction to $\check T\subset\check G$. Then there exists a spherical 
automorphic function on $G_{\AA}$ whose Langlands' parameters correspond to $E_{\check G}$.
\end{thm}

\begin{proof}

Let $\ol X\simeq X\underset{\on{Spec}(\Fq)}\times \on{Spec}(\Fqb)$ and let
$\overline E_{\check G}$ denote $E_{\check G}$, viewed as a local system just on $\ol X$.
Let $Fr$ denote the Frobenius acting on $\overline X$ and on $\overline{\Bun_G}:=
\Bun_G\underset{\on{Spec}({\mathbb F}_q)}\times \on{Spec}(\overline{{\mathbb F}_q})$.

\smallskip

To construct the sought-for automorphic function it is sufficient to produce an object 
$\Aut_{\overline E_{\check G}}\in \Sh(\overline{\Bun_G})$ together with a Weil structure on it, which 
is a Hecke eigen-sheaf with respect to $\overline E_{\check G}$ in 
a $Fr$-compatible way. In other 
words, we need that there exists an isomorphism:
$$
Fr^*(\Aut_{\overline E_{\check G}})\to \Aut_{\overline E_{\check G}}
$$ 
such that for each $\lambda\in\Lambda_G^{+}$ the diagram
$$
\CD
Fr^*\circ H_G^{\lambda}(\Aut_{\overline E_{\check G}}) @>>> 
Fr^*(\Aut_{\overline E_{\check G}})\boxtimes V^\lambda_{Fr^*(\overline E_{\check G})} \\
@VVV     @VVV  \\
H_G^{\lambda}(\Aut_{\overline E_{\check G}}) @>>> \Aut_{\overline E_{\check G}}\boxtimes
V^\lambda_{\overline E_{\check G}}
\endCD
$$
is commutative.

\medskip

Let $\overline E_{\check T}$ be a reduction of the $\check G$-local system 
$\overline E_{\check G}$ over $\overline X$
to $\check T\subset \check G$ (such a reduction exists according to our assumtion on 
$E_{\check G}$). We have now the following assertion: \footnote{We are obliged to G.~Prasad who has
explained to us the proof of this result. He has moreover authorized us to 
reproduce his argument, which we shall do in the Appendix.}

\begin{prop} \label{prasad}
\smallskip

{\em (a)}
Let $\overline E_{\check T}$ be as above and assume (as in the formulation of \thmref{langlands}) 
that $E_{\check G}$ is irreducible. Then the $\check T$-local system $\overline E_{\check T}$ 
is regular.

{\em (b)} There exists an element $\on{w}\in W$ and its lift $\wt{\on{w}}\in N(\check T)$ such that
$Fr^*(\overline E_{\check T})\simeq \overline E_{\check T}^{\on{w}}$ and such the identification
$$\on{Ind}_{\check T}^{\check G}(\overline E_{\check T})\simeq \ol E_{\check G}\simeq
Fr^*(\ol E_{\check G})\simeq \on{Ind}_{\check T}^{\check G}(Fr^*(\overline E_{\check T}))\simeq
\on{Ind}_{\check T}^{\check G}(\overline E_{\check T}^{\on{w}})$$
is induced by $\wt{\on{w}}$.

\end{prop}

We define now the object $\Aut_{\overline E_{\check G}}$ as 
$\Eis(\Aut_{\overline E_{\check T}})$. Using \thmref{exfuncteq} and \propref{prasad}
we construct an isomorphism $Fr^*(\Eis(\Aut_{\overline E_{\check T}}))\simeq \Eis(\Aut_{\overline E_{\check T}})$
that corresponds to the above $\wt{\w}$. \thmref{exfuncteq} guarantees that 
$\Aut_{\overline E_{\check G}}$ possesses all the required properties.

\end{proof}

\sssec{}

Let $E_{\check T}$ be again an arbitrary (not necessarily regular) $\check T$-local system on $X$. 
We will now investigate the connection between the function on 
$G_{\OO}\backslash G_{\AA}/G_{\K}\simeq\Bun_G({\mathbb F}_q)$ corresponding to
$\Eis(\Aut_{E_{\check T}})\in \Sh(\Bun_G)$ and the classical Eisenstein series constructed 
starting with $E_{\check T}$.

\smallskip

For $E_{\check T}$ as above and $\mu\in\Lambda$ consider the object 
$\Eis'(\Aut_{E_{\check T}}^\mu)$ defined as
$$\po_{!}\circ \qo^{!*}(\Aut_{E_{\check T}}^\mu).$$

\smallskip

We let $\on{Funct}(\Eis(\Aut_{E_{\check T}}^\mu))$ 
(resp., $\on{Funct}(\Eis'(\Aut_{E_{\check T}}^\mu))$) denote the functions on the set
$\Bun_G({\mathbb F}_q)$ corresponding to $\Eis(\Aut_{E_{\check T}})$ and 
$\Eis'(\Aut_{E_{\check T}})$, respectively. 

\smallskip

The following result is a corollary of a description of stalks of the intersection cohomology 
sheaf on $\BunBb$ (\thmref{compute1}):
\footnote{The proof of \thmref{compute1} will be given in a forthcoming paper \cite{BFGM}. 
For the case
$X={\mathbb P}^1$ this result is proved in \cite{FFKM}. Note, however, that  
all the results of this paper except for \thmref{comparefinite} and \thmref{fullcompare} 
are independent of \thmref{compute1}.}

\begin{thm} \label{comparefinite}
The function $\on{Funct}(\Eis(\Aut_{E_{\check T}}^\mu))$ equals
$$\underset{\{n_\alpha\}\in{\Delta^+}^{\NN}}\Sigma
\on{Funct}(\Eis'(\Aut_{E_{\check T}}^{\mu-\underset{\alpha\in\Delta^+}\Sigma n_\alpha\cdot\alpha}))\cdot
\underset{\alpha\in\Delta^+}\Pi\on{Tr}
(Fr,H^{\cdot}(X^{(n_{\alpha})},
(E^{\alpha}_{\check T})^{(n_{\alpha})})\otimes\Ql(n_{\alpha})),$$
where $X^{(m)}$ denotes the $m$-th symmetric power of $X$ and for a local system $E$, 
$E^{(m)}$ denotes the $m$-th symmetric power of $E$.
\end{thm}

We will now reformulate \thmref{comparefinite} by encoding the information contained in it into a 
generating series.

\smallskip

Consider the group--ring $\Ql[\Lambda]$, which is the same, of course, as the ring of functions 
on the torus $\check T$;
for $\mu\in\Lambda$ we shall denote by $t^\mu$ the corresponding element of $\Ql[\Lambda]$.
We will form a completed ring $\widehat{\Ql[\Lambda]}$ by allowing infinite expressions of the form
$$\underset{\mu}\Sigma \,a_\mu\cdot t^\mu,$$ if $\mu$ runs over a sub-set of $\Lambda$ of 
the form $\mu\geq\mu'$, where $\mu'$ is some fixed element of $\Lambda$.

\smallskip

The classical Eisenstein series can be thought of as a function on $\Bun_G({\mathbb F}_q)$ with 
values in $\widehat{\Ql[\Lambda]}$, equal to
\begin{equation*}
\Eis_{cl}(\Aut_{E_{\check T}})(t):=
\underset{\mu\in\Lambda}\Sigma \on{Funct}(\Eis'(\Aut_{E_{\check T}}^{\mu}))\cdot t^\mu.
\end{equation*}

\smallskip

Similarly, consider the ``modified'' Eisenstein series
\begin{equation*}
\Eis_{mod}(\Aut_{E_{\check T}})(t):=\underset{\mu\in\Lambda}\Sigma  
\on{Funct}(\Eis(\Aut_{E_{\check T}}^{\mu}))\cdot t^\mu,
\end{equation*}
viewed again as a function on $\Bun_G({\mathbb F}_q)$ with values in $\widehat{\Ql[\Lambda]}$ 
(cf. \lemref{finite}).

\smallskip

Finally, for $\alpha\in\Delta^+$ we introduce the (abelian) L-series 
$L(E_{\check T},\alpha,t)$ to be the element 
of $\widehat{\Ql[\Lambda]}$ equal 
to $\underset{n\in\NN}\Sigma \on{Tr}(Fr,H^{\cdot}(X^{(n)},
(E^{\alpha}_{\check T})^{(n)})\otimes\Ql(n))
\cdot t^{n\cdot\alpha}.$

\medskip

We have:

\begin{thm} \label{fullcompare}
$$\Eis_{mod}(\Aut_{E_{\check T}})(t)=\Eis_{cl}(\Aut_{E_{\check T}})(t)\cdot \underset{\alpha\in\Delta^+}
\prod L(E_{\check T},\alpha,t).$$
\end{thm}

When $G=GL(n)$, the above result has been established by G.~Laumon in \cite{La}.

\medskip

\noindent{\it Remark.}
It is well-known the the power series $\Eis_{cl}(\Aut_{E_{\check T}})(t)$ 
satisfies the functional equation of the form
\begin{equation*}
\Eis_{cl}(\Aut_{E_{\check T}})(t)=\Eis_{cl}(\Aut_{E^{\w}_{\check T}})(t)\cdot 
\frac{\prod_{\alpha\in\Delta^+} 
L(E_{\check T}^{\w},\alpha,t)}{\prod_{\alpha\in\Delta^+}L(E_{\check T},\alpha,t)},
\end{equation*}
where $L(E_{\check T},\alpha,t)$ are precisely the L-functions that appear 
in the formulation of \thmref{fullcompare}. However, \thmref{exfuncteq}
implies that $\Eis_{mod}(\Aut_{E_{\check T}})(t)=
\Eis_{mod}(\Aut_{E^{\w}_{\check T}})(t)$ (if $E_{\check T}$ is regular). Therefore,
the fact that the ratio between $\Eis_{mod}$ and $\Eis_{cl}$ is the mentioned above product of 
L-functions is very natual from
this point of view.

It is quite remarkable that the Eisenstein series $\Eis_{mod}(\Aut_{E_{\check T}})(t)$, 
which is more natural than
$\Eis_{cl}(\Aut_{E_{\check T}})(t)$ from the geometric point of view 
(e.g. it comes from an object of $\Sh(\Bun_G)$, whose 
construction is self-dual in the Verdier sense) incorporates the L-function.

\ssec{Non-principal Eisenstein series}  \label{nonprinceseries}

We will now consider the case when $B$ is replaced by a parabolic $P$.

\sssec{}

As in the previous subsection we define functor $\qw_P^{!*}:\Sh(\Bun_M)\to \Sh(\BunPbw)$ by
\begin{equation*}
\S\to \IC_{\BunPbw}\otimes\qw_P^{*}(\S)\otimes(\Ql[1](\frac{1}{2}))^{\otimes-\dim(\Bun_M)}.
\end{equation*}
Again, $\qo_P^{!*}$ will denote the functor $$\jw_P^*\circ\qw_P^{!*}:\, \Sh(\Bun_M)\to \Sh(\Bun_P).$$
Since the map $\qo_P:\Bun_P\to\Bun_M$ is smooth, $\qo_P^{!*}$ is essentially the
ordinary pull-back.

The following theorem is parallel to \thmref{goodtopullbackprince}:

\begin{thm}  \label{goodtopullbacknonprince}

\smallskip

(a) The functor $\qw_P^{!*}$ is exact and commutes with Verdier duality.

\smallskip

(b) If $\S$ is a perverse sheaf on $\Bun_M$, we have: 
$\qw_P^{!*}(\S)\simeq \jw_P{}_{!*}\circ\qo_P^{!*}(\S)$.
\end{thm}

We define the functor $\Eis_M^G:\,\Sh(\Bun_M)\to \Sh(\Bun_G)$ by 
\begin{equation*}
\S\to \pw_P{}_{!}\circ \qw_P^{!*}(\S).
\end{equation*}

As in \corref{goodEisprince}, \thmref{goodtopullbacknonprince} 
insures that the functor $\Eis^G_M$ commutes
with Verdier duality and maps pure complexes into pure ones.

\sssec{}

Let $\S$ be a perverse sheaf on $\Bun_M$. One may attempt to give an 
alternative definition of the functor $\Eis^G_M$ as follows, using the compactification $\BunPb$:

\smallskip

We can first pull back $\S$ onto $\Bun_P$ by means of $\S\to {\qo}^{!*}_P(\S)$, 
then extend it by means of $j_P{}_{!*}$ to
obtain a perverse sheaf on $\BunPb$ and then take the direct image of the resulting sheaf 
with respect to the map
$\p_P:\BunPb\to\Bun_G$. It will not be in general true that the two functors yield the same result. 
However, one can single
out a class of perverse sheaves on $\Bun_M$ for which this will be true (cf. \secref{good}):

\medskip

For a reductive group $G$ a perverse sheaf $\S\in\Sh(\Bun_G)$ is said to be ``good'' if 
for any sequence $\lambda_1,...,\lambda_n$
of elements of $\Lambda_G^{+}$, the sheaf
$$(H_G^{\lambda_n}\boxtimes\on{id}^{n-1})\circ...\circ(H_G^{\lambda_2}\boxtimes\on{id})
\circ H_G^{\lambda_1}(\S)$$
on $\Bun_G\times X^n$ is perverse as well.

\smallskip
 
For instance, when $G$ is abelian, every perverse sheaf on $\Bun_G$ is ``good''. 
For a general $G$, an automorphic sheaf 
that corresponds in the sense of Drinfeld--Langlands to a $\check G$-local system on $X$ is ``good''.

\smallskip

We have the following assertion, proved in \secref{good}:

\begin{thm}  \label{compat}
Let $\S\in \Perv(\Bun_M)$ be ``good''. Then there is a canonical isomorphism
${\r_P}_{!}\circ\qw_P^{!*}(\S)\simeq {j_P}_{!*}\circ \qo_P^{!*}(\S)$.
\end{thm}

\begin{cor}
For a ``good'' perverse sheaf $\S$ on $\Bun_M$ the sheaf 
${\p_P}_{!}\circ j_P{}_{!*}\circ\qo_P^{!*}(\S)$ on $\Bun_G$ 
is canonically isomorphic to $\Eis^G_M(\S)$.
\end{cor}

\noindent{\it Remark.}
Let $G=GL(n)$. In this case there exists an open embedding of the stack 
$\Bun_{GL(n)}$ into the stack $Coh_n$ that classifies coherent sheaves
on $X$ of generic rank $n$. Following Laumon (cf. \cite{La}) one can consider a functor
$$\Eis^n_{n_1,...,n_k}:\Sh(Coh_{n_1}\times...\times Coh_{n_k})\to \Sh(Coh_{n}),$$
for $n=n_1+...+n_k$.

\smallskip

Now let $\S$ be a ``good'' perverse sheaf on $\Bun_{GL(n_1)}\times...\times \Bun_{GL(n_k)}$ and let
$\S'$ denote its Goresky-MacPherson extension on the whole of 
$Coh_{n_1}\times...\times Coh_{n_k}$. In this case 
we conjecture that the restriction of
$\Eis^n_{n_1,...,n_k}(\S')$ to $\Bun_{GL(n)}\subset Coh_n$ is canonically isomorphic to 
$\Eis_M^{GL(n)}(\S)$, where $M=GL(n_1)\times...\times GL(n_k)$.

\sssec{}

In the case of non-principal Eisenstein series we have the following 
generalization of \thmref{commutewithHecke}:

\begin{thm} \label{commutewithHeckenonprince}
For every $\lambda\in\Lambda_G^{+}$ and $\S\in\Sh(Bun_M)$ there is a functorial isomorphism
$$H^{\lambda}_G\circ\Eis_M^G(\S)\simeq {\underset{\nu\in\Lambda_M^+}\oplus}\, (\Eis_M^G\boxtimes\on{id})\circ
H^{\nu}_M(\S)
\otimes \on{Hom}_{\check M}(U^{\nu},V^{\lambda}).$$
\end{thm}

As in the case $P=B$ we can show in addition, that the isomorphism of functors in the above theorem
is compatible with the convolution of Hecke functors.

\medskip

Let $E_{\check M}$ be an $\check M$-local system on $X$ and let $\Aut_{E_{\check M}}$ 
be an automorphic perverse 
sheaf on $\Bun_M$ that satisfies a Hecke eigen-property with respect to $E_{\check M}$. 
Let $E_{\check G}$ denote a $\check G$-local system induced from $E_{\check M}$.

\begin{cor} \label{Heckeeigennon-prince}
The complex $\Eis_M^G(\Aut_{E_{\check M}})$ is a Hecke eigensheaf with respect to 
$E_{\check G}$.
\end{cor}

\sssec{}

Finally, we have the following compatibility between the functors of Eisenstein series:

\begin{thm} \label{compos}
For $\S\in\Sh(\Bun_T)$ there is a functorial isomorphism
$$\Eis_T^G\simeq \Eis_M^G\circ \Eis^M_T.$$ Moreover, this isomorphism of functors is compatible 
with the Hecke property 
in the sense that for every $\lambda\in\Lambda_G^{+}$ the diagram
$$
\CD
H_G^\lambda\circ \Eis_T^G(\S) @>>> H_G^\lambda\circ\Eis_M^G\circ \Eis^M_T(\S) \\
@VVV       @VVV   \\
\underset{\mu}\oplus \Eis_T^G\circ H_T^\mu(\S)\otimes V^\lambda(\mu) @>>> 
\underset{\nu}\oplus\underset{\mu}\oplus\Eis^G_M\circ\Eis_T^M\circ H_T^\mu(\S)
\otimes\on{Hom}_{\check M}(U^{\nu},V^{\lambda}) \otimes U^{\nu}(\mu)   
\endCD
$$
commutes.
\end{thm}

Although the assertion of \thmref{compos} seems very natural, the proof is quite non-trivial and
it should be thought of as the main technical result of this paper. In particular, \thmref{functeq} 
will be a rather easy consequence of \thmref{compos} via an explicit calculation in the case
when $G$ has rank $1$.


\section{Action of Hecke operators (the principal case)}

\ssec{The basic diagram}

The construction described below will play an essential role in this paper.  

\sssec{}  \label{bd}

Let $\lambda$ be a dominant coweight of $G$ and let us denote by 
$\ol{Z}$ the fiber product $\Hb_G^{\lambda}{\underset{\Bun_G}\times}\BunBb$,
where $\Hb_G^{\lambda}$ is mapped to $\Bun_G$ by 
means of the projection $\hr_G$.

\smallskip

\begin{propconstr}
There exists a map of stacks $\phi:\ol{Z}\to \BunBb\times X$ 
with $(\p\times\on{id})\circ \phi=(\hl_G\times\pi)\circ{}'\p$
\end{propconstr}

$$
\CD
\BunBb\times X @<\phi<< \ol{Z} @>{'\hr_G}>>  \BunBb  \\
  @V{\p\times\on{id}}VV        @V{'\p}VV     @V{\p}VV   \\
\Bun_G\times X   @<{\hl_G\times\pi}<<   \Hb_G^{\lambda}    @>{\hr_G}>>  \Bun_G
\endCD
$$
Note that the left square of this diagram is NOT Cartesian.

\begin{proof}

An $S$-point $z$ of $\ol{Z}$ consists by definition of the following data: $x\in \Hom(S,X)$; a pair of $G$-bundles 
$\F_G$ and $\F'_G$ on $X\times S$ identified with one another outside on $X\times S-\Gamma_x$;
a collection of line bundles $\L^{\lambdach}_{\F'_T}$, each  embedded as subsheaf 
into the corresponding $\V^{\lambdach}_{\F'_G}$, such that the Pl\"ucker relations hold. 

Recall that the fact that the pair $(\F_G,\F'_G)$ belongs to $\Hb_G^{\lambda}$ means, according to
\secref{intrhecke} that for each $\lambdach\in\check\Lambda_G^{+}$ we have:
\begin{equation*}
\V^{\lambdach}_{\F'_G}\subset \V^{\lambdach}_{\F_G}(\langle -\w_0(\lambda),\lambdach\rangle\cdot \Gamma_x).
\end{equation*}

Therefore, the embedding 
$\kappa'{}^{\lambdach}:\L^{\lambdach}_{\F'_T}\hookrightarrow \V^{\lambdach}_{\F'_G}$ gives rise to an embedding
$$\kappa^{\lambdach}:\L^{\lambdach}_{\F'_T}(\langle \w_0(\lambda),\lambdach\rangle\cdot \Gamma_x)
\hookrightarrow \V^{\lambdach}_{\F_G}$$
and we set by definition $\phi(z)$ to be the object of $\BunBb$ that corresponds to \newline
$\F_G$, $\F_T:=\F'_T(\w_0(\lambda)\cdot \Gamma_x)$ and the system of embeddings $\kappa^{\lambdach}$ described above.

\end{proof}

\sssec{}

For every element $\nu\in \Lambda_G^{\on{pos}}$ consider the closed embedding 
$$i_\nu:\BunBb\times X\hookrightarrow \BunBb\times X,$$
given by sending an object 
$(\F_G,\F_T,\kappa^{\lambdach}:\L^{\lambdach}_{\F_T}\hookrightarrow \V^{\lambdach}_{\F_G},x)$ to 
$$(\F_G,\,\F_T(-\nu\cdot \Gamma_x),\,\,
\L^{\lambdach}_{\F_T}(-\langle \nu,\lambdach\rangle\cdot \Gamma_x)\hookrightarrow \L^{\lambdach}_{\F_T}\hookrightarrow 
\V^{\lambdach}_{\F_G},x).$$

\begin{thm} \label{mainHeckeprince}
$\phi_{!}(\IC_{\ol{Z}})\simeq \underset{\nu\in\Lambda_G^{\on{pos}}}\oplus i_\nu{}_{!}
(\IC_{\BunBb\times X})\otimes V^{\lambda}(\w_0(\lambda)+\nu).$
\end{thm}

\medskip

\noindent{\it Remark.} 
An analogue of this theorem has been independently obtained by 
M.~Finkelberg and I.~Mirkovic. Our proof
is slightly different from the argument of \cite{FM} and in \secref{actnonprince} 
it will be generalized to the case when $B$ is replaced by a parabolic subgroup.
 
\sssec{}

Let us first show how \thmref{mainHeckeprince} implies \thmref{commutewithHecke}.

\begin{proof}(of \thmref{commutewithHecke})

For $\nu\in\Lambda$ let us denote by $m^\nu$ the map
$$\Bun_T\times X\simeq \Hb_T^\nu \overset{\hr_T}\longrightarrow \Bun_T.$$

By definition, for $\S\in\Sh(\Bun_T)$, the sheaf $m^\nu{}^*(\S)$ is canonically the same as $H^\nu_T(\S)$.

\begin{lem} \label{simpleHeckeT}
We have:

\smallskip

{\em (a)} The maps $\q\circ{}'\hr_G$ and 
$m^{\on{w}_0(\lambda)}\circ(\q\times\on{id})\circ\phi$ from $\ol{Z}$ to $\Bun_T$ coincide.

\smallskip

{\em (b)} For every $\nu$, the maps $(\q\times\on{id})\times i_\nu$ and 
$(m^\nu\times\on{id})\circ (\q\times \on{id})$
from $\BunBb\times X$ to $\Bun_T\times X$ coincide.
\end{lem}

\medskip

For the proof of \thmref{commutewithHecke}, observe that for an object $\S\in \Sh(\Bun_T)$, the sheaf
$\hr_G{}^*\circ\Eis(\S)$ on $\ol{Z}$ identifies (by base change and the projection formula) with
\begin{equation*}
'\p_{!}({}'\hr_G{}^*\circ\q^*(\S)\otimes {}'\hr_G{}^*(\IC_{\BunBb})\otimes {}'\p^*(\IC_{\Hb^{\lambda}_G}))\otimes 
(\Ql(\frac{1}{2})[1])^{\otimes-\dim(\Bun_T)}.
\end{equation*}

\smallskip

However, the map $\hr_G:\Hb_G^{\lambda}\to\Bun_G$ decomposes locally in the smooth topology with respect to $\Bun_G$
into a direct product, which implies that
\begin{equation*}
\IC_{\ol{Z}}\simeq {}'\p^*(\IC_{\Hb^{\lambda}_G})\otimes {}'\hr_G{}^*(\IC_{\BunBb})\otimes 
(\Ql(\frac{1}{2})[1])^{\otimes-\dim(\Bun_G)}.
\end{equation*}

\smallskip

Therefore,
$$\hr_G{}^*\circ\Eis(\S)\simeq {}'\p_{!}({}'\hr_G{}^*\circ\q^*(\S) \otimes \IC_Z)\otimes 
(\Ql(\frac{1}{2})[1])^{\otimes-\dim(\Bun_T)+\dim(\Bun_G)}$$
and hence
$$H_G^\lambda\circ \Eis(\S)\simeq (\p\times\on{id})_{!}\circ\phi_{!}
({}'\hr_G{}^*\circ\q^*(\S) \otimes \IC_Z)(\Ql(\frac{1}{2})[1])^{\otimes-\dim(\Bun_T)}.$$

\smallskip

Acording to point (a) of \lemref{simpleHeckeT} and the projection formula, the last expression can be rewritten as
\begin{equation*}
(\p\times\on{id})_{!}(\phi_{!}(\IC_{\ol{Z}})\otimes (\q\times\on{id})^*\circ m^{\w_0(\lambda)}{}^*(\S))\otimes 
(\Ql(\frac{1}{2})[1])^{\otimes-\dim(\Bun_T)}
\end{equation*}
and finally, by applying \thmref{mainHeckeprince} and the projection formula we obtain that the sheaf
$H_G^\lambda\circ \Eis(\S)$ identifies with the direct sum over $\nu\in\Lambda_G^{\on{pos}}$ of the expressions
\begin{align*}
&(\p\times\on{id})_{!}\circ i_\nu{}_{!}\,(i_\nu^*\circ (\q\times\on{id})^*\circ m^{\w_0(\lambda)}{}^*(\S)\otimes 
\IC_{\BunBb\times X})\otimes \\
&\otimes V^{\lambda}(\w_0(\lambda)+\nu)\otimes  (\Ql(\frac{1}{2})[1])^{\otimes-\dim(\Bun_T)}.
\end{align*}

Now, by \lemref{simpleHeckeT}(b), 
$$i_\nu^*\circ (\q\times\on{id})^*\circ m^{\w_0(\lambda)}{}^*(\S)
\simeq (\q\times\on{id})^*\circ H_T^{\w_0(\lambda)+\nu}(\S).$$
Since the maps $(\p\times\on{id})\circ i_\nu$ and $\p\times\on{id}$ 
from $\BunBb\times X\to\Bun_G\times X$ coincide, we obtain that $H_G^\lambda\circ \Eis(\S)$
can be identified with

\begin{align*}
&\underset{\nu\in\Lambda_G^{\on{pos}}}\oplus 
(\p\times\on{id})_{!}\,((\q\times\on{id})^*\circ H_T^{\w_0(\lambda)+\nu}(\S)\otimes 
\IC_{\BunBb\times X})\otimes V^{\lambda}(\w_0(\lambda)+\nu)\otimes \\
&\otimes(\Ql(\frac{1}{2})[1])^{\otimes-\dim(\Bun_T)} 
\simeq \underset{\nu\in\Lambda_G^{\on{pos}}}\oplus (\Eis\boxtimes\on{id})\circ H_T^{\w_0(\lambda)+\nu}(\S)
\otimes V^{\lambda}(\w_0(\lambda)+\nu),
\end{align*}
as required.

\end{proof}

The fact that the isomorphism 
$H_G^\lambda\circ \Eis(\S)\simeq\underset{\nu\in\Lambda}\oplus (\Eis\boxtimes\on{id})\circ H_T^{\nu}(\S)
\otimes V^{\lambda}(\nu)$ constructed above is compatible with the convolution of Hecke functors 
will be explained in the next section in a more general context.

Thus, modulo \thmref{mainHeckeprince}, we have established \thmref{commutewithHecke} together with
\thmref{exHecke}.

\ssec{Preliminaries on the affine Grassmannian}  \label{prelaffgrass}

In this subsection we will review several facts concerning the affine Grassmannian corresponding to the 
group $G$. These facts will be needed for the proof of \thmref{mainHeckeprince} as well as for the rest of the paper;
the proofs can be found in \cite{Lu} or in \cite{BD}.

\sssec{} \label{ag}

Let $x\in X$ be a point and let $\D_x$ (resp., $\D^*_x$) be the formal disc (resp., the formal
punctured disc) around $x$. Let $G(\O_x)$ (resp., $G(\K_x)$) be the group--scheme (resp., group--indscheme) 
that classifies maps from $\D_x$ (resp., from $\D^*_x$) to $G$.

The affine Grassmannian $\Gr_G$ is, by definition, the quotient $G(\K_x)/G(\K_x)$. (Sometimes, 
in order to emphasize the dependence on $x$ we will put a subscript ${\Gr_G}_x$.) In other words,
$\Gr_G$ is an indscheme that classifies the data of $(\F_G,\beta)$, where $\F_G$ is a $G$-bundle on $\D_x$
and $\beta$ is its trivialization: $\F_G\simeq \F^0_G$ on $\D_x^*$. (We leave it to the reader to formulate what
this means on the level of $S$-points, cf. \cite{Pu}).

Being an indscheme, $\Gr_G$ is a union of its closed finite--dimensional subschemes. Therefore,
the notion of a perverse sheaf on $\Gr_G$ makes sense. (By definition, every such perverse 
is supported on a finite--dimensional subscheme of $\Gr_G$.) 

For $\lambda\in\Lambda_G^{+}$ let $\grb_G^\lambda$ denote the closed 
sub-scheme of $\Gr_G$ that corresponds 
to pairs $(\F_G,\beta)$ as above for which
$\V_{\F^0_G}(-\langle\lambda,\lambdach\rangle\cdot x)\subset 
\V^{\lambdach}_{\F_G}$
for every $G$-module $\V$ such all its weights are 
$\leq \lambdach$. (Here $x$ denotes, of course, the closed point 
of $\D_x$.) 

By construction, each $\grb_G^\lambda$ is a finite-dimensional
projective variety, stable under the left $G(\O_x)$-action. We have:
$$\grb_G^{\lambda'}\subset \grb_G^\lambda \text{ if and only if } \lambda'\leq\lambda.$$
In addition, it is known (cf. \cite{Lu}) that $\grb_G^\lambda$ 
contains a unique dense $G(\O_x)$-orbit, denoted $\Gr_G^\lambda$.
If $t_x$ is a uniformizer at $x$, then $\Gr_G^\lambda=G(\O_x)\cdot t_x^\lambda$, where $t_x^\lambda$ 
is the corresponding point of $T(\K_x)\subset G(\K_x)$. 

Thus, we obtain a stratification of $\Gr_G$ by the sub-schemes $\Gr_G^\lambda$.
\footnote{In this paper, the term ``stratification'' will be used in the following weak sense: 

Let $Y$ be a scheme (resp., indscheme, stack, indstack) and let $Y_i$ be a collection of locally closed subschemes
of $Y$. We say that they form a stratification if $Y(\Fqb)$ is a disjoint union $\cup Y_i(\Fqb)$.}

\smallskip

For $\lambda\in\Lambda_G^+$ let $\A_G^\lambda$ denote the intersection cohomology sheaf $\IC_{\grb_G^\lambda}$.
Let $\Sph_G$ denote the full abelian subcategory of $\Perv(\Gr_G)$ formed by direct sums 
of the perverse sheaves $\A_G^\lambda$. By definition, every object of $\Sph_G$ is $G(\O_x)$--equivariant.

\medskip

Note that $G(\O_x)$--equivariant perverse sheaves on $\Gr_G$ can be thought of as ``perverse sheaves'' on the stack
that classifies triples $(\F^1_G,\F^2_G,\beta)$, where $\F^i_G$ are $G$-bundles on $X$ $\D_x$ 
and $\beta$ is an isomorphism $\F^1_G|_{\D^*_x}\to \F^2_G|_{\D^*_x}$. The latter
stack carries a natural involution given by the flip $\F^1_G\leftrightarrow \F^2_G$. This involution
gives rise to a covariant functor, denoted $\ast$, from the category of $G(\O_x)$--equivariant perverse 
sheaves $\Gr_G$ to itself. We have $\ast\A_G^\lambda\simeq \A_G^{-\w_0(\lambda)}$, hence, 
$\ast$ preserves the subcategory $\Sph_G$. 

\sssec{}      \label{convolution}

Now we will recall the definition of the convolution operation on the category $\Sph_G$. Consider
the indscheme $\Conv_G:=\Gr_G \overset{G(\O_x)}\times G(\K_x)$. By definition, it classifies
the data of $(\F_G,\F'_G,\wt{\beta},\beta')$, where $\F_G$ and $\F'_G$ are a $G$--bundles
on $\D_x$, $\wt{\beta}$ is an isomorphism $\F_G|_{\D^*_x}\simeq \F'_G|_{\D^*_x}$ and 
$\beta'$ is a trivialization $\F'_G|_{\D^*_x}\simeq \F^0_G|_{\D^*_x}$. There are two projections
$pr$ and $pr'$ from $\Conv_G$ to $\Gr_G$:

\noindent For $g_1\times g_2\in G(\K_x)\times \Gr_G$,
$$pr(g_1\times g_2)=g_1\cdot g_2\,\on{mod}\, G(\O_x) \text{   and   }
pr'(g_1\times g_2)=g_2\,\on{mod}\, G(\O_x).$$

In the functorial language, $pr'$ and $pr$ act as follows:
$$pr'(\F_G,\F'_G,\wt{\beta},\beta')=(\F'_G,\beta') \text{ and }
pr(\F_G,\F'_G,\wt{\beta},\beta')=(\F_G,\wt{\beta}\circ \beta').$$

\smallskip

If $\S^1$ is a $G(\O_x)$--equivariant perverse sheaf on $\Gr_G$ and $\S_2$ is 
an arbitrary perverse sheaf on $\Gr_G$,
we can form their twisted external product $\S^1\widetilde{\boxtimes}\S^2$ 
(cf. \secref{conventions}), 
which is a perverse sheaf on $\Conv_G$.
The convolution $\S^1\star \S^2\in\Sh({}\Gr_G)$ is defined as 
$pr_{!}(\S^1\tboxtimes\S^2)$. We have the following result:

\begin{prop} \label{lusztig}
In the above situation $\S^1\star \S^2$ is a perverse sheaf. 
\end{prop}

\noindent{\it Remark.}
\propref{lusztig} has been first established in \cite{Lu} under the 
assumption that $\S_2$
is also $G(\O_x)$--equivariant. The proof in the 
general case is given in \cite{Pu}.

\medskip

Assume now that both $\S_1$ and $\S_2$ belong to $\Sph_G$. The \propref{lusztig} 
and the decomposition theorem (using the fact the each $\Gr_G^\lambda$ is simply-connected)
imply that $\S^1\star \S^2$ is also an object of $\Sph_G$. 
In addition, it follows from the definitions that there is a canonical
isomorphism $\ast(\S_1\star\S_2)\simeq (\ast\S_2)\star (\ast\S_1)$.

\sssec{}  \label{globalgrass}

We will now reformulate the definition of $\Gr_G$ using the global curve $X$. 
We claim that $\Gr_G$ is the indscheme that classifies the data of
$(\F_G,\beta)$, where $\F_G$ is a $G$--bundle on the curve $X$ and 
$\beta$ is its trivialization on $X-x$. 

Indeed, to $(\F_G,\beta)$ as above, we can attach $\F_G|_{\D_x}$, which is a 
$G$--bundle on $\D_x$
and $\beta$ defines its trivialization over $\D^*_x$. Now, a theorem 
of Beauville and Laszlo (\cite{BL}) asserts that the above map from the set
 of pairs
$(F_G,\beta)$ ``on $X$'' to that ``on $\D_x$'' is a bijection.

\medskip

Consider the stack $\Bun_G$. Let $_x\G$ denote the canonical 
$G(\O_x)$--torsor over $\Bun_G$, whose fiber over 
$\F_G\in\Bun_G$ is the set of isomorphisms $\F_G|_{\D_x}\simeq\F^0_G|_{\D_x}$.

Recall the stack $\H_G$ and let $_x\H_G$ denote its fiber over 
$\Bun_G\times x\subset \Bun_G\times X$; let $_x\Hb_G^\lambda\subset {}_x\H_G$ 
be the corresponding finite-dimensional substack.

The above mentioned result of \cite{BL} implies that both projections 
$\hl_G$ and $\hr_G$ realize $_x\H_G$
as a fibration over $\Bun_G$ with the typical fiber $\Gr_G$:
$$_x\H_G\simeq  \Gr_G\overset{G(\O_x)} \times{}_x\G.$$

We call these identifications $\on{id}^l$ and $\on{id}^r$, respectively. 
By definition, under $\on{id}^l$,
$_x\Hb_G^\lambda$ goes over to the substack
$\grb_G^\lambda\overset{G(\O_x)} \times {}_x\G \subset \Gr_G\overset{G(\O_x)} \times {}_x\G$.

We define the locally closed substack $_x\H_G^\lambda$ of $_x\Hb_G^\lambda$ as
$\Gr_G^\lambda\overset{G(\O_x)} \times {}_x\G$. In what follows, for $(\F_G,\F'_G,\beta)\in {}_x\H_G^\lambda$,
we will say that "$\F'_G$ is in position $\lambda$ with respect to $\F_G$".

\medskip

Thus, to $\S\in\Sph_G$ and an object $\T\in\Sh(\Bun_G)$ we can attach their twisted external products
$(\S\wt{\boxtimes}\T)^l$ and $(\S\wt{\boxtimes}\T)^r$. 
In particular, the intersection cohomology sheaf $\IC_{_x\Hb_G^\lambda}$ is nothing but
$(\A_G^\lambda\wt{\boxtimes}\IC_{\Bun_G})^l$.

\smallskip

We introduce the Hecke functors $_x\Hl_G(\cdot,\cdot)$ and $_x\Hr_G(\cdot,\cdot)$ from
$\Sph_G\times \Sh(\Bun_G)$ to $\Sh(\Bun_G)$ as follows:
$$_x\Hl_G(\S,\T):=\hl_G{}_{!}(\ast\S\wt{\boxtimes}\T)^r \text{ and }
_x\Hr_G(\S,\T):=\hr_G{}_{!}(\S\wt{\boxtimes}\T)^l.$$

Since every $\grb_G^\lambda$ is complete, the projections $\hl_G$ and $\hr_G$ are proper and 
the Hecke functors commute with Verdier duality. Moreover,
they are compatible with the tensor structure on $\Sph_G$ in the sense that
there exist functorial isomorphisms:
$$_x\Hl_G(\S_1,{}_x\Hl_G(\S_2,\T))\simeq {}_x\Hl_G(\S_1\star\S_2,\T) \text{ and }
_x\Hr_G(\S_1,{}_x\Hr_G(\S_2,\T))\simeq {}_x\Hr_G(\S_2\star\S_1,\T).$$

\smallskip

Finally, let us observe that there is a canonical isomorphism
$_x\Hl_G(\S,\T)\simeq {}_x\Hr_G(\ast\S,\T)$ and for a fixed $\S\in\Sph_G$, the functors
$$\T\to {}_x\Hl_G(\S,\T)\text{ and }\T\to {}_x\Hr_G({\mathbb D}(\S),\T)$$ 
are mutually (both left and right) adjoint.

\sssec{}  \label{semiinf}

In addition to the stratification by the $\Gr_G^\lambda$'s, the indscheme 
$\Gr_G$ is stratified by the so--called semiinfinite
orbits:

\smallskip

Consider the standard Borel structure on $\F^0_G$ given by the collection of the embeddings
$$\L^{\lambdach}_{\F^0_T}\hookrightarrow \V^{\lambdach}_{\F^0_G}$$
and for $\nu\in\Lambda$ let $\Sb_G^\nu$ denote the closed ind--subscheme of $\Gr_G$ corresponding to those pairs
$(\F_G,\beta)$, for which the meromorphic map
$$\L^{\lambdach}_{\F^0_T}(-\langle\nu,\lambdach\rangle\cdot x)\to 
\V^{\lambdach}_{\F^0_G}\overset{\beta}\longrightarrow \V^{\lambdach}_{\F_G}$$
does {\it not} have a pole $\forall\,\lambdach\in\check\Lambda_G^{+}$.

Inside $\Sb_G^\nu$ there is an open ind--subscheme $S_G^\nu$ corresponding to those pairs
$(\F_G,\beta)$, for which the above map
$$\L^{\lambdach}_{\F^0_T}(-\langle\nu,\lambdach\rangle\cdot x)\to \V^{\lambdach}_{\F_G}$$ 
does not have a zero either. 
Hence, $\Gr_G$ (resp., $\Sb_G^{\nu'}$) is stratified by the $S_G^\nu$, $\nu\in\Lambda_G$ (resp., by the
$S_G^\nu$, $\nu\leq\nu'$).

Consider the group--indschemes $U(\K_x)\subset U(\K_x)T(\O_x)\subset G(\K_x)$. For each $\nu\in\Lambda_G$, 
$S_G^\nu$ is stable with respect to the $U(\K_x)T(\O_x)$--action on $\Gr_G$. (Moreover, it is easy to see that
$S_G^\nu$ is isomorphic to the quotient of $U(\K_x)$ by a certain group-subscheme.)

\medskip

The following assertion (due to Mirkovic and Vilonen, \cite{MV}) will play a central role in this paper:

\begin{prop} \label{dimestimate}
Consider the intersection $\Gr_G^\lambda\cap S_G^\nu$.

\smallskip

\noindent{\em (a)} It is empty unless $\on{w}_0(\lambda)\leq \nu\leq \lambda$. For $\nu=\lambda$, it
is open inside $\Gr_G^\lambda$ (and hence is smooth and irreducible). For $\nu=\on{w}_0(\lambda)$, 
it is isomorphic to the point-scheme.

\smallskip

\noindent{\em (b)} $\dim(\Gr_G^\lambda\cap S_G^\nu)\leq \langle \lambda+\nu,\rhoch \rangle$.
\end{prop}

\sssec{}

It has been known since the works of Lusztig 
\cite{Lu}, Ginzburg \cite{Gi} and Mirkovic-Vilonen \cite{MV} that the functor 
$\star:\Sph_G\times\Sph_G\mapsto\Sph_G$ prolongs to a structure of tensor category on $\Sph_G$, and as such
it is equivalent to the category of finite-dimensional representations
of the Langlands--dual group $\check G$. We will need this statement in the following formulation (\cite{BD}):

\begin{thm} \label{dualgroup}
{\em (a)}
The convolution product $\S^1,\S^2\to \S^1\star\S^2$ admits associativity and commutativity constraints, 
i.e. $\Sph_G$ has a structure of a (symmetric) tensor category. Moreover, the functor 
$$\on{gRes}^G_T:\Sph_G\to \Lambda\text{--graded vector spaces}$$ given by
$$\on{gRes}^G_T(\S)=\underset{\nu\in\Lambda}
\oplus H_c^{\langle\nu,2\rhoch\rangle}(S_G^\nu,\S|_{S_G^\nu})\otimes\Ql(\langle\nu,\rhoch\rangle)$$
has a natural structure of a tensor functor. 

\smallskip

{\em (b)}
The algebriac group corresponding (by Tannaka's theorem) to 
$\Sph_G$ and the composition 
$$\text{Forget}\circ \on{gRes}^G_T: \Sph_G\to \Lambda\text{--graded vector spaces}\to \text{ vector spaces}$$
is canonically isomorphic to the Langlands dual group $\check G$, in such a way that the equivalence
$$F_G:\Sph_G\to \Rep(\check G)$$ has the following properties:

\begin{itemize}

\item
$\on{gRes}^G_T$ goes over to the natural restriction functor 
$\on{Res}^G_T:\Rep(\check G)\to \Rep(\check T)$, i.e.
$$\on{Res}^G_T\circ F_G\simeq \on{gRes}^G_T.$$

\item
$F_G(\A_G^\lambda)=V^{\lambda}\in \Rep(\check G)$.

\item
The countravariant self--functor $\S\to\DD(\ast\S)$ on $\Sph_G$ goes over to the dualization functor on 
$\Rep(\check G)$.

\end{itemize}

\end{thm}

The following assertion can be obtained by combining \thmref{dualgroup} and \propref{dimestimate}:

\begin{cor} \label{weightspaces}
For $\lambda\in\Lambda_G^{+}$ and $\nu\in\Lambda$ consider the $\Ql$-vector space \newline
$H^i_c(S_G^\nu,\A_G^\lambda|_{S_G^\nu})\otimes\Ql(\langle\nu,\rhoch\rangle)$. We have:

\smallskip

{\em (a)} It is zero for $i>\langle\nu,2\rhoch\rangle$.

\smallskip

{\em (b)} For $i=\langle\nu,2\rhoch\rangle$ it identifies with $V^\lambda(\nu)$. Moreover,
$$V^\lambda(\nu)\simeq 
H_c^{\langle\lambda+\nu,2\rhoch\rangle}(\Gr_G^\lambda\cap S_G^\nu,\Ql(\langle\lambda+\nu,\rhoch\rangle)).$$

\end{cor}

In particular, we infer that $V^\lambda(\nu)$ has a basis given by the irreducible components of the scheme
$\Gr_G^\lambda\cap S_G^\nu$ which have (the maximal possible) dimension $\langle\lambda+\nu,\rhoch\rangle$.

\medskip

\noindent{\it Remark.}
The paper \cite{MV} contains an even stronger assertion, namely that the cohomology group
$H^i(S_G^\nu,\A_G^\lambda|_{S_G^\nu})$ vanishes unless $i=\langle\nu,2\rhoch\rangle$ and that the scheme
$\Gr_G^\lambda\cap S_G^\nu$ is of pure dimension $\langle\nu,2\rhoch\rangle$. However, for our purposes 
the formulation of \corref{weightspaces} will be sufficient.

\ssec{Proof of \thmref{mainHeckeprince}}

\sssec{}

For a point $x\in X$, let $_x\ol{Z}$ be the pre-image $'\p^{-1}({}_x\Hb_G^\lambda)\subset \ol{Z}$. 
(By abuse of notation we will continue to denote by $\phi$ (resp., $'\hr_G$, $\hr_G$, $\hl_G$, $i_{\nu}$) 
the corresponding maps when $\BunBb\times X$ is replaced by $\BunBb=\BunBb\times x$.)

\smallskip

We have the following diagram:

$$
\CD
\BunBb @<\phi<< _x\ol{Z} @>{'\hr_G}>>  \BunBb  \\
  @V{\p}VV        @V{'\p}VV     @V{\p}VV   \\
\Bun_G   @<{\hl_G}<<   _x\Hb_G^\lambda    @>{\hr_G}>>  \Bun_G
\endCD
$$

\smallskip

To simplify the notation, we will prove the following version of \thmref{mainHeckeprince}:

\begin{thm} \label{mainHeckeprinceatx}
$\phi_{!}(\IC_{_x\ol{Z}})\simeq \underset{\nu\in\Lambda_G^{\on{pos}}}\oplus{} i_\nu{}_{!}
(\IC_{\BunBb})\otimes V^{\lambda}(\on{w}_0(\lambda)+\nu).$
\end{thm}

The proof of \thmref{mainHeckeprince} is no diffrent from the presented below proof of
\thmref{mainHeckeprinceatx}.

\medskip

For $\nu\in\Lambda_G^{\on{pos}}$, let $_{x,\nu}\BunBb$ (resp., $_{x,\geq\nu}\BunBb$)
be the closed (resp., locally closed) substack of $\BunBb$ that corresponds to those 
triples $(\F_G,\F_T,\kappa)$ for which each of the maps
$$\L^{\lambdach}_{\F_T}\to \V^{\lambdach}_{\F_G}$$
has a zero of order not less than (resp., exactly) $\langle\nu,\lambdach\rangle$ at $x$, for 
all $\lambdach\in\check\Lambda_G^{+}$. 

\smallskip

By definition, $_{x,\geq 0}\BunBb=\BunBb$ and 
$\BunBb(\Fqb)=\underset{\nu\in\Lambda_G^{\on{pos}}}\cup {}_{x,\nu}\BunBb(\Fqb)$. 
We will denote by $j_{\nu}$ the 
embedding of $_{x,\nu}\BunBb$ into $\BunBb$. It is easy to see that the map $i_{\nu}$ defines isomorphisms 
$$\BunBb\to{}_{x,\geq\nu}\BunBb\text{ and } _{x,0}\BunBb\to {}_{x,\nu}\BunBb.$$

\smallskip

Let us consider in addition the open substack 
$_{x,\nu}\Bun_B\subset {}_{x,\nu}\BunBb$ that corresponds to triples 
$(\F_G,\F_T,\kappa)$ as above for which each of the maps 
$$\L^{\lambdach}_{\F_T}\to \V^{\lambdach}_{\F_G}$$ 
has a zero of order exactly $\langle\nu,\lambdach\rangle$ at $x$ and has no zeroes at other points of $X$. We have:
$_{x,0}\Bun_B=\Bun_B$ and for $\nu\neq 0$, $i_{\nu}$ defines an isomorphism 
$\Bun_B\to{}_{x,\nu}\Bun_B$.

\medskip

By applying the decomposition theorem \cite{BBD} 
(and taking into account the fact that $\phi_{!}(\IC_{_x\ol{Z}})$ is Verdier self-dual),
we infer that \thmref{mainHeckeprinceatx} is equivalent to the following assertion:

\begin{prop} \label{propmainHeckeprinceatx} The complex $j_{\nu}^*\circ \phi_{!}(\IC_{_x\ol{Z}})$
satisfies the following:

\smallskip

{\em (a)}
It lives in non-positive cohomological degrees (in the sense of the perverse t-structure). 

\smallskip

{\em (b)}
$h^0(j_{\nu}^*\circ \phi_{!}(\IC_{_x\ol{Z}}))$ identifies with $\IC_{_{x,\nu}\BunBb}\otimes 
V^{\lambda}(\w_0(\lambda)+\nu)$.
\end{prop}

\sssec{}  \label{describeprince}

For an element $\nu\in\Lambda_G^{\on{pos}}$, let $Z^{?,\nu}$ (resp., $Z^{\nu,?}$) 
denote the pre-image in $_x\ol{Z}$
of the stratum $_{x,\nu}\BunBb\subset\BunBb$ under the map $'\hr_G$ (resp., $\phi$). For two elements 
$\nu,\nu'\in \Lambda_G^{\on{pos}}$ let $Z^{\nu,\nu'}$ denote the intersection 
$Z^{?,\nu'}\cap Z^{\nu,?}$. 
(As we shall see shortly, $Z^{\nu,\nu'}$ is empty unless $\nu\geq\nu'$.)

\smallskip

If now $\lambda'$ is an element of $\Lambda_G^{+}$ such that $\lambda'\leq\lambda$, 
let $Z^{\nu,\nu',\lambda'}$ denote the intersection
$$Z^{\nu,\nu'}\cap {}'\p^{-1}({}_x\H_G^{\lambda'}).$$

\smallskip

To prove \propref{propmainHeckeprinceatx} it is enough to prove the following 

\begin{prop} \label{ppropmainHeckeprinceatx}
For $\nu,\nu'$ and $\lambda'$ as above, let $K^{\nu,\nu',\lambda'}\in\Sh({}_{x,\nu}\BunBb)$ be defined as
$\phi_{!}(\IC_{_x\ol{Z}}|_{Z^{\nu,\nu',\lambda'}})$. We have:

\smallskip

{\em (a)}
$K^{\nu,\nu',\lambda'}$ 
lives in cohomological degrees $\leq 0$ and the equality is strict unless $\nu'=0$ and $\lambda'=\lambda$.

\smallskip

{\em (b)}
The $*$-restriction of $K^{\nu,0,\lambda}$ to 
$_{x,\nu}\BunBb-{}_{x,\nu}\Bun_B$ lives in cohomological degrees $<0$.

\smallskip

{\em (c)}
The restriction of $h^0(K^{\nu,0,\lambda})$ to $_{x,\nu}\Bun_B$ can be identified with 
$$\IC_{_{x,\nu}\Bun_B}\otimes V^{\lambda}(\on{w}_0(\lambda)+\nu).$$ 

\end{prop}

Consider the open sub-stack $_{x,0}\BunBb\subset\BunBb$.
For $(\F_G,\F_T,\kappa)\in {}_{x,0}\BunBb$, the restriction of this data to the
formal disc $\D_x$ defines a $B$-bundle over $\D_x$, or which is the same, a $B(\O_x)$-torsor. 
Therefore, there exists
a canonical $B(\O_x)$-torsor over $_{x,0}\BunBb$, which we will denote by $_x\B$. 

\smallskip

Using the fact that $i_{\nu}$ maps $_{x,0}\BunBb$ isomorphically onto $_{x,\nu}\BunBb$, 
we obtain a $B(\O_x)$-torsor $_x\B^\nu$ over each $_{x,\nu}\BunBb$. 
\footnote{Note, however, that there is no globally defined 
$B(\O_x)$-torsor over $\BunBb$ that would restrict to $_x\B^\nu$ over the corresponding stratum.}
We will denote by $_x\overset{o}{\B}{}^\nu$ 
the restriction of $_x\B^\nu$ to the open sub-set $_{x,\nu}\Bun_B\subset{}_{x,\nu}\BunBb$.
 
\medskip

The next assertion follows from the definitions:

\begin{lem} \label{descrprince}

{\em (a)} The stack $Z^{?,\nu'}$ with the projection $'\hr_G:Z^{?,\nu'}\to {}_{x,\nu'}\BunBb$ is a 
locally trivial fibration with the typical fiber $\grb_G^{-\on{w}_0(\lambda)}$. More precisely there is an isomorphism
$$Z^{?,\nu'}\simeq \grb_G^{-\on{w}_0(\lambda)} \overset{B(\O_x)}\times {}_x\B^{\nu'},$$
where $\grb_G^{-\on{w}_0(\lambda)}$ is viewed as a $B(\O_x)$--scheme via $B(\O_x)\hookrightarrow G(\O_x)$.

\smallskip

{\em (b)} The sub-stack $Z^{\nu,\nu',\lambda'}\hookrightarrow  Z^{?,\nu'}$ identifies 
(using the notation of point (a)) with
$$(\Gr_G^{-\on{w}_0(\lambda')}\cap S_G^{-\on{w}_0(\lambda)-\nu+\nu'})\overset{B(\O_x)}\times {}_x\B^{\nu'}
\subset  \grb_G^{-\on{w}_0(\lambda)} \overset{B(\O_x)}\times {}_x\B^{\nu'}\simeq Z^{?,\nu'}.$$

\smallskip

{\em (c)} The sub-stack $Z^{\nu',\nu,\lambda'}$ when viewed as a stack projecting to 
$_{x,\nu}\BunBb$ by means of $\phi$, identifies with
$$(\Gr_G^{\lambda'}\cap S_G^{\nu-\nu'+\on{w}_0(\lambda)})\overset{B(\O_x)}\times {}_x\B^{\nu}.$$

\smallskip

{\em (d)} For every $\nu,\nu',\lambda'$, the pre--image of the open substack 
$_{x,\nu}\Bun_B\subset {}_{x,\nu}\BunBb$ in $Z^{\nu',\nu,\lambda'}$ under the projection $\phi$ 
coincides with the pre-image
of $_{x,\nu'}\Bun_B\subset {}_{x,\nu'}\BunBb$ under $'\hr_G$.

\end{lem}

\sssec{}  \label{proofmainHeckeprince}

We are finally ready to prove \propref{ppropmainHeckeprinceatx} (and hence \thmref{mainHeckeprince}):

\smallskip

Since the map $\hr_G:{}_x\Hb^\lambda_G\to\Bun_G$ is a locally trivial fibration, the sheaf
$\IC_{_x\ol{Z}}|_{Z^{?,\nu'}}$ is a twisted external product 
$$\A_G^{-\w_0(\lambda)} \wt{\boxtimes} \IC_{\BunBb}|_{_{x,\nu'}\BunBb}$$
(in the above identification we have used the isomorphism of \lemref{descrprince}(a)).

\smallskip

Hence, the restriction $\IC_{_x\ol{Z}}|_{Z^{\nu,\nu',\lambda'}}$ is a twisted external product
$$\A_G^{-\w_0(\lambda)}|_{\Gr_G^{-\w_0(\lambda')}\cap
S_G^{-\w_0(\lambda)-\nu+\nu'}}\widetilde{\boxtimes}\IC_{\BunBb}|_{_{x,\nu'}\BunBb}.$$

\medskip

First of all, by the very definition of $\IC$, 
the sheaf $\IC_{\BunBb}|_{_{x,\nu'}\BunBb}$ lives in cohomological degrees $\leq 0$ 
and the equality is strict unless $\nu'=0$. 

\smallskip

The restriction $\A_G^{-\w_0(\lambda)}|_{\Gr_G^{-\w_0(\lambda')}}$ is a constant complex over 
$\Gr_G^{-\w_0(\lambda')}$ which lives in cohomological degrees $<0$ unless $\lambda'=\lambda$. Using 
\propref{dimestimate}, we obtain that the sheaf
$$\A_G^{-\w_0(\lambda)}|_{\Gr_G^{-\w_0(\lambda')}\cap S_G^{-\w_0(\lambda)-\nu+\nu'}}$$ 
lives in cohomological degrees
$$\leq-\on{codim}(\Gr_G^{-\w_0(\lambda')}\cap S_G^{-\w_0(\lambda')-\nu+\nu'},\Gr_G^{-\w_0(\lambda')})\leq
-\langle \lambda'-\lambda+\nu-\nu',\rhoch\rangle,$$
with the inequality being strict unless $\lambda'=\lambda$.

\medskip

Hence, if either $\lambda'\neq\lambda$ or $\nu'\neq 0$, the sheaf $\IC_{_x\ol{Z}}|_{Z^{\nu,\nu',\lambda'}}$
lives in cohomological degrees $<-\langle \lambda'-\lambda+\nu-\nu',\rhoch\rangle$.
However, \lemref{descrprince}(c) and \propref{dimestimate} imply that the dimension of the fibers of the projection
$$\phi:Z^{\nu,\nu',\lambda'}\to{}_{x,\nu}\BunBb$$ is $\leq \langle \lambda'-\lambda+\nu-\nu',\rhoch\rangle$
and this proves point (a) of \propref{ppropmainHeckeprinceatx}.

\medskip

Let $\overset{o}Z{}^{\nu,0,\lambda}$ denote the preimage of the open subset $\Bun_B\subset {}_{x,0}\BunBb$
in $Z^{\nu,0,\lambda}$ under the map $'\hr_G$. Considerations analogous to the ones above
show that the restriction of 
$\IC_{_x\ol{Z}}|_{Z^{\nu,0,\lambda}}$ to $Z^{\nu,0,\lambda}-\overset{o}Z{}^{\nu,0,\lambda}$
is concentrated in cohomological degrees $<0$. This implies \propref{ppropmainHeckeprinceatx}(b) in view of 
\lemref{descrprince}(d).

\medskip

To prove \propref{ppropmainHeckeprinceatx}(c), observe that the restriction of 
$\IC_{_x\ol{Z}}|_{\overset{o}Z{}^{\nu,0,\lambda}}$ is a constant sheaf. More precisely, when
we identify $\overset{o}Z{}^{\nu,0,\lambda}$ (using points (c) and (d) of \lemref{descrprince}) with
$(\Gr_G^\lambda\cap S_G^{\w_0(\lambda)+\nu})\overset{B(\O_x)}\times{}_x\overset{o}{\B}{}^\nu $,
the sheaf $\IC_{_x\ol{Z}}|_{\overset{o}Z{}^{\nu,0,\lambda}}$ identifies with the twisted external product
$(\Ql[1](\frac{1}{2}))^{\otimes \langle \nu,2\rhoch\rangle}_{\Gr_G^\lambda\cap S_G^{\w_0(\lambda)+\nu}}\wt\boxtimes
\IC_{_{x,\nu}\Bun}$.

\smallskip

As the group $B(\O_x)$ is connected, this implies that the $0$-th=top perverse cohomology of its direct
image on $_{x,\nu}\Bun$ can be identified with
$$\IC_{_{x,\nu}\Bun}\otimes H^0_c(\Gr_G^\lambda\cap S_G^{\w_0(\lambda)+\nu},
(\Ql[1](\frac{1}{2}))^{\otimes \langle \nu,2\rhoch\rangle}).$$

Now, the assertion of \propref{ppropmainHeckeprinceatx}(c) follows from \corref{weightspaces}(b).


\section{Action of Hecke operators (the non-principal case)} 
\label{actnonprince}

\ssec{Variant of the basic diagram}

\sssec{}  \label{infinity}

Let $P$ be a parabolic sub-group of $G$. 
Consider the stack $_{x,\infty}\BunPbw$ that classifies triples
$$(\F_G,\F_M,\widetilde\kappa_P:\F_G|_{X-x}\to 
\GUPb\underset{M}\times \F_M|_{X-x}).$$ 

In other words, a point of $_{x,\infty}\BunPbw$ is a pair 
$(\F_G,\F_M)\in\Bun_G\times\Bun_M$ and a 
compatible system of maps
$$\widetilde\kappa_P^{\V}:
(\V^{U(P)})_{\F_M}\to \V_{\F_G}(\infty\cdot x),$$
for every $G$-module $\V$, which satisfy the Pl\"ucker 
relations in the same sense as in the definition of $\BunPbw$ 
(cf. \secref{BunPbw}).

By bounding degrees of poles of the maps 
$\widetilde\kappa_P^{\V}$ we can represent
$_{x,\infty}\BunPbw$ as a union of its finite--dimensional closed substacks:

\smallskip

Let $\nu\in\Lambda$ be dominant with respect to $M$, i.e. 
$\nu\in\Lambda_M^+$. We define a closed
finite--dimensional substack $_{x,\geq\nu}\BunPbw$ as follows: 
it corresponds to those triples
$(\F_G,\F_M,\widetilde\kappa_P)$ for which 
$\widetilde\kappa_P^{\V}$ maps
$$(\V^{U(P)})_{\F_M}\to \V_{\F_G}
(\langle -\w_0^M(\nu),\lambdach\rangle\cdot x),$$
if $\V$ is a $G$-module whose weights are $\leq \lambdach$.
Let $i_{\nu}$ denote the corresponding closed embedding.

In particular, $\BunPbw\subset{}_{x,\infty}\BunPbw$ coincides with 
$_{x,\geq 0}\BunPbw$. It is clear that if for $\nu'\in\Lambda^+_M$, 
$\w_0^M(\nu'-\nu)\in\Lambda_G^{\on{pos}}$, 
then $_{x,\geq\nu'}\BunPbw$ is contained in 
$_{x,\geq\nu}\BunPbw$ and that $_{x,\infty}\BunPbw$ is 
the inductive limit of the $_{x,\geq\nu}\BunPbw$'s. 

\smallskip

We will abuse the notation and continue to denote by $\pw_P$ and $\qw_P$  
the natural projections from $_{x,\infty}\BunPbw$ to $\Bun_G$ and to $\Bun_M$, respectively. 
Finally, note that when $P=G$, the stack $_{x,\infty}\BunPbw$ is identified with the Hecke 
stack $_x\H_G$.

\sssec{}  \label{variantM}

Let $_{x,\infty}Z_{P,M}$ denote the Cartesian product:
$$_{x,\infty}Z_{P,M}={}_x\H_M{\underset{\Bun_M}\times}{}_{x,\infty}\BunPbw,$$
where $_x\H_M$ maps to $\Bun_M$ by means of $\hr_M$.

By definition, the stack $_{x,\infty}Z_{P,M}$ classifies quadruples $(\F_G,\F_M,\F'_M,\wt\kappa_P)$, where \newline
$(\F_G,\F'_M,\wt\kappa_P)$ is a point of $_{x,\infty}\BunPbw$ and $\F_M$ is an $M$-bundle over $X$ identified with 
$\F'_M$ over $X-x$. However,
to a quadruple $(\F_G,\F_M,\F'_M,\wt\kappa_P)$ as above one can attach another point of $_{x,\infty}\BunPbw$, namely
$(\F_G,\F_M,\wt\kappa_P)$. Thus, as in \secref{bd} we obtain a second projection 
$'\hl_M:{}_{x,\infty}Z_{P,M}\to\Bun_M$.

\smallskip

It is easy to see that we obtain in fact two Cartesian squares:
$$
\CD
_{x,\infty}\BunPbw  @<{'\hl_M}<<  _{x,\infty}Z_{P,M} @>{'\hr_M}>>  _{x,\infty}\BunPbw  \\
@V{\qw_P}VV     @V{'\qw_P}VV    @V{\qw_P}VV    \\
\Bun_M  @<{\hl_M}<<  _x\H_M  @>{\hr_M}>>   \Bun_M
\endCD
$$

\smallskip

Therefore, for $\S\in \Sph_M$ and $\T\in \Sh({}_{x,\infty}\BunPbw)$ 
we can define their
twisted external products $(\S\tboxtimes\T)^l,\,(\S\tboxtimes\T)^r\in 
\Sh({}_{x,\infty}\BunPbw)$.
We introduce the Hecke functors $_x\Hl_{P,M}(\cdot,\cdot)$ and 
$_x\Hr_{P,M}(\cdot,\cdot)$
from $\Sph_M\times \Sh({}_{x,\infty}\BunPbw)$ to $\Sh({}_{x,\infty}\BunPbw)$ 
by setting
$$_x\Hl_{P,M}(\S,\T)={}'\hl_M{}_{!}(\ast\S\tboxtimes\T)^r \text{ and }
_x\Hr_{P,M}(\S,\T)={}'\hr_M{}_{!}(\S\tboxtimes\T)^l.$$

These functors have the properties exactly analogous to those of 
$_x\Hl_G(\cdot,\cdot)$ and 
$_x\Hr_G(\cdot,\cdot)$, described in \secref{globalgrass}.

For $\nu\in\Lambda_M^{+}$, we will denote by $_xH^\nu_{P,M}$ 
the self-functor on $\Sh({}_{x,\infty}\BunPbw)$ defined as
$\T\to {}_x\Hl_{P,M}(\A^\nu_M,\T)$. 

We have the following assertion:

\begin{thm}    \label{mainHeckenonprinceM}
$_xH_{P,M}^\nu(\IC_{_{x,\geq 0}\BunPbw})$ is canonically isomorphic to the intersection
cohomology sheaf of $_{x,\geq -\w_0^M(\nu)}\BunPbw$.
\end{thm}

\smallskip

A proof of this theorem will be given in the next subsection.

\sssec{}  \label{variantG}

Now let us consider the stack
$$_{x,\infty}Z_{P,G}={}_x\H_G{\underset{\Bun_G}\times}{}_{x,\infty}\BunPbw,$$
where $_x\H_G$ projects to $\Bun_G$ by means of $\hr_G$.

Similarly to what we had above, one can define a second projection 
$'\hl_G:{}_{x,\infty}Z_{P,G}\to\Bun_G$ and 
we obtain a diagram

$$
\CD
_{x,\infty}\BunPbw  @<{'\hl_G}<<  _{x,\infty}Z_{P,G} @>{'\hr_G}>>  
_{x,\infty}\BunPbw  \\
@V{\pw_P}VV     @V{'\pw_P}VV    @V{\pw_P}VV    \\
\Bun_G  @<{\hl_G}<<  _x\H_G  @>{\hr_G}>>   \Bun_G
\endCD 
$$
in which both squares are Cartesian.

\smallskip

Proceeding as above, we obtain two functors $_x\Hl_{P,G}(\S,\T)$ and 
$_x\Hr_{P,G}(\S,\T)$
from $\Sph_G\times \Sh({}_{x,\infty}\BunPbw)\mapsto 
\Sh({}_{x,\infty}\BunPbw)$. For $\lambda\in\Lambda_G^+$
and $\T\in \Sh({}_{x,\infty}\BunPbw)$ we set 
$_xH_{P,G}^\lambda(\T)={}_x\Hl_{P,G}(\A^\lambda_G,\T)$.

It is easy to see that for $?=\to$  or  $\leftarrow$  and  $?'=\to$  
or  $\leftarrow$
the actions $_xH_{P,G}^?$ and $_xH_{P,M}^{?'}$ 
commute in the sense that there is a functorial isomorphism 
$$_xH_{P,G}^?(\S_G,{}_xH_{P,M}^{?'}(\S_M,\T))\simeq 
{}_xH_{P,M}^?(\S_M,{}_xH_{P,G}^{?'}(\S_G,\T)),$$
where $\S_G\in\Sph_G$, $\S_M\in\Sph_M$ and $\T\in \Sh({}_{x,\infty}\BunPbw)$.

\smallskip

\begin{thm}  \label{mainHeckenonprinceG}

For $\lambda\in\Lambda_G^{+}$ there is a canonical isomorphism:
$$_xH_{P,G}^\lambda(\IC_{_{x,\geq 0}\BunPbw})\simeq \underset{\nu\in\Lambda_M^{+}}\oplus 
\IC_{_{x,\geq \nu}\BunPbw}\otimes \on{Hom}_{\check M}(U^{\nu},V^{\lambda}).$$

\end{thm}

\sssec{}  

The combination of \thmref{mainHeckenonprinceG} and \thmref{mainHeckenonprinceM} yields the following:

\smallskip

Let $\on{gRes}^G_M:\Sph_G\to\Sph_M$ denote the tensor functor corresponding to the restriction
functor $\on{Res}^G_M:\on{Rep}(\check G)\to \on{Rep}(\check M)$ under the equivalence of 
\thmref{dualgroup}.

\begin{cor}  \label{compwithtensorstructure}
The two functors $\Sph_G\to \Sh({}_{x,\infty}\BunPbw)$:
$$\S\mapsto {}_x\Hl_{P,G}(\S,\IC_{_{x,\geq 0}\BunPbw})\text{  and  } 
\S\mapsto {}_x\Hr_{P,M}(\on{gRes}^G_M(\S),\IC_{_{x,\geq 0}\BunPbw})$$ 
are canonically isomorphic. 
\end{cor}

Moreover, from the proof of the above theorems, it will follow that
the isomorphism of functors of the above corollary is compatible with the tensor 
structure in the following sense:

\smallskip

Let $\S_1,\S_2\in \Sph_G$ the two isomorphisms
\begin{align*}
&_x\Hl_{P,G}(\S_1,{}_x\Hl_{P,G}(\S_2,\IC_{\BunPbw}))
\simeq {}_x\Hl_{P,G}(\S_1\star\S_2,\IC_{\BunPbw})\to \\
&{}_x\Hr_{P,M}(\on{gRes}^G_M(\S_1\star\S_2),\IC_{\BunPbw}) \to
{}_x\Hr_{P,M}(\on{gRes}^G_M(\S_1)\star \on{gRes}^G_M(\S_2),\IC_{\BunPbw}) \to \\
&{}_x\Hr_{P,M}(\on{gRes}^G_M(\S_2),{}_x\Hr_{P,M}(\on{gRes}^G_M(\S_1),\IC_{\BunPbw})) 
\end{align*}
and
\begin{align*}
&{}_x\Hl_{P,G}(\S_1,{}_x\Hl_{P,G}(\S_2,\IC_{\BunPbw}))\to 
{}_x\Hl_{P,G}(\S_1,{}_x\Hr_{P,M}(\on{gRes}^G_M(\S_2),\IC_{\BunPbw})) \to \\
&{}_x\Hr_{P,M}(\on{gRes}^G_M(\S_2),{}_x\Hl_{P,G}(\S_1,\IC_{\BunPbw})) \to \\
&{}_x\Hr_{P,M}(\on{gRes}^G_M(\S_2),{}_x\Hr_{P,M}(\on{gRes}^G_M(\S_1),\IC_{\BunPbw}))
\end{align*}
coincide. We leave the verification to the reader.

\sssec{}   \label{prcomheckenonprince}

Now let us explain how \thmref{mainHeckenonprinceG} implies \thmref{commutewithHeckenonprince}. 

\smallskip

To simplify the notation, we shall prove the following statement:

\smallskip

\noindent {\it For $\S\in\Sph_G$ and $\T\in\Sh(\Bun_M)$ there is a functorial isomorphism 
\begin{equation} \label{hnp}
_x\Hl_G(\S,\Eis_M^G(\T))\simeq \Eis^G_M({}_x\Hl_M(\on{gRes}^G_M(\S),\T)).
\end{equation}
Moreover, this isomorphism is compatible with the tensor structures on 
$\Sph_G$ and on $\Sph_M$.}

\bigskip

The proof of \thmref{commutewithHeckenonprince} in its original form (i.e. when $x$ is not fixed) is 
absolutely analogous.

\medskip

For $\T\in\Sh(\Bun_M)$ the LHS of \eqref{hnp} can be rewritten using the stack $_{x,\infty}\BunPbw$ 
instead of $\BunPbw$: 
$$_x\Hl_G(\S,\Eis_M^G(\T)):={}_x\Hl_G(\S,\pw_P{}_{!}(\qw_P^*(\T)\otimes \IC_{_{x,\geq 0}\BunPbw}))\otimes 
(\Ql[1](\frac{1}{2}))^{\otimes -\dim(\Bun_M)}.$$

\smallskip

The maps $\qw_P\circ{}'\hl_G$ and $\qw_P\circ{}'\hr_G$ from $_{x,\infty}Z_{P,G}$ to $\Bun_M$ coincide.
Therefore, using the projection formula and the base change on the diagram definining 
$_{x,\infty}Z_{P,G}$, the above expression can be rewritten as
\begin{equation} \label{manone}
\pw_P{}_{!}(\qw_P^*(\T)\otimes {}_x\Hl_{P,G}(\S,\IC_{_{x,\geq 0}\BunPbw}))\otimes 
(\Ql[1](\frac{1}{2}))^{\otimes -\dim(\Bun_M)}.
\end{equation}

\smallskip

Now, let us apply \thmref{mainHeckenonprinceG} (in the form of \corref{compwithtensorstructure})
and rewrite \eqref{manone} as
\begin{equation*} \label{mantwo}
\pw_P{}_{!}(\qw_P^*(\T)\otimes {}_x\Hr_{P,M}(\on{gRes}^G_M(\S),\IC_{_{x,\geq 0}\BunPbw}))
\otimes (\Ql[1](\frac{1}{2}))^{\otimes -\dim(\Bun_M)}.
\end{equation*}

\smallskip

Again, since $\pw_P\circ{}'\hl_M=\pw_P\circ{}'\hr_M$ as maps $_{x,\infty}Z_{P,M}\to\Bun_G$,
the projection formula implies that \eqref{mantwo} is the same as 
$$\pw_P{}_{!}(\qw_P^*({}_x\Hl_M(\on{gRes}^G_M(\S),\T))\otimes \IC_{_{x,\geq 0}\BunPbw}) 
\otimes (\Ql[1](\frac{1}{2}))^{\otimes -\dim(\Bun_M)},$$
which is what we had to prove.

The fact that this system of isomorphisms is compatible with the tensor structure follows
from corresponding property of the isomorphism of \corref{compwithtensorstructure}.

\ssec{Proof of \thmref{mainHeckenonprinceM}} 

\sssec{}  \label{stratinf}

Let $_{x,0}\BunPbw$ denote the open substack of $_{x,\geq 0}\BunPbw=\BunPbw$ that corresponds to 
those triples $(\F_G,\F_M,\wt\kappa_P)$ for which the maps
$$\widetilde\kappa_P^{\V}:(\V^{U(P)})_{\F_M}\to \V_{\F_G}$$
have no zero at $x$.  

\begin{prop}
The composition
$$'\hl_M{}^{-1}(_{x,0}\BunPbw)\hookrightarrow {}_{x,\infty}Z_{P,M}
\overset{'\hr_M}\longrightarrow {}_{x,\infty}\BunPbw$$
is a locally closed embedding.
\end{prop}

\begin{proof}

It is easy to see that the assertion of the proposition amounts to the following:

\medskip

Let $(\F_G,\F_M,\wt\kap_P)$ be as $S$-point of $\BunPbw$, such that for every $G$-module $\V$, the sheaf embedding
$$\wt\kap^{\V}:(\V^{U(P)})_{\F_M}\to \V_{\F_G}$$
has no zero along the divisor $x\times S\subset X\times S$. Suppose that $(\F_M,\F'_M,\beta_M)$ is
an $S$-point of $_x\H_M$ such that the {\it a priori} meromorphic maps
$$\wt\kap'{}^{\V}:(\V^{U(P)})_{\F'_M}|_{(X-x)\times S}\to \V_{\F_G}|_{(X-x)\times S}$$
extend to regular maps on $X\times S$, which do not have zeroes along $x\times S$.

In this case we have to show that $\F'_M$ in fact equals $\F_M$, i.e. that $\beta_M$ induces regular
maps $\beta_M^{\U}:\U_{\F'_M}\to \U_{\F_M}$ for every $M$-module $\U$.

\medskip

First, let us assume that $\U$ is isomorphic to $\V^{U(P)}$ for some $G$-module $\V$. In this case,
the above map $\wt\kap'{}^{\V}$ must factor as
$$(\V^{U(P)})_{\F'_M}\overset{\beta_M^{\U}} \longrightarrow
(\V^{U(P)})_{\F_M}\overset{\wt\kap^{\V}}\longrightarrow \V_{\F_G},$$
because $(\V^{U(P)})_{\F_M}$ is maximal at $x\times S$ in $\V_{\F_G}$, by assumption.
In particular, this means that $\beta_M^{\U}$ is regular.

\medskip

The same argument shows that if $\U$ is $1$-dimensional, corresponding to a character 
$\lambdach\in\check\Lambda_{G,P}\cap \Lambda^+_G$, the map $\beta_M^{\U}$ is an isomorphism. Since
the semigroup $\lambdach\in\check\Lambda_{G,P}$ generates $\check\Lambda_{G,P}$, we obtain that
$\beta_M^{\U}$ is an isomorphism for all $1$-dimensional representations.

To finish the proof, it remains to observe that every $M$-module can be embedded into a tensor
product of a $1$-dimensional module and a one of the form $\V^{U(P)}$.

\end{proof}

For $\nu\in\Lambda_M^+$, let $_{x,\nu}\BunPbw$ denote the image under 
$'\hl_M{}^{-1}(_{x,0}\BunPbw)\overset{'\hr_M}\longrightarrow {}_{x,\infty}\BunPbw$
of $'\hl_M{}^{-1}(_{x,0}\BunPbw)\cap {}'\qw_P^{-1}({}_x\H_M^\nu)$.
Let $j_\nu$ denote the locally closed embedding of $_{x,\nu}\BunPbw$ into $_{x,\infty}\BunPbw$.

\begin{prop}
The locally closed substacks $_{x,\nu'}\BunPbw$, $\on{w}_0^M(\nu'-\nu)\in\Lambda^{\on{pos}}_G$
form a stratification of $_{x,\geq\nu}\BunPbw$, with $_{x,\nu}\BunPbw$ being the biggest stratum.
\end{prop}

\begin{proof}

The fact that $_{x,\nu}\BunPbw$ is contained in $_{x,\geq\nu}\BunPbw$ follows from the definitions.
Therefore, to prove the proposition we have to show that every $\Fqb$-point 
$(\F_G,\F_M,\wt\kap_P)$ of $_{x,\infty}\BunPbw$ is contained is some $_{x,\geq\nu}\BunPbw$.

The latter amounts to showing that there exists an $M$-bundle $\F'_M$ and an isomorphism
$\beta_M:\F'_M|_{X-x}\to \F_M|_{X-x}$, such that the composition
$$(\V^{U(P)})_{\F'_M}|_{X-x}\overset{\beta_M^{\U}} \longrightarrow
(\V^{U(P)})_{\F_M}|_{X-x}\overset{\wt\kap^{\V}}\longrightarrow \V_{\F_G}|_{X-x}$$
extends to a regular map on the entire $X$, which has no zeroes at $x$. 

The construction of such $\F'_M$ repeats the proof of \propref{pointsofBunBb} in its parabolic variant.

\end{proof}

\medskip

Since the complex $_xH_{P,M}^\nu(\IC_{\BunPbw})$ is self-dual, in order to prove
\thmref{mainHeckenonprinceM} it suffices to prove the following assertion:
 
\begin{prop} \label{propmainHeckenonprinceatx} For any $\nu'\in\Lambda_M^{+}$ the complex 
$$j_{\nu'}^*\circ{}'\hl_M{}_{!}(\A^{-\on{w}_0^M(\nu)}_M\tboxtimes \IC_{\BunPbw})^r$$
lives in cohomological degrees $\leq 0$ and the inequality is strict unless 
$\nu'=-\on{w}_0^M(\nu)$. Moreover, for $\nu'=-\on{w}_0^M(\nu)$
its $0$-th cohomology can be identified with $\IC_{_{x,-\on{w}_0^M(\nu)}\BunPbw}$
\end{prop}

\sssec{}

For $\nu',\eta,\eta'\in\Lambda_M^{+}$ let us consider the following locally closed sub-stacks of 
$_{x,\infty}Z_{P,M}$:

\begin{align*}
&Z_{P,M}^{\nu',?,?}:=({}'\hl_M)^{-1}({}_{x,\nu'}\BunPbw),\,\,  
Z_{P,M}^{?,\eta,?}:=({}'\hr_M)^{-1}({}_{x,\eta}\BunPbw),\\
&Z_{P,M}^{?,?,\eta'}:=({}'\qw_P)^{-1}({}_x\H^{\eta'}_M).
\end{align*}
Let $Z_{P,M}^{\nu',\eta,\eta'}$ denote their intersection and similarly for $Z_{P,M}^{?,\eta,\eta'}$, etc.

We will denote by $K_M^{\nu',\eta,\eta'}$ the complex 
$$'\hl_M{}_{!}((\A^{-\w_0^M(\nu)}_M\tboxtimes \IC_{\BunPbw})^r|_{Z_{P,M}^{\nu',\eta,\eta'}})$$
on $_{x,\nu'}\BunPbw$.

To prove the proposition it is enough to show that $K_M^{\nu',\eta,\eta'}$ lives in negative cohomological 
degrees unless $\nu'=-\w_0^M(\nu)$, $\eta=0$, $\eta'=\nu$ and that 
$$h^0(K_M^{-\w_0^M(\nu),0,\nu})\simeq \IC_{_{x,-\w_0^M(\nu)}\BunPbw}.$$

Recall the indscheme $\Conv_M$ (cf. \secref{convolution}). As we have seen before, 
the map $pr'$ makes it a fibration
over $\Gr_M$ with the typical fiber $\Gr_M$. For $\nu_1,\nu_2\in\Lambda^+_M$ let $\Conv_M^{\nu_1,\nu_2}$
be the corresponding locally closed subscheme of $\Conv_M$, which is fibered over $\Gr_M^{\nu_1}$ with the
typical fiber being $\Gr_M^{\mu_2}$. If $\mu_3$ is a third element of $\Lambda_M^{+}$, we will denote by
$\Conv_M^{\mu_1,\mu_2;\mu_3}$ the intersection $\Conv_M^{\nu_1,\nu_1}\cap\,pr^{-1}(\Gr_M^{\nu_3})$.

We will need the following fact, which follows easily from \propref{lusztig} (or alternatively, 
from \propref{dimestimate}):
\begin{equation}  \label{convestimate}
\on{dim}(\Conv_M^{\mu_1,\mu_2;\mu_3})=\langle\mu_1+\mu_2+\mu_3,\rhoch_M\rangle
\end{equation}

Consider the composition
\begin{equation*}  \label{morestacksleft}
Z_{P,M}^{\nu',?,\eta'}\overset{'\hl_M}\longrightarrow
{}_{x,\nu'}\BunPbw\simeq \BunPbw\underset{\Bun_M}\times {}_x\H^{\nu'}_M\to\BunPbw.
\end{equation*}
It follows from the definitions that $Z_{P,M}^{\nu',?,\eta'}$ is a locally trivial fibration
over $\BunPbw$ with the typical fiber $\Conv_M^{\nu',\eta'}$. More precisley,
we obtain the following identifications of stacks:
$$
\CD
Z_{P,M}^{\nu',?,\eta'}  @>{\sim}>> \Conv_M^{\nu',\eta'} \overset{M(\O_x)}\times 
(\Bun_P\underset{\Bun_M}\times {}_x\M) \\
@V{'\hl_M}VV        @V{\on{id}\times pr'}VV    \\
_{x,\nu'}\BunPbw  @>{\sim}>>  \Gr_M^{\nu'} \overset{M(\O_x)}\times 
(\Bun_P\underset{\Bun_M}\times {}_x\M),
\endCD
$$
where $_x\M$ denotes the canonical $M(\O_x)$--torsor over $\Bun_M$. Similarly, we have the identifications
$$
\CD
Z_{P,M}^{?,\eta,\eta'}  @>{\sim}>> \Conv_M^{\eta,-\w_0^M(\eta')} \overset{M(\O_x)}\times 
(\Bun_P\underset{\Bun_M}\times {}_x\M) \\
@V{'\hr_M}VV        @V{\on{id}\times pr'}VV    \\
_{x,\eta}\BunPbw  @>{\sim}>>  \Gr_M^{\eta} \overset{M(\O_x)}\times
(\Bun_P\underset{\Bun_M}\times {}_x\M).
\endCD
$$

\begin{lem} \label{descrnonprince}
Under the above identifications, the stack $Z_{P,M}^{\nu',\eta,\eta'}$
admits the following description:

{\em (a)} 
Inside $Z_{P,M}^{\nu',?,\eta'}$, it is identified with
$\Conv_M^{\nu',\eta';\eta} \overset{M(\O_x)}\times 
(\Bun_P\underset{\Bun_M}\times {}_x\M)$.

\smallskip

{\em (b)}
Inside $Z_{P,M}^{?,\eta,\eta'}$, it is identified with
$\Conv_M^{\eta,-\w_0^M(\eta');\nu'} \overset{M(\O_x)}\times (\Bun_P\underset{\Bun_M}\times {}_x\M)$.

\end{lem}

\smallskip

Let us consider the $*$--restriction of $(\A^{-\w_0^M(\nu)}_M\tboxtimes \IC_{\BunPbw})^r$ to 
$Z_{P,M}^{?,\eta,\eta'}$. By definition, this complex
lives in cohomological degrees $\leq 0$ (and the inequality is strict unless $\eta'=\nu$ and $\eta=0$). 
Moreover, it is constant along the fibers of the projection
$$Z_{P,M}^{?,\eta,\eta'}\to {}_{x,\eta}\BunPbw.$$

\smallskip 

Therefore, from the description of $Z_{P,M}^{\nu',\eta,\eta'}$ given in \lemref{descrnonprince}(b), we infer that
when we restrict this complex further to $Z_{P,M}^{\nu',\eta,\eta'}$ it will live in cohomological degrees
$$\leq -\on{codim}(\Conv_M^{\eta,-\w_0^M(\eta');\nu'},\Conv_M^{\eta,-\w_0^M(\eta')})=
-\langle \eta+\eta'-\nu',\rhoch_M\rangle,$$
where the last inequality follows from \eqref{convestimate}. Moreover, 
the inequality is strict unless $\eta=0$ and $\eta'=\nu$.

\smallskip

Let us observe now that according to \lemref{descrnonprince}(a), the fibers of the map
$$Z_{P,M}^{\nu',\eta,\eta'}\to {}_{x,\nu'}\BunPbw$$ can be identified with the fibers of the projection
$$pr': \Conv_M^{\nu',\eta';\eta}\to \Gr_M^{\nu'}.$$ 
By applying again \eqref{convestimate} we conclude that
those have dimension $\leq\langle \eta+\eta'-\nu',\rhoch_M\rangle$.

\medskip

This implies that when $\eta\neq 0$ or $\eta'\neq\nu$ the complex $K_M^{\nu',\eta,\eta'}$ lives in negative
cohomological degrees and it remains to analyze the case $\eta=0$, $\eta'=\nu$. 

\smallskip

However, in the latter case $Z_{P,M}^{?,0,\nu}$ coincides with 
$Z_{P,M}^{-\w_0^M(\nu),0,\nu}$ and its projection
onto $_{x,\nu}\BunPbw$ is an isomorphism. This implies the last statement of the proposition, since
the $*$--restriction of $(\A^{-\w_0^M(\nu)}_M\tboxtimes \IC_{\BunPbw})^r$
to $Z_{P,M}^{?,0,\nu}$ obviously coincides with $\IC_{Z_{P,M}^{?,0,\nu}}$.

\ssec{Proof of \thmref{mainHeckenonprinceG}}

\sssec{}

For the proof of \thmref{mainHeckenonprinceG} we need to recall several additional facts about affine Grassmannians.

\smallskip

Recall the lattice $\Lambda_{G,P}$ (cf. \secref{BunPb}), which by definition can be identified with 
the lattice of characters of the torus $Z(\check M)$. Let us observe, that
$\Lambda_{G,P}$ can be identified, in addition, with the set of connected components of the affine Grassmannians
$\Gr_{M/[M,M]}$ and $\Gr_M$.

\smallskip

For $\theta\in \Lambda_{G,P}$ consider the closed sub-scheme $\overline S^\theta_P$ of $\Gr_G$ 
that corresponds to those pairs $(\F_G, \beta:\F^0_G|_{\D^*_x}\to\F_G|_{\D^*_x})$, for which the meromorphic map 
$$\L^{\lambdach}_{\F^0_T}\to \V^{\lambdach}_{\F^0_G}\overset{\beta}\longrightarrow 
\V^{\lambdach}_{\F_G}$$ 
has a pole of order $\leq \langle\theta,\lambdach\rangle$ for every 
$\lambdach\in \check\Lambda_{G,P}$.

Similarly, let $S^\theta_P$ be an open subscheme of $\overline S^\theta_P$ that corresponds to 
the pairs
$(\F_G, \beta)$ as above for which the maps
$$\L^{\lambdach}_{\F^0_T}(-\langle\theta,\lambdach\rangle\cdot x)\to \V^{\lambdach}_{\F_G},
\,\,\lambdach\in \check\Lambda_{G,P}$$
have no zero either.

\smallskip

The subscheme $S^\theta_P$ is stable under the action of the group $[P,P](\K_x)M(\O_x)$
on $\Gr_G$. Moreover, the action of an even smaller group, namely, of $[P,P](\K_x)$ on $S^\theta_P$
is transitive.

\medskip

The following assertion is well-known: 

\begin{prop}
For each $\theta\in\Lambda_{G,P}$ there exists a unique map of indschemes
$$\tf^\theta_P:S^\theta_P\to \Gr_M$$ with the following property:

\noindent For $(\F_G, \beta:\F^0_G|_{\D^*_x}\to\F_G|_{\D^*_x})\in S^\theta_P$, its image
$(\F_M,\beta_M)$ is the unique point of $\Gr_M$, for which the composition
$$(\V^{U(P)})_{\F_M}\overset{\beta_M}\longrightarrow (\V^{U(P)})_{\F^0_M}\to 
\V_{\F^0_G}\overset{\beta}\longrightarrow \V_{\F_G}$$
has neither zero nor pole $\forall\lambdach\in\check\Lambda_G^{+}$.
\end{prop}

By definition, $\tf^\theta_P$ takes values in the connected component of $\Gr_M$ corresponding to $\theta$. 
Let $\nu$ be an element of $\Lambda_M^{+}$ and let $\theta$ be its image under the projection 
$\Lambda\to\Lambda_{G,P}$. In what follows, we will denote by
$S^\nu_P$ the pre-image $\tf_P^\theta{}^{-1}(\Gr_M^\nu)\subset S^\theta_P$. 
We will denote by $\tf_P^\nu$ the corresponding map $S^\nu_P\to\Gr_M^\nu$.

\medskip

The next result is a consequence of \propref{dimestimate}:

\begin{prop} \label{estimnonprince}

{\em (a)}
Let $\nu$ (resp., $\lambda$) be a dominant integral coweight of $M$ (resp., of $G$) The intersection
$S^\nu_P\cap\Gr_G^\lambda$ has dimension $\leq \langle\nu+\lambda,\rhoch\rangle$ and the fibers 
of the projection
$$\tf^\nu_P:S^\nu_P\cap\Gr_G^\lambda\to\Gr_M^\nu$$
are of dimension $\leq \langle\nu+\lambda,\rhoch\rangle-\langle\nu,2\rhoch_M\rangle$.

\smallskip

{\em (b)}
For $\lambda\in\Lambda_G^{+}$ and $\theta\in\Lambda_{G,P}$, the direct image 
$$\tf_P^\theta{}_{!}(\A^\lambda_G|_{S_P^\theta})\otimes 
(\Ql[1](\frac{1}{2}))^{\otimes \langle \theta,2(\rhoch-\rhoch_M)\rangle}$$
lives in cohomological degrees $\leq 0$. (In the above formula we have used the fact that 
$2(\rhoch-\rhoch_M)\in\check\Lambda_{G,P}$.) Its $0$-th cohomology belongs to $\Sph_M$.

\smallskip

{\em (c)}
The functor $\Sph_G\to \Sph_M$ given by
$$\S\to \underset{\theta}\oplus h^0(\tf_P^\theta{}_{!}(\S|_{S_P^\theta})\otimes 
(\Ql[1](\frac{1}{2}))^{\otimes \langle \theta,2(\rhoch-\rhoch_M)\rangle})$$
has a natural structure of a tensor functor.

\end{prop}

\medskip

Finally, we will use the following theorem, which follows from \thmref{dualgroup} (cf. \cite{BD}):

\begin{thm} \label{dualgroupnonprince}
The tensor functor 
$$\S\to \underset{\theta}\oplus h^0(\tf_P^\theta{}_{!}(\S|_{S_P^\theta})\otimes 
(\Ql[1](\frac{1}{2}))^{\otimes \langle \theta,2(\rhoch-\rhoch_M)\rangle})$$
of \propref{estimnonprince}(c) is canonically isomorphic to the functor $\on{gRes}^G_M$ that makes the following 
diagram of categories commutative:
$$
\CD
\Sph_G @>{\on{gRes}^G_M}>> \Sph_M \\
@V{F_G}VV    @V{F_M}VV  \\
\on{Rep}(\check G)  @>{\on{Res}^G_M}>> \on{Rep}(\check M)
\endCD
$$
\end{thm}

\begin{cor} \label{nonprinceestimatecor}
The complex 
$$\tf_P^\nu{}_{!}(\Ql|_{\Gr_G^\lambda})\otimes 
(\Ql[1](\frac{1}{2}))^{\otimes \langle\lambda,2\rhoch\rangle+\langle\nu,2(\rhoch-\rhoch_M)\rangle}$$
lives in the cohomological degrees $\leq 0$ and its $0$-th perverse cohomology is canonically isomorphic
to $\IC_{\Gr^\nu_M}\otimes \on{Hom}_{\check M}(U^\nu,V^\lambda)$.
\end{cor}

\sssec{} \label{fulldescrG}

Now we will proceed to the proof of \thmref{mainHeckenonprinceG}. The argument will be a direct 
generalization of the ones of \thmref{mainHeckeprinceatx} and \thmref{mainHeckenonprinceM}.

\medskip

As the complex $_xH_{P,G}^\lambda(\IC_{\BunPbw})$ is self-dual, in order to prove that it is isomorphic 
to the direct sum 
$\underset{\nu\in\Lambda_M^{+}}\oplus \IC_{_{x,\geq\nu}\BunPbw}\otimes \on{Hom}_{\check M}(U^{\nu},V^{\lambda})$,
it is enough to show that the complex
\begin{equation} \label{insteadofprop}
j_{\nu}^*\circ {}'\hl_G{}_{!}(\A_G^{-\w_0(\lambda)}\tboxtimes\IC_{\BunPbw})^r 
\end{equation}
lives in non-positive cohomological degrees and that its $0$-th cohomology can be identified with
$\IC_{_{x,\nu}\BunPbw}\otimes\on{Hom}_{\check M}(U^{\nu},V^\lambda)$.

\medskip

For $\nu,\nu'\in\Lambda_M^{+}$ and $\lambda'\in\Lambda_G^{+}$, let $Z_{P,G}^{\nu,?,\lambda{}'}$ 
and $Z_{P,G}^{?,\nu',\lambda{}'}$
denote the locally closed sub-stacks of $_{x,\infty}Z_{P,G}$ equal to
\begin{equation*}
'\hl_G{}^{-1}({}_{x,\nu}\BunPbw)\cap{}'\p^{-1}({}_x\H_G^{\lambda{}'}) \text{ and }
'\hr_G{}^{-1}({}_{x,\nu'}\BunPbw)\cap{}'\p^{-1}({}_x\H_G^{\lambda{}'}), \text{ respectively}.
\end{equation*}

Let $_x\P$ denote the canonical $P(\O_x)$-torsor over $_{x,0}\BunPbw$ (cf. \secref{describeprince}) 
and let $_x\overset{o}\P$ denote its restriction to the open substack $\Bun_P\subset{}_{x,0}\BunPbw$.

\smallskip

We have the following identifications of stacks:
$$
\CD
Z_{P,G}^{\nu,?,\lambda{}'} @>{'\hl_G}>> {}_{x,\nu}\BunPbw @>>> {}_{x,0}\BunPbw  \\
@V{\simeq}VV        @V{\simeq}VV        @V{\on{id}}VV  \\
(\Gr_G^{\lambda'}\times\Gr_M^{\nu}) \overset{P(\O_x)}\times {}_x\P @>>>
\Gr_M^{\nu}\overset{P(\O_x)}\times {}_x\P@>>> {}_{x,0}\BunPbw
\endCD
$$
and
$$
\CD
Z_{P,G}^{?,\nu',\lambda'} @>{'\hr_G}>> {}_{x,\nu'}\BunPbw @>>> {}_{x,0}\BunPbw  \\
@V{\simeq}VV        @V{\simeq}VV        @V{\on{id}}VV  \\
(\Gr_G^{-\w_0(\lambda')}\times \Gr_M^{\nu'})\overset{P(\O_x)}\times {}_x\P @>>>
\Gr_M^{\nu'}\overset{P(\O_x)}\times {}_x\P @>>> {}_{x,0}\BunPbw.
\endCD
$$

\medskip

Now we need to introduce one more piece of notation. Let $\mu_1,\mu_2,\mu_3$ be elements of $\Lambda_M^{+}$ 
and let $\xi$ be an element of $\Lambda_G^{+}$. Let us denote by $W^{\xi,\mu_1,\mu_2;\mu_3}$ the following
locally closed subscheme of $\Gr_G\times\Gr_M$:

\smallskip 

\noindent $W^{\xi,\mu_1,\mu_2;\mu_3}$ consists of those quadruples 
$$(\F_G,\beta:\F_G|_{\D^*_x}\to\F^0_G|_{\D^*_x},
\F_M,\beta_M:\F_M|_{\D^*_x}\to\F^0_M|_{\D^*_x})$$
for which 
\begin{itemize}

\item
$(\F_G,\beta)\in S^{\mu_1}_P\cap \Gr_G^{\xi}$. 

\smallskip

\item
$(\F_M,\beta_M)\in\Gr_M^{\mu_2}$.

\smallskip

\item
$((\F_M,\beta_M)\times \tf_P^{\mu_1}(\F_G,\beta))\in \Gr_M^{\mu_2}\times \Gr_M^{\mu_1}$
lies in the image of $\Conv_M^{\mu_1,\mu_3;\mu_2}$ under the projection
$pr\times pr':\Conv_M\to\Gr_M\times\Gr_M$.

\end{itemize}

Note that one can rephrase the third condition in the following way: 
the point $(\F_M,\beta_M)\in\Gr_M$ is in position $\mu_3$ with respect to $\tf_P^{\mu_1}(\F_G,\beta)$.

\medskip

Finally, for $\nu,\nu',\eta\in\Lambda_M^{+}$ and $\lambda'\in\Lambda_G^{+}$ we define the locally closed
sub-stacks $Z_{P,G}^{\nu,\nu',\eta,\lambda'}$ and $\widetilde Z_{P,G}^{\nu,\nu',\eta,\lambda'}$ of
$_{x,\infty}Z_{P,G}$ as follows:

\noindent

Using the identification of the first of the above commutative diagrams,
$Z_{P,G}^{\nu,\nu',\eta,\lambda'}$ is equal to
\begin{equation*}
W^{\lambda',\eta,\nu;\nu'} \overset{P(\O_x)}\times {}_x\P\subset 
(\Gr_G^{\lambda'}\times \Gr_M^{\nu})\overset{P(\O_x)}\times {}_x\P.
\end{equation*}
The stack $\widetilde Z_{P,G}^{\nu,\nu',\eta,\lambda'}$ is defined via the second diagram as
\begin{equation*}
W^{-\w_0(\lambda'),-\w_0^M(\eta),\nu';\nu} \overset{P(\O_x)}\times {}_x\P \subset 
(\Gr_G^{-\w_0(\lambda')}\times\Gr_M^{\nu'}) \overset{P(\O_x)}\times {}_x\P.
\end{equation*}

\begin{lem} \label{descrG}
The sub-stacks $Z_{P,G}^{\nu,\nu',\eta,\lambda'}\subset Z_{P,G}$ and 
$\widetilde Z_{P,G}^{\nu,\nu',\eta,\lambda'}\subset Z_{P,G}$ coincide.
\end{lem}

\bigskip

We define the complex $K^{\nu,\nu',\eta,\lambda'}$ on $_{x,\nu}\BunPbw$ as the direct image under the map
$'\hl_G:{}Z_{P,G}^{\nu,\nu',\eta,\lambda'}\to {}_{x,\nu}\BunPbw$ of the $*$-restriction of 
$(\A_G^{-\w_0(\lambda)}\tboxtimes\IC_{\BunPbw})^r$ 
to $Z_{P,G}^{\nu,\nu',\eta,\lambda'}$.

\sssec{}

To establish the required isomorphism, it is enough to show that 

\smallskip

\noindent a) The complex $K^{\nu,\nu',\eta,\lambda'}$ lives in cohomological degrees $\leq 0$ and
the inequality is strict unless $\nu'=0$, $\lambda'=\lambda$ and $\eta=\nu$.

\smallskip

\noindent b) The $*$--restriction of $h^0(K^{\nu,0,\nu,\lambda})$ to 
$_{x,\nu}\BunPbw-{}_{x,\nu}\Bun_P$
vanishes.

\smallskip

\noindent c) $h^0(K^{\nu,0,\nu,\lambda})$ restricted to $_{x,\nu}\Bun_P$ is isomorphic to
$\IC_{_{x,\nu}\Bun_P}\otimes\on{Hom}_{\check M}(U^{\nu},V^\lambda)$.

\bigskip

First of all, as in the proofs of \thmref{mainHeckeprinceatx} and \thmref{mainHeckenonprinceM}, the $*$--restriction
of $(\A_G^{-\w_0(\lambda)}\tboxtimes\IC_{\BunPbw})^r$ to the substack
$Z_{P,G}^{\nu,\nu',\eta,\lambda'}=\widetilde Z_{P,G}^{\nu,\nu',\eta,\lambda'}$
lives in cohomological degrees 
$$\leq -\on{codim}(W^{-\w_0(\lambda'),-\w_0^M(\eta),\nu';\nu},\Gr_G^{-\w_0(\lambda')}\times\Gr_M^{\nu'})$$
and it follows from \eqref{convestimate} and \propref{estimnonprince}(a) that
$$\on{codim}(W^{-\w_0(\lambda'),-\w_0^M(\eta),\nu';\nu},\Gr_G^{\lambda'}\times\Gr_M^{\nu'})=
\langle\lambda',\rhoch\rangle+\langle \w_0^M(\eta),\rhoch\rangle+\langle\nu'-\w_0^M(\eta)-\nu,\rhoch_M\rangle.$$

\smallskip

Now, the fibers of the map $'\hl_G:{}Z_{P,G}^{\nu,\nu',\eta,\lambda'}\to{}_{\nu,x}\BunPbw$ 
have the same dimension as the fibers of the map
$W^{\lambda',\eta,\nu;\nu'}\to \Gr_M^{\nu}$, and the latter equals
$\langle\lambda',\rhoch\rangle+\langle \eta,\rhoch\rangle+\langle\nu'-\nu-\eta,\rhoch_M\rangle$.

\smallskip

As $\langle\eta,\rhoch-\rhoch_M\rangle=\langle \w_0^M(\eta),\rhoch-\rhoch_M\rangle$, point (a) above follows. 
Point (b) follows in the same way as point (b) of \propref{ppropmainHeckeprinceatx}.

\smallskip

Let us consider now the case $\lambda'=\lambda$, $\nu'=0$. (In this case $\eta$ is automatically equal to $\nu$.) 
Let $\overset{o}Z{}^{\nu,0,\nu,\lambda}$ denote the pre--image of $_{\nu,x}\Bun_P$ in 
$Z_{P,G}^{\nu,0,\nu,\lambda}$ under the map $'\hl_G$. Then $\overset{o}Z{}^{\nu,0,\nu,\lambda}$ 
is also the pre-image of $\Bun_P$ under the map $'\hr_G$.

We have the following commutative diagram of stacks:
$$
\CD
\overset{o}Z{}^{\nu,0,\nu,\lambda}  @>{\sim}>> (S_P^{\nu}\cap\Gr_G^\lambda) \overset{P(\O_x)}\times 
{}_x\overset{o}\P \\
@V{'\hl_G}VV       @V{\on{id}\times \tf^{\nu}_P}VV   \\
_{\nu,x}\Bun_P  @>{\sim}>>  \Gr_M^{\nu} \overset{P(\O_x)}\times{}_x\overset{o}\P.
\endCD
$$

Moreover, as $\Bun_P$ is smooth, the restriction of $(\A_G^{-\w_0(\lambda)}
\tboxtimes\IC_{\BunPbw})^r$
to $\overset{o}Z_{P,G}{}^{\nu,0,\nu,\lambda}$ is the twisted external product
\begin{equation*}
(\Ql[1](\frac{1}{2}))^{\otimes \langle\lambda,2\rhoch\rangle 
+2\langle \nu,\rhoch-\rhoch_M\rangle}_{S_P^{\nu}\cap\Gr_G^\lambda} \widetilde{\boxtimes} \IC_{\Bun_P}.
\end{equation*}

\smallskip

Therefore, since the group $P(\O_x)$ is connected, 
\corref{nonprinceestimatecor} implies that
the restriction of $h^0(K^{\nu,0,\nu,\lambda})$ to $_{x,\nu}\Bun_P$ can be 
identified with 
$\IC_{_{x,\nu}\Bun_P}\otimes \on{Hom}_{\check M}(U^{\nu},V^\lambda)$.

\medskip

Hence, we have established the isomorphism
$$_xH_{P,G}^\lambda(\IC_{\BunPbw})\simeq 
\underset{\nu\in\Lambda_M^{+}}\oplus 
\IC_{_{x,\geq\nu}\BunPbw} \otimes \on{Hom}_{\check M}(U^{\nu},V^{\lambda}).$$ 

The fact that this isomorphism is compatible with the tensor structure, is a corollary of point (c) of 
\propref{estimnonprince} combined with \thmref{dualgroupnonprince}.

\section{The acyclicity theorem}

\ssec{The notion of local acyclicity}

Let $f:Y_1\to Y_2$ be a map between smooth algebraic varieties and let $\S$ be an object of $\Sh(Y_1)$. In
\cite{SGA41/2}, there was introduced the notion of universal local acyclicity (we will abbreviate it to ULA)
of $\S$ with respect to $f$, which we are now going to review.

\sssec{}  \label{intracyc}

Let $g:Y'\to Y''$ be a map between algebraic varieties. If $\S_1,\S_2$ are two objects of $\Sh(Y'')$, there
exists a natural morphism in $\Sh(Y')$:
\begin{equation}  \label{intmap}
g^*(\Hhom(\S_1,\S_2))\to \Hhom(g^*(\S_1),f^*(\S_2)),
\end{equation}
where $\Hhom$ is the internal $\on{Hom}$.

\smallskip

Assume now that both $Y'$ and $Y''$ are smooth and take $$\S_2:={\mathbb D}_{Y''}\simeq 
(\Ql(1)[2])^{\otimes \on{dim}(Y'')}.$$
Then \eqref{intmap} yields us a functorial map
$$g^*({\mathbb D}(\S))\otimes(\Ql(1)[2])^{\otimes \on{dim}(Y')-\on{dim}(Y'')} \to {\mathbb D}(g^*(\S)),$$
and by replacing $\S$ by ${\mathbb D}(\S)$ we obtain a functorial map
$$can_g:g^*(\S)\otimes (\Ql(\frac{1}{2})[1])^{\otimes \on{dim}(Y')-\on{dim}(Y'')}\to
g^{!}(\S)\otimes (\Ql(\frac{1}{2})[1])^{\otimes \on{dim}(Y'')-\on{dim}(Y')}.$$

\smallskip

The above natural transformation $can_g$ is clearly an isomorphism when $g$ is smooth.
Let us point out that in general, the map $can_g$ is not expressible via the standard six functors. 
As an incarnation of this fact, it is rather hard to define it on the level of D-modules. 

\sssec{} \label{propacyc}

Now let $f:Y_1\to Y_2$ be a map between smooth algebraic varieties. 
Take $Y'=Y_1$, $Y''=Y_1\times Y_2$ and $g=\Gamma_f$.
According to the above discussion, for any $\S\in\Sh(Y_1)$ and $\T\in\Sh(Y_2)$, we obtain a canonical map
\begin{align*}
can_{\Gamma_f}:\S\otimes f^*(\T)\otimes (\Ql(\frac{1}{2})[1])^{\otimes -\on{dim}(Y_2)} \to
\S\overset{!}\otimes f^{!}(\T)\otimes(\Ql(\frac{1}{2})[1])^{\otimes \on{dim}(Y_2)},
\end{align*}
where $\overset{!}\otimes$ is the ${\mathbb D}$-conjugate of $\otimes$, i.e.
$\S_1\overset{!}\otimes\S_2:={\mathbb D}({\mathbb D}(\S_1)\otimes{\mathbb D}(\S_2))$.

\medskip

\noindent{\bf Definition.} 
{\it An object $\S\in\Sh(Y_1)$ is said to be locally acyclic with respect to $f$ if $can_{\Gamma_f}$ 
is an isomorphism $\forall \,\T\in\Sh(Y_2)$.
A sheaf $\S$ is universally locally acyclic (ULA) with respect to $f$, if the above property 
holds after any smooth base change $Y'_2\to Y_2$.}

\smallskip

For instance, any sheaf on $Y_1$ is ULA with respect to the projection $Y_1\to \on{pt}$. 
When $f$ is the identity map $Y_1\to Y_1$,
one can show that $\S$ is ULA if and only of it is locally constant.

\smallskip

\noindent{\it Remark.}
The above definition of universal local acyclicity is {\it a priori} weaker than the one given
in \cite{SGA41/2}. (In other words, if a complex $\S$ is ULA in the sense of {\it loc. cit.},
then it is ULA in our sense as well.) We conjecture that the two notions are in fact equivalent. 

However, for the purposes of this paper, it does not matter which definition to use:
the proof of the main result, \thmref{acycthm}, works for both of them. 

\medskip

Here are the most immediate properties of the ULA condition:

\begin{enumerate}  

\item{}

Let $\S$ be an ULA sheaf on $Y_1$ with respect to $f:Y_1\to Y_2$. 
Let $Y'$ be another smooth variety and $\S'$ be an arbitrary sheaf on it.
Then the sheaf $\S\boxtimes\S'$ on $Y_1\times Y'$ is ULA with respect to the map 
$Y_1\times Y'\to Y_1\overset{f}\to Y_2$.
(In most applications we will have $Y_2=Y_1$, $\S:={\Ql}_{Y_1}$.)

\item{}

The ULA condition is local in the smooth topology on the source. In other words, if $s:Y_1\to Y'_1$ 
is a smooth (resp., smooth and surjective)
map and $\S$ is a sheaf on $Y'_1$, then the sheaf $s^*(\S)$ on $Y_1$ is ULA with respect 
to $f:=f'\circ s$ if (resp., if and only if) the  
sheaf $\S$ on $Y'_1$ is ULA with respect to $f'$.

\item{}

If $s:Y_1\to Y'_1$ is a proper map (resp., closed embedding) and $\S$ is a sheaf on $Y_1$, then the sheaf 
$s_{*}(\S)$ on $Y'_1$ is ULA with respect to $f'$ (we are assuming that $f$ factorizes as $f=f'\circ s$) if
(resp., if and only if) the sheaf $\S$ is ULA with respect to $f$. \newline
This property allows to formulate the ULA condition in the situation when $Y_1$ is not necessarily smooth: 
it is enough to (locally) embed
$Y_1$ as a closed sub--scheme into some $Y'_1$ (with $f=f'\circ s$ and $Y'_1$ being smooth) and 
to require that the sheaf 
$s_*(\S)$ on $Y'_1$ is ULA with respect to $f'$.

\item{}

A sheaf $\S$ is ULA if and only ${\mathbb D}(\S)$ is.

\item{}

Let $\S$ is a ULA sheaf on $Y_1$ with respect to $f:Y_1\to Y_2$ and consider the functor 
$\Sh(Y_2)\to \Sh(Y_1)$
given by $\T\to 
\S\otimes f^*(\T)\otimes (\Ql(\frac{1}{2})[1])^{\otimes -\on{dim}(Y_2)}$.
Then it commutes with the Verdier duality in the sense that there is an isomorphism of functors:
$${\mathbb D}(\S\otimes f^*(\T)\otimes (\Ql(\frac{1}{2})[1])^{\otimes -\on{dim}(Y_2)})\simeq
\DD(\S)\otimes f^*(\DD(\T))\otimes (\Ql(\frac{1}{2})[1])^{\otimes -\on{dim}(Y_2)}.$$
Moreover, when $\S$ is concentrated in non-positive (resp., non-negative) cohomological degrees,
the above functor is right (resp., left) exact. 

\item{}

If $t:Y_2\to Y'_2$ is a smooth map and $\S$ is a sheaf on $Y_1$ which is ULA with respect to 
$f:Y_1\to Y_2$, then $\S$ is ULA with respect to $t\circ f$.

\end{enumerate}

\bigskip

An important technical tool is provided by the following theorem (cf. \cite{SGA41/2}):

\begin{thm} \label{Deligne}
Let $f:Y_1\to Y_2$ be a map, where $Y_1$ is a scheme of finite type and
$Y_2$ is a smooth variety. Let $\S$ be an object of $\Sh(Y_1)$. 
Then there exists a non-empty open subvariety $Y^0_1\subset Y_1$, 
such that $\S$ is ULA when restricted to $f^{-1}(Y^0_1)\subset Y_2$.
\end{thm}

\sssec{}  \label{proofofgoodtopullback}

The main result of this section is the following theorem:

\begin{thm} \label{acycthm}

The sheaves $\IC_{\BunPbw}$ and $\jw_P{}_{!}(\IC_{\Bun_P})$ on $\BunPbw$ are ULA with respect to 
the map $\qw_P:\BunPbw\to\Bun_M$.

\end{thm}

The proof will be given in the next two subsections. However, we will present the proof only of the first 
part of \thmref{acycthm}, 
namely that $\IC_{\BunPbw}$ is ULA with respect to $\qw_P$. The proof for $\jw_P{}_{!}(\IC_{\Bun_P})$ 
is absolutely analogous.

\smallskip

Let us show how \thmref{acycthm} implies \thmref{goodtopullbacknonprince} and, in particular, 
\thmref{goodtopullbackprince}.

\begin{proof} (of \thmref{goodtopullbacknonprince})

Point (a) of \thmref{goodtopullbacknonprince} follows immediately from Property 5 of \secref{propacyc}.
Moreover, by the same reason, for a perverse sheaf $\T$ on $\Bun_M$, the sheaf $\qw_P^{!*}(\T)$
is perverse.

\smallskip

Therefore, to prove point (b) of the Theorem, it is enough to show that whenever $\T$ is irreducible, 
$\qw_P^{!*}(\T)$ is irreducible as well. First of all, since the map ${\mathfrak q}:\Bun_P\to\Bun_M$
is smooth (and has connected fibers), it is clear that 
${\mathfrak q}_P^{!*}(\T)\simeq\qw_P^{!*}(\T)|_{\Bun_P}$
is irreducible. Since the situation is Verdier self-dual, it is enough to show, therefore, that the 
$*$-restriction of $\qw_P^{!*}(\T)$ to $\BunPbw-\Bun_P$ lives in the cohomological degrees $<0$. 

Let $K$ be the cone of the map $\jw_P{}_{!}(\IC_{\Bun_P})\to \IC_{\BunPbw}$. By definition,
$K$ lives in the cohomological degress $<0$. 
\footnote{The last assertion uses, of course, \propref{densenonprince}, which
will be independently proven later (cf. \secref{stra1} and \secref{stra2}).} 
Moreover, \thmref{acycthm} implies that the sheaf $K$ is also ULA with respect to $\qw_P$. 

However, $$\qw_P^{!*}(\T)|_{\BunPbw-\Bun_P}\simeq 
K\otimes \qw_P^*(\T)\otimes (\Ql(\frac{1}{2})[1])^{\otimes -\on{dim}(\Bun_M)}$$
and the required assertion follows from the Property 5 of \secref{propacyc}.

\end{proof}

\ssec{Proof of \thmref{acycthm} in the Borel case}

\sssec{}

We will use the following observation: 

Let $H$ be a group (or a group-stack), let $Y\overset{f}\longrightarrow H$ be
an arbitrary map and let $\S$ be a sheaf on $Y$. Consider the composition 
$$m_f:H\times Y\overset{f\times\on{id}}\longrightarrow H\times H\overset{\on{mult}}\longrightarrow H.$$

We claim that the sheaf ${\Ql}_H\boxtimes\S$ is always ULA with respect to $m_f$. Indeed, 
the automorphism $(h,y)\mapsto (h\cdot f(y),y)$ transforms the data of $(m_f,{\Ql}_H\boxtimes\S)$ into a 
direct product situation (cf. Property 1 of \secref{propacyc}).

\medskip

The proof of \thmref{acycthm} will be essentially a reduction of our situation ($\q:\BunBb\to\Bun_T$) to the
one mentioned above, where the role of $H$ will be played by $\Bun_T$. Namely, we will construct a stack $Z$ that
fits into a commutative diagram:

$$
\CD
Z  @>>>  \BunBb \\
@VVV   @V{\q}VV    \\
\Bun_T\times\BunBb  @>{m_{\q}}>> \Bun_T
\endCD
$$
where the upper horizontal arrow is surjective and smooth and the left vertical arrow is just smooth. 
This will prove \thmref{acycthm} in view of Property 2 of \secref{propacyc}.

\sssec{}

Let us choose elements $\lambda_1,...,\lambda_r\in\Lambda_G^{+}$ in such a way that 
they form a basis for $\Lambda\underset{\ZZ}\otimes{\mathbb Q}$. Let $m$ be an integer $\geq 2g-1$ and consider the
product $X^{m\cdot r}-\Delta$, where $\Delta$ denotes the divisor of diagonals.

\smallskip

Let $\H^{?}_G$ denote the following version of the Hecke stack:
$\H^{?}_G$ classifies the data of 
$$(\{x_{1,1},...,x_{1,m},x_{2,1},...,x_{r,1},...,x_{r,m}\}\in X^{m\cdot r}-\Delta,\,\,\F_G,\F'_G,\,\beta),$$
where $\beta$ is an isomorphism 
$$\beta:\F_G\simeq\F'_G|_{X-\{x_{1,1},...,x_{1,m},x_{2,1},...,x_{r,1},...,x_{r,m}\}}$$
such that for every $i$ and $j$, $\F_G$ is in position $\lambda_i$ with respect to $\F'_G$ at $x_{i,j}$.

We let $\hl_G$ and $\hr_G$ denote the projections from $\H^{?}_G$ to $\Bun_G$ that send the above point 
of $\H^{?}_G$ to $\F_G$ and $\F'_G$, respectively. The projection from $\H^{?}_G\to X^{m\cdot r}-\Delta$
will be denoted by $\pi$.

\smallskip

Let us denote by $\overline{Z}$ the fiber product 
$\H_G^{?}\underset{\Bun_G}\times\BunBb$, where $\H_G$ is mapped  
to $\Bun_G$ by means of $\hr_G$. As in \secref{bd}, we have the second projection 
$\phi:\overline{Z}\to\BunBb\times (X^{m\cdot r}-\Delta)$ and we obtain a commutative diagram:
$$
\CD
\BunBb\times  (X^{m\cdot r}-\Delta) @<{\phi}<< \overline{Z} @>{'\hr_G}>> \BunBb  \\
@V{\p\times\on{id}}VV      @VVV       @V{\p}VV    \\
\Bun_G\times (X^{m\cdot r}-\Delta)  @<{\hl_G\times\pi}<<   \H'_G  @>{\hr_G}>>  \Bun_G.
\endCD
$$

\medskip

In addition, we have the Abel-Jacobi map $\on{AJ}: (X^{m\cdot r}- \Delta)\to\Bun_T$ that maps
$$(x_{1,1},...,x_{1,m},x_{2,1},...,x_{r,1},...,x_{r,m})\to \F_T^0(\underset{i,j}\Sigma\, \lambda_i\cdot x_{i,j})$$
and the map 
$\overline{Z}\overset{'\hr_G}\longrightarrow \BunBb\overset{\q}\longrightarrow\Bun_T$
coincides with the composition
$$\overline{Z}\overset{\phi}\longrightarrow \BunBb\times (X^{m\cdot r}-\Delta)
\overset{\on{id}\times\on{AJ}}\longrightarrow \BunBb\times\Bun_T\overset{m_{\q}}\longrightarrow\Bun_T.$$

\medskip

The sought-for stack $Z$ is defined as an open sub-stack of $\overline{Z}$:

\smallskip

\noindent A point 
$(\{x_{i,j}\},\F_G,\beta,\F'_G,\F'_T,\{\kappa'{}^{\lambdach}:
\L_{\F'_T}^{\lambdach}\hookrightarrow \V^{\lambdach}_{\F'_G}\})$
belongs to $Z$ if the following holds:

a) The map $\kappa'{}^{\lambdach}:\L_{\F'_T}^{\lambdach}\hookrightarrow \V^{\lambdach}_{\F'_G}$ has no zero 
at any of the points $x_{i,j}$, $\forall\lambdach\in\check\Lambda_G^{+}$.

b) The map $\kappa^{\lambdach}:=\L_{\F'_T}^{\lambdach}
(-\underset{i,j}\Sigma\, \langle\lambda_i,\lambdach\rangle\cdot x_{i,j})
\hookrightarrow \V^{\lambdach}_{\F_G}$ has no zero at any of the points $x_{i,j}$ either, 
$\forall\lambdach\in\check\Lambda_G^{+}$.

\medskip

It remains, therefore, to check that $Z$ defined in the above way satisfies all the requirements. 

\smallskip

The projection $'\hr_G:\overline{Z}\to\BunBb$ is smooth by definition, hence, so is the map $Z\to\BunBb$. The fact
that this map is surjective follows from \lemref{descrprince}(b). 

\smallskip

Now, the restriction of $\phi$ onto $Z$ is an isomorphism, as follows from \lemref{descrprince}(c),
and the map $\on{AJ}:(X^{m\cdot r}- \Delta)\to\Bun_T$ is smooth due to the condition that $m\geq 2g-1$.
This finishes the proof of \thmref{acycthm} when $P=B$.

\ssec{Proof of \thmref{acycthm} in the general case}

In principle, it is not difficult to generalize the proof given in the previous subsection to treat the case
of an arbitrary parabolic $P$. However, we will present here a different argument, based on \thmref{Deligne}.

\sssec{}  \label{gencaseacyc}

As a first step, we will exhaust the stack $\BunPbw$ by open substacks which would be of finite type over $\Bun_M$.
This can be done as follows:

\medskip

Consider the stack $\BunPb$ and for an element $\theta\in\Lambda_{G,P}^{\on{pos}}$ consider the open substack 
$\BunPb^{\leq \theta}\subset\BunPb$ that corresponds to parabolic Drinfeld's structures 
$$(\F_G,\L_{\F_{M/[M,M]}}^{\lambdach}\hookrightarrow \V^{\lambdach}_{\F_G},\,\,\lambdach\in\check\Lambda_{G,P}\cap
\check\Lambda_G^{+})$$
whose "total singularity" does not exceed $\theta$. This means that for every $\lambdach\in \check\Lambda_{G,P}\cap
\check\Lambda_G^{+}$, the coherent sheaf $\V^{\lambdach}_{\F_G}/\L_{\F_{M/[M,M]}}^{\lambdach}$ has no torsion
subsheaves of length $>\langle\theta,\lambdach\rangle$.

\smallskip

Let $\BunPbw^{\leq \theta}\subset\BunPbw$ denote the pre-image of 
$\BunPb^{\leq \theta}\subset\BunPb$ under the map $\r_P$.
For instance, when $\theta=0$, $\BunPbw^{\leq 0}$ coincides with $\Bun_P$. 

\begin{lem} \label{ft}
For every $\theta$, the stack $\BunPbw^{\leq \theta}$ is of finite type over $\Bun_M$.
\end{lem}

\smallskip

It is clear that $\BunPbw=\underset{\theta}\cup\,\BunPbw^{\leq \theta}$ and 
as the assertion of \thmref{acycthm} is local on $\BunPbw$, it would be sufficient to prove that 
$\forall\theta\in\Lambda_{G,P}^{\on{pos}}$, $\IC_{\BunPbw^{\leq \theta}}$ is locally acyclic with respect to $\qw_P$. 
Therefore,
in order not to overload the notation, we will now fix some $\theta$ and until the end of this section, we will
replace the notation "$\BunPbw^{\leq \theta}$" simply by "$\BunPbw$". The reader will readily check that this open 
substack is stable under all the manipulations that we are about to perform.

\sssec{}

Let $\overset{o}{\Bun}_M$ denote the maximal open substack of $\Bun_M$, over which $\IC_{\BunPbw}$ is ULA with
respect to $\qw_P$. Due to the finite type property above (\lemref{ft}) and \thmref{Deligne}, 
$\overset{o}{\Bun}_M$ is non-empty.

\medskip

Let $\F_M$ and $\F'_M$ be two 
$\ol{\Fq}$-points of $\Bun_M$. We will write that $\F_M\prec\F'_M$ if the following
condition holds:

\smallskip

\noindent There exists a {\it $G$-dominant} coweight $\lambda\in\Lambda_M^{+}$ such that
the pair $(\F_M,\F'_M)$ is the image under the map
$$\H_M^\lambda\overset{\hl_M\times\hr_M}\longrightarrow \Bun_M\times\Bun_M$$
of some $\Fqb$-point of $\H_M^\lambda$.

\smallskip

Let $\sim$ be the equivalence relation on the set of isomorphism classes of $\overline{\Fq}$-points of $\Bun_M$,
generated by $\prec$. \thmref{acycthm} clearly follows from the next two assertions:

\begin{prop} \label{eqacyc}
If $\F_M\prec\F'_M$ then $\F'_M\in \overset{o}{\Bun}_M$ if and only if $\F_M$ belongs to 
$\overset{o}{\Bun}_M$ too.
\end{prop}

\begin{prop} \label{alleq}
All $\ol{\Fq}$-points of $\Bun_M$ are $\sim$-equivalent.
\end{prop}

We will first prove \propref{alleq}.

\begin{proof}

The assertion of the proposition is in fact an easy corollary of the following fact proven
in \cite{DS}:

\begin{lem}  \label{affinetrivial}
Let $\widetilde M$ be a semi-simple group and let $\F_{\widetilde M}$ and $\F'_{\widetilde M}$ be two
$\widetilde M$-bundles on $\ol{X}$. Then for any $x\in X(\Fqb)$, the restrictions $\F_{\widetilde M}|_{\ol{X}-x}$
and $\F'_{\widetilde M}|_{\ol{X}-x}$ become isomorphic.
\end{lem}

\smallskip

Let $\F_M$ and $\F'_M$ be two $\overline{\Fq}$-points of $\Bun_M$ and let $\F_{\widetilde M}$ and 
$\F'_{\widetilde M}$ denote their reductions to $\widetilde M:=M/Z^0(M)$, where $Z^0(M)$ is the connected 
component of the identity of $Z(M)$.

\smallskip

The map $\Lambda_G^{+}\to \Lambda_{\widetilde M}^{+}$ is surjective, therefore one can find 
$\lambda\in \Lambda_G^{+}$
and a point of $\H_M^\lambda$ which projects to a point $(\F_M,\F''_M)\in\Bun_M\times\Bun_M$, where 
$\F'_M$ and $\F''_M$ have the same reduction to $\widetilde M$.

\smallskip

We have: $\F''_M=\F'_M\otimes \F_{Z^0(M)}$, where $\F_{Z^0(M)}$ is some $Z^0(M)$-bundle. Now the proof follows from
the fact that $\Lambda_G^{+}\cap \Lambda_{Z^0(M)}$ spans $\Lambda_{Z^0(M)}\underset{\ZZ}\otimes{\mathbb Q}$.

\end{proof}

\sssec{}

\begin{proof} (of \propref{eqacyc}.)

Let us first prove the "if" part of the proposition. We will 
fix $\lambda\in\Lambda_G^{+}$ and let us denote by $\overline Z$ the
fiber product $\H^{\lambda}_G\underset{\Bun_G}\times\BunPbw$, where $\H^{\lambda}_G$
is mapped to $\Bun_G$ by means of the projection $\hr_G$.

\smallskip

Let also $_{rel,\infty}\BunPbw\to X$ (resp., $_{rel,\geq\nu}\BunPbw$, $_{rel,\nu}\BunPbw$)
denote the relative version of the stack $_{x,\infty}\BunPbw$ (resp., $_{x,\geq\nu}\BunPbw$, $_{x,\nu}\BunPbw$)
introduced in \secref{infinity}. In other words, the fiber of $_{rel,\infty}\BunPbw$ over $x\in X$ 
is $_{x,\infty}\BunPbw$ and similarly for $_{rel,\geq\nu}\BunPbw$ and $_{rel,\nu}\BunPbw$.

\smallskip

As in \secref{variantG}, there exists a second projection 
$'\hl_G:\overline Z\to{}_{rel,\infty}\BunPbw$:

$$
\CD
_{rel,\infty}\BunPbw @<{'\hl_G}<< \overline Z @>{'\hr_G\times\pi}>>  
\BunPbw\times X\\
@VVV   @VVV     @VVV  \\
\Bun_G\times X  @<{\hl_G\times\pi}<<  \H_G^{\lambda}  @>{\hr_G\times\pi}>> 
\Bun_G\times X.
\endCD
$$

\medskip

Let $Z\subset \overline Z$ be the following locally closed substack:
$$Z:={}'\hl_G{}^{-1}({}_{rel,\lambda}\BunPbw)\cap 
{}'\hr_G{}^{-1}({}_{rel,0}\BunPbw)$$
($_{rel,0}\BunPbw$ is, according to our conventions, an open substack of 
$\BunPbw\times X$.)

Thus, $Z$ classifies the data of 
$$z=(x,\F_G,\F_M,\wt\kappa_P,\F'_G,\F'_M,\wt\kappa'_P,\beta,\beta_M),$$
where $(x,\F_G,\F_M,\wt\kappa_P)\in{}_{x,0}\BunPbw$, 
$(x,\F'_G,\F'_M,\wt\kappa'_P)\in{}_{x,0}\BunPbw$,
$(x,\F_G,\F'_G,\beta)\in\H^\lambda_G$, $(x,\F_M,\F'_M,\beta_M)\in\H^\lambda_M$ and
the data of $\beta,\beta_M,\wt\kappa_P,\wt\kappa'_P$ are compatible in the sense that over $X-x$, 
$\beta\circ\wt\kappa_P=\wt\kappa'_P\circ\beta_M$.

We have the projections
$$Z\to {}_{rel,0}\BunPbw\underset{\Bun_M\times X}\times \H_M^{\lambda}
\text{ and } Z\to {}_{rel,0}\BunPbw\underset{\Bun_M\times X}\times \H_M^{-\w^M_0(\lambda)},$$
that send a point $z\in Z$ as above 
to $$((x,\F_G,\F_M,\wt\kappa_P),\F'_M,\beta_M) \text{ and }
((x,\F'_G,\F'_M,\wt\kappa'_P),\F_M,\beta^{-1}_M),$$ respectively.
 
If follows from \secref{fulldescrG} that the first of the above projections is smooth, since the map 
$\tf_P^{\lambda}:S_P^{\lambda}\cap \Gr_G^{\lambda}\to \Gr_M^{\lambda}$ is smooth and that the second projection
is an isomorphism, since the map 
$\tf_P^{-\w^M_0(\lambda)}:S_P^{-\w^M_0(\lambda)}\cap \Gr_G^{-\w_0(\lambda)}\to \Gr_M^{-\w^M_0(\lambda)}$
is an isomorphism.

Let us denote by $\varphi_1$ and $\varphi_2$ the maps $Z\to\Bun_M$
that send $z\in Z$ as above to $\F_M$ and $\F'_M$, respectively. Let $\overset{o}Z$ be 
the preimage of $\overset{o}{\Bun}_M$ under $\varphi_1$.

\smallskip
 
\begin{lem}
{\it $\IC_{\overset{o}Z}$ is ULA with respect to $\varphi_2:\overset{o}Z\to\Bun_M$.}
\end{lem}

\begin{proof}

According to Property 2 of \secref{propacyc},
$\IC_{\overset{o}Z}$ is ULA with respect to the map \newline
$\overset{o}Z\to {}_{rel,0}\BunPbw\underset{\Bun_M\times X}\times \H_M^{\lambda}\to \H_M^{\lambda}$. 

Now, $\varphi_2$ is the composition
$\overset{o}Z\to \H_M^\lambda\overset{\hr_M}\longrightarrow\Bun_M$
and the assertion follows from Property 6 of \secref{propacyc}.

\end{proof}

Consider the map
$$\overset{o}Z\hookrightarrow\ol{Z}\overset{'\hr_G\times\pi}\longrightarrow {}_{rel,0}\BunPbw\to\BunPbw.$$
It is smooth and its composition with $\qw_P:\BunPbw\to\Bun_M$ equals $\varphi_2$. Therefore, 
in order to prove that $\IC_{\BunPbw}$ is ULA with respect to $\qw_P$ in a neighbourhood of 
$\F'_M\in\Bun_M(\Fqb)$, it is enough to show that
the preimage $\qw_P^{-1}(\F'_M)\subset \BunPbw(\Fqb)$ is contained in the image 
of $\overset{o}Z$ in $\BunPbw$
under the above map.

Thus, let 
$(\F'_G,\F'_M,\{\widetilde\kappa'_P{}^{\V}:
(\V^{U(P)})_{\F'_M}\hookrightarrow \V_{\F'_G}\})$
be a $\ol{\Fq}$-point of $\qw_P^{-1}(\F'_M)$.
As the map $\H_M^\lambda\overset{\hr_M}\longrightarrow 
\Bun_M\times X$ is smooth and, in particular, open, the fact 
that $\F_M\prec\F'_M$ for some $\F_M\in \overset{o}{\Bun}_M(\Fqb)$,
implies that there exists a triple $(\F''_M,x,\beta_M)$ such that

a) $(\F''_M,\F'_M,x,\beta_M)\in \H^\lambda_M(\Fqb)$.

b) The embeddings 
$\wt\kappa'_P{}^{\lambdach}:
(\V^{(U(P)})_{\F'_M}\hookrightarrow 
\V_{\F'_G}$ have no zero at $x$.

c) $\F''_M\in \overset{o}{\Bun}_M(\Fqb)$.

\smallskip

But this exactly means that the point
$$(x,\F'_G,\F'_M,\widetilde\kappa'_P,\F_M'',\beta^{-1}_M)\in 
{}_{rel,0}\BunPbw\underset{\Bun_M\times X}\times 
\H_M^{-\w^M_0(\lambda)}\simeq Z$$ belongs to $\overset{o}Z$.

\medskip

The "only if" part of \propref{eqacyc} follows from similar considerations by interchanging
right and left.

\end{proof}

\section{The structure of Drinfeld's compactifications}  

\ssec{Stratifications-I}  \label{stra1}

\sssec{}  \label{strataofBunBb}

Consider the set $\syminfty(\Lambda)$ whose elements are 
finite unordered collections of elements of $\Lambda$ with possible repetitions. For
$$\ol{\lambda}=\{\underset{n_1\text{ times}}{\underbrace{\lambda_1,...,\lambda_1}},...,
\underset{n_k\text{ times}}{\underbrace{\lambda_k,...,\lambda_k}}\}\in\syminfty(\Lambda),$$
we define the corresponding partially symmetrized power of $X$ with all the diagonals
removed as 
$$X^{\ol\lambda}:=X^{(n_1)}\times...\times X^{(n_k)}-\Delta.$$ 
To a point $\ol{x}\in X^{\ol\lambda}$, where $\ol{x}=x_{1,1},...,x_{1,n_1},x_{2,1},...,x_{k,n_k}$
we will attach a $\Lambda$--valued divisor $\ol{\lambda}\cdot\ol{x}$ equal to
$\Sum\, \lambda_i\cdot x_{i,j}$. We will denote by $|\ol{\lambda}|$ the element of $\Lambda$ equal to 
$\Sum\, n_i\cdot\lambda_i$.

\smallskip

Let $\ol{\lambda}$ belong to $\syminfty(\Lambda_G^{\on{pos}}-0)$. Consider the map
$$j_{\ol\lambda}:\Bun_B\times X^{\ol\lambda}\to\BunBb,$$ which sends a point
$(\F_G,\F_T,\kappa)\times \ol{x}$ to $(\F_G,\F'_T,\kappa')$, where
$\F'_T:=\F_T(-\ol{\lambda}\cdot\ol{x})$ and $\kappa'$ corresponds to the composition:
$$\L^{\lambdach}_{\F'_T}\hookrightarrow \L^{\lambdach}_{\F_T}
\overset{\kappa^{\lambdach}}\longrightarrow \V^{\lambdach}_{\F_G}.$$

\begin{prop}
For $\ol\lambda\in \syminfty(\Lambda_G^{\on{pos}}-0)$, the map $j_{\ol\lambda}$ 
a locally closed embedding.
\end{prop}

\begin{proof}

The fact that $j_{\ol\lambda}$ is representable is evident, since both the source and the target
are representable over $\Bun_G$.

Hence, for an $S$-point of $\BunBb$, its preimage under $j_{\ol\lambda}$ is a {\it set}, rather
than a category, and to prove that $j_{\ol\lambda}$ is a locally closed embedding, one has to show
that this set consists of at most one element. The latter is, however, obvious from the definitions.

\end{proof}

Let $_{\ol\lambda}\Bun_B$ denote the locally closed substack of $\BunBb$ equal to the image of 
$j_{\ol\lambda}$. For a fixed point $\ol{x}\in X^{\ol\lambda}$, let $_{\ol{x},\ol{\lambda}}\Bun_B$ 
denote the image of $\Bun_B\times \ol{x}$ under $j_{\ol{\lambda}}$.

\medskip

Therefore, by combining the above proposition with \propref{pointsofBunBb} we obtain:

\begin{prop}  \label{stra}
The locally closed substacks
$_{\ol{\lambda}}\Bun_B$ form a stratification of $\BunBb$ as $\ol\lambda$ runs over 
$\syminfty(\Lambda_G^{\on{pos}}-0)$.
\end{prop}

\sssec{} \label{symgr}

Now let us give a proof of \propref{denseprince} which was announced in \secref{intrcompact}. For that end
we need to introduce one more piece of notation. For $\lal\in \syminfty(\Lambda_G^+)$, let
$\H_G^{\lal}$ denote the relative version of the Hecke stack:

By definition, $\H_G^{\lal}$ is endowed with a map $\pi: \H_G^{\lal}\to X^{\lal}$ and its fiber over
$\xl=x_{1,1},...,x_{1,n_1},x_{2,1},...,x_{k,n_k}$ is the stack of triples $(\F_G,\F'_G,\beta)$, where
$\F_G$ and $\F'_G$ are $G$--bundles on $X$ and $\beta$ is an identification between them on
$X-\{\xl\}$ such that $\F'_G$ is in position $\lambda_i$ with respect to $\F_G$ at $x_{i,j}$. 

Let $\hl_G$ and $\hr_G$ be the projections $\H_G^{\lal}\to\Bun_G$ that send $(\xl,\F_G,\F'_G,\beta)$ to
$\F_G$ and $\F'_G$, respectively. We will denote by $_{\xl}\H_G^{\lal}$ the fiber of $\H_G^{\lal}$ over
$\xl\in X^{\lal}$. 

In a similar way we define the stacks $\Hb_G^{\lal}$ and $_{\xl}\Hb_G^{\lal}$.

\begin{proof}(of \propref{denseprince})

Let $(\F_G,\F_T,\kappa)$ be an $\Fqb$-point of $\BunBb$. To prove the proposition, it is enough to construct
an irreducible stack $Z$ with a map $\phi:Z\to \BunBb$ such that

\smallskip

\noindent{\em (a)}
$\on{Im}(\phi)\cap \Bun_B$ is non-empty.

\smallskip

\noindent{\em (b)}
$(\F_G,\F_T,\kappa)\in \on{Im}(\phi)$.

\medskip

Let $x_1,...x_n$ be the set of points where $(\F_G,\F_T,\kappa)$ has singularities and let
$\nu_1,...,\nu_n\in\Lambda_G^{\on{pos}}$ be the corresponding defects. Let us choose elements
$\lambda_1,...,\lambda_n\in\Lambda^+_G$ in such a way that the weight spaces 
$V^{\lambda_i}(\w_0(\lambda_i)+\nu_i)$ are nonzero. The collection $\{(x_1,\lambda_1)...,(x_n,\lambda_n)\}$
corresponds to a unique $\ol{\lambda}\in \syminfty(\Lambda_G^+)$ and $\ol{x}\in X^{\lal}$.

We define $Z'$ as a fiber product $Z':={}_{\xl}\H_G^{\lal}\underset{\Bun_G}\times\Bun_B$, where
$_{\xl}\H_G^{\lal}$ is mapped to $\Bun_G$ by means of $\hr_G$.

The stack $Z'$ splits into connected components (which are in a bijection with the connected components of $\Bun_B$,
and hence of $\Bun_T$) and we take $Z$ to be the connected component of $Z'$ corresponding to the coweight 
$\on{deg}(\F_T)-\w_0(|\lal|)$.

\smallskip

We define the map $\phi: Z'\to \BunBb$ as in \secref{bd}. Now, it follows from \lemref{descrprince}
that $Z$ satisfies conditions (a) and (b) above.

\end{proof}

\sssec{}  

For an element $\lambda\in\Lambda_G^{\on{pos}}$ let us consider the following (semi-simple) 
complex over $\on{Spec}(\Fq)$:
$$\on{Kost}_\lambda:=\underset{\lambda=\Sigma \,m_{\alpha}\cdot \alpha,\, \alpha\in\Delta^+}\oplus 
(\Ql(1)[2])^{\otimes \Sigma \,m_{\alpha}},$$
where the sum is taken over all the possible ways to represent $\lambda$ as a sum of positive coroots with non-negative
coefficients. (This is a geometric counterpart of the $q$-analog of Kostant's partition function (cf. \cite{Lu}).)

\begin{thm} \label{compute1}
For $\lal=\{\underset{n_1\text{ times}}{\underbrace{\lambda_1,...,\lambda_1}},...,
\underset{n_k\text{ times}}{\underbrace{\lambda_k,...,\lambda_k}}\}
\in\syminfty(\Lambda^{\on{pos}}-0)$, the complex \newline
$j_{\lal}^*(\IC_{\BunBb})$ is isomorphic to 
$(\IC_{\Bun_B}\boxtimes {\Ql}_{X^{\lal}})\underset{i=1}{\overset{k}\otimes}
\on{Kost}_{\lambda_i}^{\otimes n_i}$.

\end{thm}

This result (along with its generalization, \thmref{compute3}) will be proven in a subsequent
publication. Now, let us show how \thmref{compute1} implies \thmref{comparefinite} (and hence \thmref{fullcompare}):

\smallskip

\begin{proof}(of \thmref{comparefinite})

For every collection of non-negative integers $\{m_\alpha\}$, $\alpha\in\Delta^+$, 
consider the corresponding map
$$i_{\{m_\alpha\}}:\Bun_B\times \underset{\alpha\in\Delta^+}\Pi X^{(m_{\alpha})}\to\BunBb.$$

\thmref{compute1} implies that the function on $\BunBb(\Fq)$ corresponding to the sheaf 
$\q^{!*}(\Aut^\mu_{E_{\check T}})$ equals 
$$\underset{\{m_{\alpha}\}}\Sigma\,
i_{\{m_\alpha\}}{}_{!}(\on{Funct}({\mathfrak q}^{!*}
(\Aut^{\mu-\Sigma m_{\alpha}\cdot\alpha}_{E_{\check T}}))\boxtimes
\underset{\alpha\in\Delta^+}\Pi\, \on{Funct}((E^{\alpha}_{\check T})^{(m_\alpha)}\otimes\Ql(m_{\alpha}))).$$

\smallskip

By summing up along the fibers of the projection $\p$, we derive the formula of \thmref{comparefinite}.

\end{proof}

\ssec{Stratifications-II}   \label{stra2}

\sssec{}  \label{stratBunP}

Let us generalize the above discussion to the case of a parabolic subgroup $P$.

Conisider the set $\syminfty(\Lambda_{G,P}^{\on{pos}}-0)$
and for every element $\thl=\{\underset{n_1\text{ times}}{\underbrace{\theta_1,...,\theta_1}},...,
\underset{n_k\text{ times}}{\underbrace{\theta_k,...,\theta_k}}\}$ in it, consider the corresponding
variety $X^{\thl}=X^{(n_1)}\times...\times X^{(n_k)}-\Delta$

As in \secref{strataofBunBb}, we have the natural locally closed embeddings
$$j_{\thl}:\Bun_P\times X^{\thl}\hookrightarrow \BunPb,$$
and for a fixed $\thl\in \syminfty(\Lambda_{G,P}^{\on{pos}}-0)$ (resp., $\xl\in X^{\thl}$)
we will denote by $_{\thl}\Bun_P$ (resp., $_{\xl,\thl}\Bun_P$) the corresponding stratum (resp., locally
closed subset) of $\BunPb$.

Thus, for $\thl$ and $\xl$ as above we obtain the locally closed substacks $\r_P^{-1}({}_{\thl}\Bun_P)$
and $\r_P^{-1}({}_{\xl,\thl}\Bun_P)$ of $\BunPbw$. However, it will not be true that $\IC_{\BunPbw}$
is smooth when restricted to the strata of the form $\r_P^{-1}({}_{\thl}\Bun_P)$. Our next goal is 
to define a suitable refinement of this stratification.

\sssec{} \label{pospart}

Consider the affine Grassmannian of the group $M$ and let $\Gr_M^+\subset\Gr_M$ be a closed subscheme definied
by the following condition:

\smallskip

$\F_M\in\Gr_M^+$ if for every $G$-module $\V$, the map
$$\beta^{\V^{(U(P)}}_M: (\V^{(U(P)})_{\F_M}|_{\D^*_x}\to 
(\V^{(U(P)})_{\F^0_M}|_{\D^*_x}$$
is regular on $\D_x$.

It is clear from the definition, that $\Gr_M^+$ is stable under the 
$M(\O_x)$-action on $\Gr_M$.

\begin{prop}

For an $M$ co-weight $\nu\in\Lambda_M^{+}$ the following conditions are equivalent:

\smallskip

\noindent{\em (a)} $\grb_M^\nu\subset \Gr_M^+$.

\smallskip

\noindent{\em (b)} $\on{w}_0^M(\nu)\in \Lambda_G^{\on{pos}}$.

\end{prop}

\begin{proof}

Assume that $\grb_M^\nu\subset \Gr_M^+$. In particular, the corresponding point 
$\F_M:=t_x^{\nu}\in M(\K_x)/M(\O_x)=\Gr_M$ belongs to $\Gr_M^+$. By definition, we must have that
for $\lambdach\in\Lambda^+_G$ the map 
$$\U^{\lambdach}_{\F_M}|_{\D^*_x}\to \U^{\lambdach}_{\F^0_M}|_{\D^*_x}$$ is regular on $\D_x$,
which implies that $\langle \nu, \w_0^M(\lambdach)\rangle \geq 0$. The latter inequality means exactly
that $\w_0^M(\nu)\in \Lambda_G^{\on{pos}}$.

Conversely, let us assume that $\on{w}_0^M(\nu)\in \Lambda_G^{\on{pos}}$.
Let $\U$ be an $M$-module of the form $\V^{U(P)}$, where $\V$ is a $G$-module. Without restricting
the generality, we way assume that the weights of $\V$ are $\leq \lambdach$ for some $\lambdach\in\Lambda^+_G$.
Then the weights of $\U$ are $\underset{M}\leq \lambdach$. Hence, by definition, $\F_M\in \grb_M^{\nu}$ means that
$$\beta^{\U}_M:\U_{\F_M}\to \U_{\F^0_M}(-\langle \w_0^M(\nu),\lambdach\rangle).$$
However, $\langle \w_0^M(\nu),\lambdach\rangle\geq 0$, by assumption.

\end{proof}

In what follows we will denote by $\Lambda_{M,G}^{+}\subset\Lambda_M^{+}$ 
the above sub-semigroup, i.e. 
$\Lambda_{M,G}^{+}=\Lambda_M^{+}\cap \w_0^M(\Lambda_G^{\on{pos}})$.
For example, for $M=T$, $\Lambda_{T,G}^{+}=\Lambda_G^{\on{pos}}$. In general, the projection
$\Lambda\to\Lambda_{G,P}$ sends $\Lambda_{M,G}^{+}$ to $\Lambda_{G,P}^{\on{pos}}$.

\medskip

Let now $\theta$ be an element of $\Lambda_{G,P}^{\on{pos}}$. We define the element 
$\flat(\theta)\in \Lambda_{M,G}^{+}$ as follows: 

\smallskip

\noindent By definition, $\theta$ is the projection under $\Lambda\to\Lambda_{G,P}$ of some 
$\underset{\i\in\I-\I_M}\Sigma \, b_\i\cdot \alpha_\i$. We set $\flat(\theta)=\w_0^M(\underset{\i\in\I-\I_M}\Sigma \,
b_\i\cdot \alpha_\i)$. It belongs to $\w_0^M(\Lambda_G^{\on{pos}})$ by construction, and to $\Lambda_M^{+}$, because
every $\alpha_\i,\i\notin \I_M$ is $M$-antidominant.

For $\theta$ as above let $\Gr_M^\theta$ denote the corresponding connected
component of $\Gr_M$. Set $\Gr_M^{+,\theta}:=\Gr_M^+\cap \Gr_M^\theta$.

From the above proposition it follows that $\flat(\theta)$ is the maximal element
in the set of $\nu\in \Lambda^+_M$ such that $\grb_M^\nu\subset \Gr_M^+\cap \Gr_M^\theta$. Hence, we obtain:
\begin{equation}
(\Gr_M^{+,\theta})_{red}=(\ol{\Gr}_M^{\flat(\theta)})_{red},.
\end{equation}
where the subscript ``$red$'' means ``the corresponding reduced scheme''.

\sssec{}  \label{hereBunPbw}

Fix $\thl=\{\underset{n_1\text{ times}}
{\underbrace{\theta_1,...,\theta_1}},...,
\underset{n_k\text{ times}}{\underbrace{\theta_k,...,\theta_k}}\}\in\syminfty(\Lambda_{G,P}^{\on{pos}}-0)$.
Let us denote by $\H_M^{+,\thl}$ a version of the Hecke stack, which is fibered over $X^{\thl}$
with the fiber over a point $\xl\in X^{\thl}$ given by $x_{1,1},...,x_{1,n_1},x_{2,1},...,x_{k,n_k}$ being the product
$$_{x_{1,1}}\H_M^{+,\theta_1}\underset{\Bun_M}\times...\underset{\Bun_M}\times {}_{x_{k,n_k}}\H_M^{+,\theta_k}.$$

\begin{prop}  \label{BunPbwfibered}
There is a canonical isomorphism 
$\H_M^{+,\thl}\underset{\Bun_M}\times\Bun_P\simeq \r_p^{-1}({}_{\thl}\Bun_P)$
that fits into the commutative diagram
$$
\CD
\Bun_P\underset{\Bun_M}\times \H_M^{+,\thl} @>{\sim}>>    \r_P^{-1}({}_{\thl}\Bun_P)\\
@V{\hl_M}VV     @V{\r_P}VV    \\
\Bun_P\times X^{\thl} @>{j_{\thl}}>>   _{\thl}\Bun_P,
\endCD
$$
where we have used the projection $\hl_M:\H_M^{+,\thl}\to\Bun_M$ to define the fiber product.
\end{prop} 

\begin{proof}

Let $(\xl,\F_P,\F'_M,\beta_M)$ be a point of $\Bun_P\underset{\Bun_M}\times \H_M^{+,\thl}$,
where $\F_P$ is a $P$--bundle on $X$ and $\beta_M$ is an isomorphism 
between $\F_M:=U(P)\backslash \F_P$
and $\F'_M$ defined outside of $\{\xl\}\subset X$.

We attach to it a point of $\r_P^{-1}({}_{\thl}\Bun_P)$ as follows: 

\smallskip

\noindent The corresponding $G$-bundle is induced from $\F_P$, i.e. $\F_G:=\F_P\overset{P}\times G$ and 
the corresponding $M$-bundle is $\F'_M$. Now, the data of 
$\widetilde\kappa'_P{}^{\V}: 
(\V^{U(P)})_{\F'_M}\to \V_{\F_G}$ is obtained as a composition

$$(\V^{U(P)})_{\F'_M} \overset{\beta^{\V^{U(P)}}_M}\hookrightarrow 
(\V^{U(P)})_{\F_M}\to\V_{\F_G}.$$

Conversely, let $(\F_G,\F'_M,\wt\kap_P)$ be an ($S$)-point of $\r_P^{-1}({}_{\thl}\Bun_P)$. 
Let $(\xl,\F_P)$ be the corresponding point of $\Bun_P\times X^{\thl}$, in particular,
let $\F_M$ be the corresponding $M$-bundle. By definition, $\F_M$ and $\F'_M$ are identified
outside $\xl$.

By assumption, for each $G$-module $\V$, we have the embedding:
$$\wt\kap'{}^{\V}_P:(\V^{U(P)})_{\F'_M} \hookrightarrow \V_{\F_G}.$$
and the {\it maximal} embedding
$$\wt\kap^{\V}_P: (\V^{U(P)})_{\F_M} \hookrightarrow \V_{\F_G}.$$
Hence, the (a priori) meromorphic map map $(\V^{U(P)})_{\F'_M}\to (\V^{U(P)})_{\F_M}$ is regular.

\end{proof}

\sssec{}

If $\nul$ is an element of $\syminfty(\Lambda_{M,G}^+-0)$, we obtain a locally closed embedding
$$j_{\nul}: \Bun_P\underset{\Bun_M}\times\H_M^{\nul}\hookrightarrow \BunPbw.$$
Its image will be denoted by $_{\nul}\Bun_P$, and for $\xl\in X^{\nul}$ we let 
$_{\xl,\nul}\Bun_P$ denote the corresponding closed substack of $_{\nul}\Bun_P$.

Thus, we see that the defect of a point $(\F_G,\F'_M,\wt\kappa_P)\in\BunPbw(\Fqb)$ at a
point of $X$ where it has a singularity is naturally an element of $\Lambda^+_{M,G}$.
However, unlike the case of $\BunBb$ and $\BunPb$, fixing the locus of singularities of an enhanced parabolic 
Drinfeld's structure, together with the defects and the saturated $P$-bundle, does not determine the
corresponding point of $\BunPbw$ uniquely: we have the remaining freedom of choosing a point in the product
of the corresponding $\grb_M^{+,\theta}$'s.

At this point we are ready to give a proof of \propref{densenonprince}. 

\begin{proof} (of \propref{densenonprince})

As we have seen above, the fibers of the map $\r_P:\BunPbw\to\BunPb$ are never empty. Therefore,
it is enough to prove that $\Bun_P$ is dense in $\BunPbw$. Let $(\F_G,\F_M,\wt\kap_P)$ be an $\Fqb$-point
of $\BunPbw$.

As in the case of \propref{denseprince}, it is enough to construct an irreducible
stack $Z$ with a map $\phi:Z\to \BunPbw$, whose image contains $(\F_G,\F_M,\wt\kap_P)$ and such that
$\on{Im}(\phi)\cap \Bun_P\neq \emptyset$. 

To simplify the notation, we will assume that $(\F_G,\F_M,\wt\kap_P)$ has a singularity only at one point,
call it $x$. Let the defect be $\nu\in \Lambda^+_{M,G}$. We have: 
$(\F_G,\F_M,\wt\kap_P)\subset {}_{x,\nu}\BunPbw$, in the terminology of \secref{stratinf}.

\smallskip

Let $\lambda\in\Lambda^+_G$ be such that: (a) $\langle \w_0(\lambda),\alphach_\i\rangle=0$ for $\i\in\I_M$;
(b) $\Hom_{\check M}(U^{\w_0(\lambda)+\nu},V^\lambda)\neq 0$. It is easy to see that such 
$\lambda$ indeed exists (cf. \lemref{soccer}).
Note that condition (a) means that $\w_0(\lambda)$ lies in the group of cocharacters
of $Z^0(M)$.

\medskip

Set $Z':={}_x\H_G^\lambda\underset{\Bun_G}\times\Bun_P$, where $_x\H_G^\lambda$ is mapped to $\Bun_G$
by means of $\hr_G$. We define the map $\varphi:Z'\to\BunPbw$ as follows:

\smallskip

For a point $(\F_G,\beta,\F'_G,\F'_M,\wt\kap'_P)\in Z'$, the resulting $G$-bundle is $\F_G$, the $M$-bundle
is $\F'_M\otimes \F^0_{Z^0(M)}(\w_0(\lambda)\cdot x)$ and the new $\wt\kap_P$ are obtained
from $\wt\kap'_P$ as in \secref{bd}.

We define $Z$ as the preimage in $Z'$ of the appropriate connected component of $\Bun_P$. The fact that
$Z$ satisfies the required properties follows from \secref{fulldescrG}.

\end{proof}

\sssec{}

Here we will describe the behaviour of $\IC_{\BunPbw}$ along the strata $_{\nul}\Bun_P$
introduced above.

Recall (cf. \thmref{dualgroup}) that to every representaion $U$ of the group $\check M$ we can attach 
in a canonical way a $M(\O_x)$-equivariant perverse sheaf $F_M^{-1}(U)$ on the affine Grassmannian $\Gr_M$. 
For $\theta\in\Lambda_{G,P}^{\on{pos}}$ we introduce the following semi--simple complex on $\Gr_M$:
$$\on{Kost}^P_{\theta}:=\underset{i\geq 0}\oplus\, F_M^{-1}(\on{Sym}^i(\check 
{\mathfrak u}_P)_{\theta})\otimes\Ql(i)[2i],$$
where $\check {\mathfrak u}_P$ denotes the unipotent radical of the corresponding parabolic in $\check G$ and
where for an $\check M$--representation $U$, $U_\theta$ denotes its piece on which $Z(\check M)$ acts by the
character $\theta$.

\begin{lem}  \label{soccer}
For an $M$-dominant coweight $\nu$,  
$\on{Hom}_{\check M}(U^\nu, \on{Sym}^i(\check {\mathfrak u}_P)_\theta)\neq 0$
if and only if $\nu\in \Lambda_{M,G}^{+}$.
Moreover, the above space is non-zero only for finitely many integers $i$.
\end{lem}

Therefore, $\on{Kost}^P_{\theta}$ is supported on $\grb_M^{\flat(\theta)}\subset \Gr_M^{+,\theta}$. 
Analogously, for $\thl\in \syminfty(\Lambda_{G,P}^{\on{pos}}-0)$ we define the semi--simple
complex $\on{Kost}^P_{\thl}$ on $\H_M^{+,\thl}$. 

\smallskip

\begin{thm} \label{compute3}
Under the identification 
$$\Bun_P\underset{\Bun_M}\times\H_M^{+,\thl}$$
the sheaf $\IC_{\BunPbw}|_{_{\thl}\BunPbw}$ goes over to 
$\IC_{\Bun_P}\boxtimes\on{Kost}^P_{\thl}$.
\end{thm}

This theorem will be neither used nor proven in this paper (the proof will appear in \cite{BFGM}). 
For our purposes the following weaker result will be sufficient:

\smallskip

\begin{thm}  \label{smoothonstrata}
For fixed $\nul\in\syminfty(\Lambda_{M,G}^+-0)$ and $\xl\in X^{\nul}$, the $*$--restriction
of $\IC_{\BunPbw}$ to $_{\xl,\nul}\Bun_P$ is locally constant.
\end{thm}

The proof will occupy the next subsection. The reader will notice that it is very
similar to the proof of \thmref{acycthm}.

\ssec{Proof of \thmref{smoothonstrata}}

\sssec{}
Consider the fiber product
$\ol{Z}={}_y\H_G^\lambda\underset{\Bun_G}\times \Bun_P$, where
$y$ is a point in $X-\{\xl\}$, $\lambda$ is some element of $\Lambda_G^+$
and the fibered product is formed using the projection $\hr_G:{}_y\H_G^\lambda\to\Bun_G$.

As before, there is a projection $'\hl_G:\ol{Z}\to {}_{y,\infty}\BunPbw$ and we will denote by
$Z$ the locally closed substack of $\ol{Z}$ defined as
$$Z:={}'\hl_G{}^{-1}({}_{y,\lambda}\BunPbw)\cap{}'\hr_G{}^{-1}({}_{y,0}\BunPbw).$$

According to \secref{stratinf}, there is an isomorphism 
$_{y,\lambda}\BunPbw\simeq {}_{y,0}\BunPbw\underset{\Bun_M}\times{}_y\H_M^\lambda$. Hence we obtain
a commutative diagram:
$$
\CD
_{y,0}\BunPbw  @<{\varphi}<< Z  @>{\varphi^1}>> _{y,0}\BunPbw  \\
@VVV        @V{\psi}VV     @VVV    \\
\Bun_M  @<{\hl_M}<<   _y\H_M^{\lambda}   @>{\hr_M}>>   \Bun_M,
\endCD
$$
where both upper horizontal arrows are smooth. Moreover, the preimages of 
$_{\xl,\nul}\Bun_P$ in $Z$ under the maps $\varphi$ and $\varphi^1$ coincide.

\smallskip

Consider two $\overline{\Fq}$-points of the stack $_{\xl,\nul}\Bun_P$. Using 
\propref{BunPbwfibered}, they can be represented 
by two triples $(\F_P,\F'_M,\beta)$ and $(\F^1_P,{\F_M^1}',\beta^1)$, where $\beta$ 
is an isomorphism 
$\F_M|_{X-\{\xl\}}\to \F'_M|_{X-\{\xl\}}$ (here $\F_M$ is the $M$-bundle induced from $\F_P$),
such that at every $x_{i,j}$, $\F'_M$ is in position $\nu_i$ with respect to $\F_M$;
and similarly for $\beta^1$.

We will write that $(\F_P,\F'_M,\beta)\succ (\F^1_P,{\F_M^1}',\beta^1)$ if there exist $y\in X-\{\xl\}$,
$\lambda\in\Lambda_G^{+}$ and $z\in Z$ as above such that
$(\F_P,\F'_M,\beta)=\varphi(z),\, (\F^1_P,{\F^1_M}',\beta^1)=\varphi^1(z)$.
Consider the equivalence relation $\sim$ on the set of $\overline{\Fq}$-points of $_{\xl,\nul}\Bun_P$
generated by $\succ$.

We will deduce the assertion of \thmref{smoothonstrata} from the following proposition:

\begin{prop} \label{propsmoothonstrata}
All $\overline{\Fq}$-points of $_{\xl,\nul}\Bun_P$ are $\sim$--equivalent.
\end{prop}

\begin{proof}

Set $\widetilde{M}:=M/Z^0(M)$ and for a an $M$--bundle, we will denote by a subscript 
$\widetilde{M}$ the corresponding induced $\widetilde{M}$--bundle.

\smallskip

\noindent{\bf Step 1.} Let $(\F_P,\F'_M,\beta)$ and $(\F^1_P,{\F^1_M}',\beta^1)$ be two $\ol{\Fq}$--points of 
$_{\xl,\nul}\Bun_P$. First we will show that $(\F^1_P,{\F^1_M}',\beta^1)$ is $\sim$--equivalent 
to another $\ol{\Fq}$--point $(\F^2_P,{\F^2_M}',\beta^2)$, such that 
$\F_{\widetilde{M}}\simeq\F^2_{\widetilde{M}}$ and the induced isomorphism
$${\F_{\widetilde{M}}}'|_{X-\{\xl\}}\simeq {\F^2_{\widetilde{M}}}'|_{X-\{\xl\}}$$
is regular on the whole of $X$.

\smallskip

Pick $y\notin \{\xl\}$. As $\widetilde M$ is semi-simple, the $\widetilde M$-bundles
$\F_{\widetilde{M}}$ and $\F^1_{\widetilde{M}}$ are isomorphic over $X-y$ (cf. \lemref{affinetrivial}).
Moreover, since $X-y$ is affine,
the group of automorphisms of $\F_{\widetilde{M}}$ over $X-y$ is dense in $\underset{i,j}\Pi\, 
G(\O_{x_{i,j}})(\ol{\Fq})$. Therefore, since at every $x_{i,j}$ the relative position of $\F'_M$ with respect 
to $\F_M$ is the same as the relative position of ${\F^1_M}'$ with respect to $\F^1_M$, the above 
isomorphism $\F_{\widetilde{M}}|_{X-y}\to \F^1_{\widetilde{M}}|_{X-y}$ can be chosen in such a way that
the isomorphism between ${\F_{\widetilde{M}}}'|_{X-\{\xl,y\}}$ and ${\F^1_{\widetilde{M}}}'|_{X-\{\xl,y\}}$
is regular on $X-y$.

\smallskip

Consider the corresponding point $(\F_{\widetilde M},\F^1_{\widetilde M},\beta_{\wt M})\in {}_y\H^\nu_{\widetilde M}$.
Let $\lambda\in\Lambda_G^+$ be such that it projects onto $\nu$ under $\Lambda\to \Lambda_{\widetilde M}$
and consider the corresponding stack $Z$. The composition
$$Z\overset{\psi\times\varphi_1}\longrightarrow 
{}_y\H_M^{\lambda}\underset{\Bun_M}\times {}_{y,0}\BunPbw\simeq
{}_y\H_{\widetilde M}^{\nu}\underset{\Bun_{\wt M}}\times{}_{y,0}\BunPbw$$
is surjective and let $z\in Z$ be such that under the above composition it 
maps to the point 
$$(\F_{\widetilde M},(\F^1_P,{\F^1_M}',\beta^1),
\beta_{\wt M}:\F_{\widetilde M}|_{X-y}\simeq \F^1_{\widetilde M}|_{X-y})\in {}_y\H_{\widetilde
M}^{\nu}\underset{\Bun_{\wt M}}\times{}_{y,0}\BunPbw.$$ 
We define $(\F^2_P,{\F^2_M}',\beta^2)$ as the image of $z$ in
under the projection $\varphi$. By construction, it satisfies the required condition.

\smallskip

\noindent{\bf Step 2.} Using the assertion of {\bf Step 1}, we can assume that there exists an isomorphism 
$\F_{\widetilde M}\to \F^1_{\widetilde M}$ which gives
rise to an isomorphism $\F'_{\widetilde M}\to {\F^1_{\widetilde M}}'$.

For a $\Lambda_{Z^0(M)}$--divisor $D$ and an $M$--bundle $\F_M$, let $\F_M(D)$ denote the new
$M$--bundle $\F_M\otimes \F^0_{Z^0(M)}(D)$. 

It is clear that we can choose $\Lambda_{Z^0(M)}\cap \Lambda_G^{+}$-valued divisors $D$ and $D_1$ on $X-\{\xl\}$,
so that there exists an isomorphism $\F_M(D)\simeq \F^1_M(D_1)$. In this case, the
identification $\F'_M(D)|_{X-\{\xl\}}\simeq {\F^1_M}'(D_1)|_{X-\{\xl\}}$ is automatically regular on the whole
of $X$.

\smallskip

Therefore, by arguing as above, we can replace our two points 
$(\F_P,\F'_M,\beta)$ and $(\F^1_P,{\F^1_M}',\beta_1)$ 
by $\sim$--equivalent points (we will abuse the notation and denote the latter by the same characters) such that
there exists an isomorphism
$\F_M\to \F^1_M$, for which the induced meromorphic map $\F'_M|_{X-\{\xl\}}\to {\F^1_M}'|_{X-\{\xl\}}$ is a global
isomorphism.

\smallskip

Moreover, by replacing $D$ and $D_1$ by $D+D'$ and $D_1+D'$, where $D'\in \Lambda_{Z^0(M)}\cap \Lambda_G^{+}$
is sufficiently large, we can ensure that the cohomology
\begin{equation} \label{vanish}
H^1(X,\U_{\F_M})=0,
\end{equation}
for all irreducible $M$--modules $\U$ which appear in the Jordan--Holder series of 
the Lie algebra ${\mathfrak u}(P)$, viewed as a $P$--module via the adjoint action.

However, \eqref{vanish} implies that any two $P$-bundles, whose reductions modulo $U(P)$ are isomorphic
to $\F_M$, are necessarily isomorphic. This implies that our two points 
$(\F_P,\F'_M,\beta)$ and $(\F_P^1,{\F^1_M}',\beta_1)$ are isomorphic and, in particular, 
$\sim$--equivalent.

\end{proof}

\sssec{}

Now let us prove \thmref{smoothonstrata}:

\begin{proof} 

Let $(\F_P,\F_M,\beta)$ and $(\F^1_P,{\F^1_M}',\beta^1)$ be two $\Fqb$-points of $_{\xl,\nul}\Bun_P$.
Let us denote by $\mathi:\Spec(\ol{\Fq})\to {}_{\xl,\nul}\Bun_P$ (resp., $\mathi_1$) the corresponding map.
Consider the complexes $$\mathi^*(\IC_{\BunPbw}|_{_{\xl,\nul}\Bun_P}),
\mathi_1{}^*(\IC_{\BunPbw}|_{_{\xl,\nul}\Bun_P}),\mathi^{!}(\IC_{\BunPbw}|_{_{\xl,\nul}\Bun_P}),
\mathi_1{}^{!}(\IC_{\BunPbw}|_{_{\xl,\nul}\Bun_P})$$
over $\Spec(\Fqb)$. We will normalize (by making a cohomological shift and Tate's twist) 
in such a way that their highest (resp., lowest) cohomology is $\Ql$ in degree $0$.

To prove the theorem, it is enough to show that for any $(\F_P,\F_M,\beta)$ and $(\F^1_P,{\F^1_M}',\beta^1)$
as above,
\begin{align*}
&\mathi^*(\IC_{\BunPbw}|_{_{\xl,\nul}\Bun_P})\simeq \mathi_1{}^*(\IC_{\BunPbw}|_{_{\xl,\nul}\Bun_P}) \\
&\mathi^{!}(\IC_{\BunPbw}|_{_{\xl,\nul}\Bun_P})\simeq \mathi_1{}^{!}(\IC_{\BunPbw}|_{_{\xl,\nul}\Bun_P}).
\end{align*}

Using \propref{propsmoothonstrata}, we can assume that $(\F_P,\F'_M,\beta) \succ (\F^1_P,{\F_M^1}',\beta^1)$
and let $\mathi_z: \Spec(\ol{\Fq})\to Z$ be a point such that $\varphi\circ \mathi_z=\mathi$,
$\varphi^1\circ\mathi_z=\mathi_1$. Let $_{\xl,\nul}Z$ denote the (common) preimage of 
$_{\xl,\nul}\Bun_P$ in $Z$ under $\varphi$ or $\varphi^1$. 

Since both maps $\varphi$ and $\varphi_1$ are smooth we have:
\begin{align*}
&\mathi^*(\IC_{\BunPbw}|_{_{\xl,\nul}\Bun_P})\simeq \mathi_z^*(\IC_Z|_{_{\xl,\nul}Z})\simeq
\mathi_1{}^*(\IC_{\BunPbw}|_{_{\xl,\nul}\Bun_P}) \text{ and } \\
&\mathi^{!}(\IC_{\BunPbw}|_{_{\xl,\nul}\Bun_P})\simeq \mathi_z^{!}(\IC_Z|_{_{\xl,\nul}Z})\simeq
\mathi_1{}^{!}(\IC_{\BunPbw}|_{_{\xl,\nul}\Bun_P}),
\end{align*}
where $\mathi_z^*(?)$ and $\mathi_z^{!}(?)$ are also normalized in the above way.

This proves our assertion. 

\end{proof}

\propref{propsmoothonstrata} implies also the following strengthening of \thmref{smoothonstrata}:

\begin{cor} \label{corsmoothonstrata}
Let $K$ be an irreducible subquotient of an $i$-th perverse cohomology sheaf
$h^i(\IC_{\BunPbw}|_{_{\nul}\Bun_P})$. Then $K$ is constant along the fibers of the projection
$_{\nul}\Bun_P\to X^{\nul}$.
\end{cor}

\begin{proof}

It is enough to prove that for each $\xl\in X^{\nul}$, every irreducible subquotient of
$h^i(\IC_{\BunPbw}|_{_{\xl,\nul}\Bun_P})$ is a constant sheaf on
$_{\xl,\nul}\Bun_P$. Therefore, it suffices to show
that $\on{Funct}(h^i(\IC_{\BunPbw}|_{_{\xl,\nul}\Bun_P}))$ is a constant function on 
$_{\xl,\nul}\Bun_P({\mathbb F}_{q'})$ for all ${\mathbb F}_{q'}\supset\Fq$. 

For that purpose, consider two points
$\mathi',\mathi'_1: \Spec({\mathbb F}_{q'})\to {}_{\xl,\nul}\Bun_P$ and for an
extension ${\mathbb F}_{q''}\supset {\mathbb F}_{q'}$, let $\mathi''$ and $\mathi''_1$ denote the corresponding
${\mathbb F}_{q''}$--points of $_{\xl,\nul}\Bun_P$.

As in the above proof of \thmref{smoothonstrata}, we infer from
\propref{propsmoothonstrata} that when ${\mathbb F}_{q''}$ is sufficiently large, the sheaves
$$\mathi''{}^*(h^i(\IC_{\BunPbw}|_{_{\xl,\nul}\Bun_P})) \text{ and }
\mathi_1''{}^*(h^i(\IC_{\BunPbw}|_{_{\xl,\nul}\Bun_P}))$$ over $\Spec({\mathbb F}_{q''})$
are isomorphic. 

Since this is true for infinitely many $q''$, we obtain that
$$\on{Tr}(Fr, \mathi'{}^*(h^i(\IC_{\BunPbw}|_{_{\xl,\nul}\Bun_P})))=
\on{Tr}(Fr, \mathi_1'{}^*(h^i(\IC_{\BunPbw}|_{_{\xl,\nul}\Bun_P}))),$$
which is what we had to prove.

\end{proof}

\ssec{``Good'' perverse sheaves and \thmref{compat}} \label{good}

\sssec{}

In this subsection our goal is to prove \thmref{compat}, which was announced in \secref{nonprinceseries}.
We shall start with the following general observation:

\smallskip

\begin{prop} \label{novancycl}
For $\S\in\Bun_G$ and $H^\lambda_G$ as in \secref{intrhecke}, the sheaf $H^\lambda_G(\S)$ on $\Bun_G\times X$ is 
ULA with respect to the projection $\Bun_G\times X\to X$.
\end{prop}

\begin{proof}

We have:
$$H^\lambda_G(\S):=(\hl_G\times\pi)_{!}(\hr_G{}^*(\S)\otimes \IC_{\Hb^\lambda_G})
\otimes(\Ql(\frac{1}{2})[1])^{-\on{dim}(\Bun_G)}.$$

As the projection $\hl_G\times\pi:\Hb^\lambda\to \Bun_G\times X$ is proper, 
it follows from Property 3 of \secref{propacyc}
that it is enough to show that the sheaf 
$\hr_G{}^*(\S)\otimes \IC_{\Hb^\lambda_G}$ on 
$\Hb_G^\lambda$ is ULA with repect to the projection $\pi:\Hb_G^\lambda\to X$.

\smallskip

However, since $\Hb_G^\lambda\overset{\hr_G\times\pi}\to\Bun_G\times X$ is a fibration 
(locally trivial in the smooth topology),
the needed assertion follows from Property 1 of \secref{propacyc}.

\end{proof}

In the same way as above, it is easy to prove that for any $\S\in \Sh(\Bun_G)$ and a sequence 
$\lambda_1,...,\lambda_n$
of elements of $\Lambda_G^{+}$, the sheaf
$$(H_G^{\lambda_n}\boxtimes\on{id}^{n-1})\circ...\circ(H_G^{\lambda_2}\boxtimes\on{id})\circ H_G^{\lambda_1}(\S)$$
on $\Bun_G\times X^n$ is ULA with respect to the projection $\Bun_G\times X^n\to\Bun_G$.

\medskip

Now let $\lal$ be an element of $\syminfty(\Lambda_G^+-0)$ equal to
$\{\underset{n_1\text{ times}}{\underbrace{\lambda_1,...,\lambda_1}},...,
\underset{n_k\text{ times}}{\underbrace{\lambda_k,...,\lambda_k}}\}$. Set $N=\Sum\, n_i$ and
consider the functor $\overset{o}H{}^{\lal}_G:
\Sh(\Bun_G)\to \Sh(\Bun_G\times X^{\lal})$ given by
$$\S\to (\hl_G\times\pi)_{!}\circ\hr_G{}^*(\S)[\langle|\lal|,2\rhoch\rangle+N],$$
defined using the stack $\H_G^{\lal}$.

\begin{prop} \label{goodestimate}
Let $\S$ be a ``good'' perverse sheaf on $\Bun_G$. Then $\overset{o}H{}^{\lal}_G(\S)$ lies in cohomological degrees
$\leq 0$ and is ULA with respect to the projection $\Sh(\Bun_G\times X^{\lal})\to X^{\lal}$.
\end{prop}

\begin{proof}

Let $\lal$ be equal to $\{\underset{n_1\text{ times}}{\underbrace{\lambda_1,...,\lambda_1}},...,
\underset{n_k\text{ times}}{\underbrace{\lambda_k,...,\lambda_k}}\}$. 
Consider the non-symmetrized Hecke stack $'\H^{\lal}_G$, defined as
$'\H^{\lal}_G:={}'X^{\lal}\underset{X^{\lal}}\times \H_G^{\lal}$, where
$'X^{\lal}:=X^{n_1}\times...\times X^{n_k}-\Delta$.

Let $'\hr_G$, $'\hl_G$ and $'\pi$ denote the projections from $'\H^{\lal}_G$ to $\Bun_G$ and $'X^{\lal}$,
respectively. The map $'X^{\lal}\to X^{\lal}$ is an \'etale covering. Therefore, it is enough to show that
the complex
$$'\overset{o}H{}^{\lal}_G(\S):=({}'\hl_G\times{}'\pi)_{!}\circ{}'\hr_G{}^*(\S)
[\langle|\lal|,2\rhoch\rangle+N]\in \Sh(\Bun_G\times {}'X^{\lal})$$
lives in non-positive cohomological degrees and is ULA with respect to the projection
$\Bun_G\times {}'X^{\lal}\to{}'X^{\lal}$.

Observe that the sheaf ${\Ql}_{\Gr^\lambda}[\langle\lambda,2\rhoch\rangle]$ on the affine Grassmannian is 
an extention of sheaves of the form
$\A_G^{\lambda'}[m]$, $\lambda'\leq\lambda$, $m\geq 0$. Hence, the complex $'\overset{o}H{}^{\lal}_G(\S)$ is an 
extention of complexes of the form
$(H_G^{\lambda'_N}\boxtimes\on{id}^{N-1})\circ...\circ(H_G^{\lambda'_2}\boxtimes\on{id})
\circ H_G^{\lambda'_1}(\S)[m]$
for finitely many values of $\lambda'_1,...,\lambda'_N$ and $m\geq 0$.

This proves the required assertion in view of the definition of ``good''' perverse sheaves and
the above \propref{novancycl}. 

\end{proof}

\sssec{}

Now we can prove \thmref{compat}.

\begin{proof}

It is clear that if $\S$ is a ``good'' perverse sheaf on $\Bun_M$, then so is $\DD(\S)$. Therefore,
to prove the theorem we have to show that the restriction of $\r_P{}_{!}\circ\qw_P^{!*}(\S)$ to
$\BunPb-\Bun_P$ lives in negative cohomological degrees if $S$ is ``good''. In other words,
we must show that for every $\nul\in \syminfty(\Lambda_{M,G}^+-0)$, the sheaf
$$\r_P{}_{!}\circ j_{\nul}{}_{!}\circ j_{\nul}^*\circ \qw_P^{!*}(\S)$$
lives in negative cohomological degrees.

For $\nul$ as above, let $\thl$ be the corresponding element
of $\syminfty(\Lambda_{G,P}^{\on{pos}}-0)$. The map $\r_P\circ j_{\nul}:
{}_{\nul}\Bun_P\simeq \Bun_P\underset{\Bun_M}\times\H_M^{\nul}\to {}_{\thl}\Bun_P$ is the long vertical map in
the diagram
$$
\CD
\Bun_P\underset{\Bun_M}\times\H_M^{\nul} @>>> \H_M^{\nul} @>{\hr_M}>> \Bun_M  \\
  @VVV         @V{\hl_M}VV   \\
\Bun_P\times X^{\nul}  @>{{\mathfrak q}_P\times\on{id}}>> \Bun_M\times X^{\nul}  \\
@V{\on{id}\times\on{sym}}VV \\
\Bun_P\times X^{\thl},
\endCD
$$
where the map $X^{\nul}\to X^{\thl}$ is a finite (symmetrization) map. 
Moreover, the long horizontal map in this diagram is nothing but $\qw_P\circ j_{\nul}$.

Let $K$ be an irreducible subquotient of $h^i(\IC_{\BunPbw}|_{_{\nul}\Bun_P})$ for some $i\in \NN$
(of course, $i< 0$). \corref{corsmoothonstrata} implies that every such $K$ is a pull--back
of $K'[N'+\langle|\nul|,2\rhoch_M\rangle]$, where $K'$ is a perverse sheaf on $X^{\nul}$
and $N'$ is the dimension of the corresponding connected component of $\Bun_P$.

By applying the projection formula, we obtain that 
$\r_P{}_{!}\circ j_{\nul}{}_{!}\circ j_{\nul}^*\circ \qw_P^{!*}(\S)$ is an extention of sheaves of the form
$$(\on{id}\times\on{sym})_{!}\circ ({\mathfrak q}_P\times\on{id})^*
(({\Ql}_{\Bun_M}\boxtimes K')\otimes \overset{o}H{}^{\nul}_M(\S))[-N+N'-\on{dim}(\Bun_M)+i],$$
where $N$ is as in \propref{goodestimate} above and $i>0$.

Therefore, it suffices to show that
$({\Ql}_{\Bun_M}\boxtimes K')\otimes \overset{o}H{}^{\nul}_M(\S)[-N]$ lives in non-positive cohomological degrees.
However, this follows immediately from \propref{goodestimate} and Property 5 of \secref{propacyc}.

\end{proof}


\section{The Functional equation}

\ssec{\thmref{restriction}}

In this subsection we will formulate and prove \thmref{restriction} which is one of the main results of this paper.

\sssec{}

Let us consider the following stack
$$\BunBP:=\BunPbw\underset{\Bun_M}\times\overline{\Bun}_{B(M)},$$
where $B(M)$ is the Borel subgroup of the group $M$:

$$
\CD
\BunBP @>{\qw'_P}>> \overline{\Bun}_{B(M)}  \\
@V{\p'_M}VV    @V{\p_M}VV   \\
\BunPbw  @>{\qw_P}>>  \Bun_M
\endCD
$$
and let $\CalC(\IC_{\BunPbw},\IC_{\ol{\Bun}_{B(M)}})$ denote the object
$$\p'_M{}^*(\IC_{\overline{\Bun}_{B(M)}})\otimes \qw'_P{}^*(\IC_{\BunPbw})\otimes 
(\Ql(\frac{1}{2})[1])^{\otimes-\dim(\Bun_M)}\in \Sh(\BunBP).$$

\smallskip

\begin{thm} \label{IC}
There sheaf $\CalC(\IC_{\BunPbw},\IC_{\ol{\Bun}_{B(M)}})$ is canonically isomorphic to the intersection cohomology
sheaf of $\BunBP$.
\end{thm}

\begin{proof}

Consider the following general set-up:

\medskip

Let $f:Y_1\to Y_2$ be a map of algebraic stacks with $Y_2$ smooth. Let $j:Y^0_1\hookrightarrow Y_1$ 
be an open substack
such that the map $f\circ j:Y^0_1\to Y_2$ is smooth as well. Assume, in addition, that the sheaves $\IC_{Y_1}$ and
$j_{!}(\IC_{Y^0_1})$ are ULA with respect to the map $f$.

\smallskip

Now let $Y_3$ be another algebraic stack mapping to $Y_2$. Let $Y'$ be the Cartesian product:
$$
\CD
Y' @>{f'}>> Y_3 \\
@V{g'}VV  @V{g}VV  \\
Y_1 @>f>> Y_2.
\endCD
$$

\begin{lem}
In the above situation $$\IC_{Y'}\simeq f'{}^*(\IC_{Y_3})\otimes g'{}^*(\IC_{Y_1})\otimes
(\Ql(\frac{1}{2})[1])^{\otimes-\dim(Y_2)}.$$
\end{lem}

The proof of this lemma follows immediately from the definition of the ULA property 
(cf. proof of \thmref{goodtopullbacknonprince} in \secref{proofofgoodtopullback}).

\smallskip

We apply this lemma in the situation when $Y_2:=\Bun_M$, $Y_1:=\BunPbw$, $Y_3:=\overline{\Bun}_{B(M)}$ 
and the assertion of our theorem follows from \thmref{acycthm}.

\end{proof}

\sssec{}

Observe that $\Bun_P\underset{\Bun_M}\times\Bun_{B(M)}$ is contained in 
$\BunBP$ as an open substack. There is a natural identification
$\Bun_P\underset{\Bun_M}\times\Bun_{B(M)}\simeq\Bun_B$ and the next proposition says that the
map in one direction extends to the whole of $\BunBP$:

\begin{prop}
There exists a natural map $\k_P:\BunBP\to\BunBb$, which is representable and proper. 
\end{prop}

\begin{proof}

A point of the stack $\BunBP$ is by definition a triple of bundles $(\F_G,\F_M,\F_T)$ plus a collection of
embeddings
\begin{align*}
&\widetilde\kappa_P^{\V}:(\V^{U(P)})_{\F_M} \hookrightarrow 
\V_{\F_G} \text{ for all $G$-modules } \V, \\
&\kappa_M^{\nuch}:\L^{\nuch}_{\F_T}\to \U^{\nuch}_{\F_M},\,\,\forall\nu\in\Lambda_M^{+},
\end{align*}
which satisfy the Pl\"ucker relations.

\smallskip

We set:  $\k_P(\F_G,\F_M,\F_T,\widetilde\kappa_P,\kappa_M)=(\F_G,\F_T,\kappa)$, 
where for $\lambdach\in\check\Lambda_G^{+}$, the map $\kappa^{\lambdach}$ is the composition:
$$\L_{\F_T}^{\lambdach}\overset{\kappa_M^{\lambdach}}\longrightarrow \U^{\lambdach}_{\F_M} \to
((\V^{\lambdach})^{U(P)})_{\F_M}
\overset{\widetilde\kappa_P^{\lambdach}}\longrightarrow \V^{\lambdach}_{\F_G}.$$

\end{proof}

\begin{thm} \label{restriction}
$\k_P{}_{!}(\CalC(\IC_{\BunPbw},\IC_{\ol{\Bun}_{B(M)}}))\simeq \IC_{\BunBb}$.
\end{thm}

The rest of this subsection will be devoted to the proof of this theorem.

\sssec{}

It is obvious that over $\Bun_B\subset\BunBb$, the sheaves 
$\k_P{}_{!}(\CalC(\IC_{\BunPbw},\IC_{\ol{\Bun}_{B(M)}}))$ and $\IC_{\BunBb}$ coincide. 
Since the map $\k_P$ is proper, it follows from \thmref{IC} that in order to prove \thmref{restriction},
it suffices to check that the $*$-restriction of the sheaf
$\k_P{}_{!}(\CalC(\IC_{\BunPbw},\IC_{\ol{\Bun}_{B(M)}}))$ 
to $\BunBb-\Bun_B$ lives in negative cohomological degrees. 

\medskip

Consider the set $\syminfty((\Lambda^{\on{pos}}_M\times\Lambda^{+}_{M,G})-0)$, whose elements
we will symbolically denote by $\ol{\mu,\nu}$, each being of the form
$$\{\underset{n_1\text{ times}}{\underbrace{(\mu_1\times\nu_1),...,(\mu_1\times\nu_1)}},...,
\underset{n_k\text{ times}}{\underbrace{(\mu_k\times\nu_k),...,(\mu_k\times\nu_k)}}\},$$
where in every pair $(\mu_1\times\nu_i)$ one of the elements is non-zero. We set
$|\ol{\mu,\nu}|:=\underset{i=1,...,k}{\Sum} n_i\cdot (\nu_i-\mu_i)\in\Lambda$.

Let $X^{\ol{\mu,\nu}}$ be the corresponding partially symmetrized power of $X$ and let
$\H_M^{\ol{\mu,\nu}}\overset{\pi}\to X^{\ol{\mu,\nu}}$ be the appropriate version of the Hecke stack.
The constructions of \secref{stra1} and \secref{stra2} yield naturally defined maps:
$$\Bun_{B(M)}\times X^{\ol{\mu,\nu}}\to \overline{\Bun}_{B(M)} \text{ and }
\Bun_P\underset{\Bun_M}\times \H_M^{\ol{\mu,\nu}}\to \BunPbw,$$
each being a composition of a smooth map followed by a locally closed embedding. (The smooth part appears
due to the fact that some of the $\mu_i$'s and $\nu_i$'s may vanish or coincide.)
We will denote their images by $_{\ol{\mu,\nu}}\Bun_{B(M)}$ and $_{\ol{\mu,\nu}}\Bun_P$, respectively.
Finally, we will denote by $_{\ol{\mu,\nu}}\Bun_{B,P}$ the intersection 
$$\p'_M{}^{-1}({}_{\ol{\mu,\nu}}\Bun_P)\cap \qw'_P{}^{-1}({}_{\ol{\mu,\nu}}\Bun_{B(M)})\subset \BunBP.$$
By definition, $_{\ol{\mu,\nu}}\Bun_{B,P}$ is the image of a locally closed embedding
$$\Bun_P\underset{\Bun_M}\times \H_M^{\ol{\mu,\nu}}\underset{\Bun_M}\times \Bun_{B(M)}\hookrightarrow \BunBP,$$
where the fiber product is formed using the map $(\hl_M\times\hr_M):\H_M^{\ol{\mu,\nu}}\to \Bun_M\times\Bun_M$.

Now, for $\lambda\in\Lambda_G^{\on{pos}}$, let $\BunBb^\lambda$ be the locally closed substack of
$\BunBb$ defined as a union of all $_{\lal}\Bun_B$ with $|\lal|=\lambda$. For $\ol{\mu,\nu}$ as above, we set
$_{\ol{\mu,\nu},\lambda}\Bun_{B,P}$ to be the intersection 
$_{\ol{\mu,\nu}}\Bun_{B,P}\cap \,\k_P^{-1}(\BunBb^\lambda)$.

To prove \thmref{restriction}, it suffices to check the following:

\begin{prop} \label{proprestr}
For $\ol{\mu,\nu}$ and $\lambda$ as above the following holds:

\smallskip

\noindent{\em (a)}
The $*$-restriction of $\C(\IC_{\BunPbw},\IC_{\ol{\Bun}_{B(M)}})$ to
$_{\ol{\mu,\nu},\lambda}\Bun_{B,P}$ lives in the cohomological degrees strictly smaller than 
$-\langle \lambda+|\ol{\mu,\nu}|,\rhoch_M\rangle$.

\smallskip

\noindent{\em (b)}
The fibers of the map $\k_P:{}_{\ol{\mu,\nu},\lambda}\BunBP\to {}_{\lambda}\BunBb$
are of dimension $\leq\langle \lambda+|\ol{\mu,\nu}|,\rhoch_M\rangle$.
\end{prop}

\sssec{Proof of \propref{proprestr}}

For $\ol{\mu,\nu}$ as above set $N=\Sum \, n_i$ and consider the corresponding symmetrization map
$X^N-\Delta\simeq{}'X^{\ol{\mu,\nu}}\to X^{\ol{\mu,\nu}}$. We will introduce an order in our $\ol{\mu,\nu}$
and write it as $(\mu_1,\nu_1),...,(\mu_j,\nu_j),...,(\mu_N,\nu_N)$. Let $'\H_M^{\ol{\mu,\nu}}$ be the corresponding
desymmetrized Hecke stack, i.e. its fiber over $(x_1,...,x_N)\in {}'X^{\ol{\mu,\nu}}$ is the fiber product
of $_{x_j}\H_M^{\nu_j}$, $j=1,...,N$ over $\Bun_M$, with respect to the projections 
$\hl_M:{}_{x_j}\H_M^{\nu_j}\to\Bun_M$.

For a collection of elements $\lambda_1,...,\lambda_N$ of $\Lambda_G^{+}$
let $Z^{\{\mu_j\},\{\nu_j\},\{\lambda_j\}}$ denote the following
locally closed substack of
$\Bun_P\underset{\Bun_M}\times {}'\H_M^{\ol{\mu,\nu}}
\underset{\Bun_M}\times \Bun_{B(M)}$:

A point $((x_1,...,x_N),\F_P,\F'_M,\beta_M,\F'_T,\kappa'_M)$ belongs to $Z^{\{\mu_j\},\{\nu_j\},\{\lambda_j\}}$
if for every $j$ and for every $\lambdach\in\check\Lambda_G^+$, the composition
$$\L_{\F'_T}^{\lambdach}\overset{\kappa'_M{}^{\lambdach}}\longrightarrow \U^{\lambdach}_{\F'_M}
\overset{\beta_M}\longrightarrow \U^{\lambdach}_{\F_M}$$ has a zero of order 
$\langle\lambda_j-\mu_j,\lambdach\rangle$ at $x_j$. 

It is clear that the preimage of $\BunBb^{\lambda}$ under
$$\Bun_P\underset{\Bun_M}\times {}'\H_M^{\ol{\mu,\nu}}
\underset{\Bun_M}\times \Bun_{B(M)}\to {}_{\ol{\mu,\nu}}\BunBP\overset{\k_P}\to\BunBb$$
is exactly the union of those $Z^{\{\mu_j\},\{\nu_j\},\{\lambda_j\}}$'s for which
$\Sum\, \lambda_j=\lambda$.

\smallskip

\thmref{smoothonstrata} implies that the $*$--restriction of $\CalC(\IC_{\BunPbw},\IC_{\ol{\Bun}_{B(M)}})$ to
$_{\ol{\mu,\nu}}\Bun_{B,P}$ is smooth along the fibers of the projection 
$_{\ol{\mu,\nu}}\Bun_{B,P}\to X^{\ol{\mu,\nu}}$. Therefore, to prove point (a) of the 
proposition, it suffices to show that for any $(x_1,...,x_N)\in X^N-\Delta$, the codimension of
$(\Bun_P\underset{\Bun_M}\times \Pi\,{}_{x_j}\H_M^{\nu_j}
\underset{\Bun_M}\times \Bun_{B(M)})\,\cap \,Z^{\{\mu_j\},\{\nu_j\},\{\lambda_j\}}$ inside the stack
$\Bun_P\underset{\Bun_M}\times \Pi\,{}_{x_j}\H_M^{\nu_j}
\underset{\Bun_M}\times \Bun_{B(M)}$ is $\geq \langle \lambda+|\ol{\mu,\nu}|,\rhoch_M\rangle$.

However, as in \secref{proofmainHeckeprince} we obtain that the latter codimension equals
$$\underset{j}{\Sum}\, \on{codim}(S_M^{\mu_j-\lambda_j}\cap \Gr_M^{-\w_0^M(\nu_j)},\Gr_M^{-\w_0^M(\nu_j)})$$
(cf. \secref{semiinf}) which is greater or equal than
$\underset{j}{\Sum}\, \langle \nu_j+\lambda_j-\mu_j,\rhoch_M\rangle=
\langle \lambda+|\ol{\mu,\nu}|,\rhoch_M\rangle$.

\medskip

Now let us prove point (b) of the proposition. To obtain the desired estimate on the dimension, 
it suffices to analyze the fibers of the map from $Z^{\{\mu_j\},\{\nu_j\},\{\lambda_j\}}$ to $\BunBb$. Note, that
if $Z^{\{\mu_j\},\{\nu_j\},\{\lambda_j\}}$ is non-empty, then none of the $\lambda_j$'s can be zero.
Let us denote by $\lal$ the element of $\syminfty(\Lambda^{\on{pos}}_G-0)$ corresponding to the collection
$\lambda_1,...,\lambda_N$.

We have a natural map of stacks: 
$$\phi:Z^{\{\mu_j\},\{\nu_j\},\{\lambda_j\}}\to (X^N-\Delta)\times \Bun_P\underset{\Bun_M}\times 
\Bun_{B(M)}\simeq (X^N-\Delta)\times\Bun_B,$$
that sends $$((x_1,...,x_N),\F_P,\F'_M,\beta_M,\F'_T,\kappa'_M)\, \mapsto \,
((x_1,...,x_N),\F_P,\F''_T,\kappa''_M),$$ where $\F''_T=\F'_T(\Sum\, (\lambda_j-\mu_j)\cdot x_j)$ and  
$\kappa''_M$ is the unique map that makes the square
$$
\CD
\L_{\F'_T}^{\lambdach} @>{\kappa'_M{}^{\lambdach}}>> \U^{\lambdach}_{\F'_M} \\
@VVV     @V{\beta_M}VV    \\
\L_{\F''_T}^{\lambdach}  @>{\kappa''_M{}^{\lambdach}}>> \U^{\lambdach}_{\F_M}
\endCD
$$
commute for each $\lambdach\in\check\Lambda_G^+$.

Now, the map $Z^{\{\mu_j\},\{\nu_j\},\{\lambda_j\}}\to\BunBb$ is a composition of the above map $\phi$
followed by 
$$(X^N-\Delta)\times \Bun_B\to\BunBb\to X^{\lal}\times\Bun_B\overset{j_{\lal}}\longrightarrow\BunBb.$$

Since the map $(X^N-\Delta)\to X^{\lal}$ is \'etale and $j_{\lal}$ is a locally closed embedding, it suffices
to analyze the fibers of $\phi$. However, as in
\secref{proofmainHeckeprince}, we obtain that they are isomorphic to
$\underset{j}\Pi\, (\Gr_M^{\nu_j}\cap S_M^{\lambda_j-\mu_j})$ and from \propref{dimestimate}
we obtain that their dimension is $\leq \underset{j}{\Sum}\, \langle \nu_j+\lambda_j-\mu_j,\rhoch_M\rangle=
\langle \lambda+|\ol{\mu,\nu}|,\rhoch_M\rangle$, which is what we had to prove.

\ssec{Composing Eisenstein series}

\sssec{}

At this point we are ready to prove the first part of \thmref{compos}, i.e. the existence
of the isomorphism of functors: $\Eis^G_T\simeq\Eis^G_M\circ\Eis^M_T$:

\begin{proof}

By definition, for $\S\in\Sh(\Bun_T)$, 
$$\Eis^G_M\circ\Eis^M_T(\S)\simeq \pw_P{}_{!}(\qw_P^*(\p_M{}_{!}\circ \q^{!*}_M(\S))\otimes \IC_{\BunPbw})
\otimes (\Ql(\frac{1}{2})[1])^{\otimes-\dim(\Bun_M)},$$
which, according to the projection formula, can be rewritten as
\begin{equation} \label{prcompone}
(\pw_P\circ \p'_M)_{!}((\q_M\circ \qw'_P)^*(\S)\otimes \CalC(\IC_{\BunPbw},\IC_{\ol{\Bun}_{B(M)}}))
\otimes (\Ql(\frac{1}{2})[1])^{\otimes-\dim(\Bun_T)}.
\end{equation}

However, the maps $\pw_P\circ \p'_M$ and $\p\circ \r_P$ from $\BunBP$ to $\Bun_G$ coincide and so
do the maps $\q_M\circ \qw'_P$ and $\q\circ\r_P$ from $\BunBP$ to $\Bun_T$. Hence, \eqref{prcompone}
can be identified, using once again the projection formula, with
$$\p_{!}(\q^*(\S)\otimes \r_P{}_{!}(\CalC(\IC_{\BunPbw},\IC_{\ol{\Bun}_{B(M)}})))
\otimes (\Ql(\frac{1}{2})[1])^{\otimes-\dim(\Bun_T)},$$
which, by \thmref{restriction} is the same as
$$\p_{!}(\q^*(\S)\otimes \IC_{\BunBb})\otimes (\Ql(\frac{1}{2})[1])^{\otimes-\dim(\Bun_T)}\simeq \Eis^G_T(\S).$$

\end{proof}

\sssec{}

Now, we need to show that the constructed above isomorphism of functors $\Eis^G_T\simeq\Eis^G_M\circ\Eis^M_T$
is compatible with the action of Hecke functors. We will consider the local Hecke functors, i.e.
$_xH_G$, $_xH_M$ and $_xH_T$, where $x$ is some fixed point of $X$. 

Consider the fiber product 
$$_{x,\infty}\BunBP:={}_{x,\infty}\BunPbw\underset{\Bun_M}\times {}_{x,\infty}\ol{\Bun}_{B(M)},$$
and for $\T_1\in\Sh({}_{x,\infty}\BunPbw)$, $\T_2\in \Sh({}_{x,\infty}\ol{\Bun}_{B(M)})$ we set
$$\CalC(\T_1,\T_2):=\p'_M{}^*(\T_1)\otimes \qw_P'{}^*(\T_2)\otimes (\Ql(\frac{1}{2})[1])^{\otimes-\dim(\Bun_M)}
\in \Sh({}_{x,\infty}\BunBP).$$

\begin{prop}  \label{intermedHecke}
For $\S\in\Sph_M$ there is a functorial isomorphism
$$\k_P{}_{!}\circ\CalC({}_x\Hr_{P,M}(\S,\T_1),\T_2)\simeq \k_P{}_{!}\circ \CalC(\T_1,{}_x\Hl_{B(M),M}(\S,\T_2)).$$
\end{prop}

\begin{proof}

Consider the fiber product
$$_{x,\infty}Z_{B,P}:={}
_{x,\infty}\BunPbw\underset{\Bun_M}\times{}_x\H_M\underset{\Bun_M}\times {}_{x,\infty}\ol{\Bun}_{B(M)}.$$

It classifies the data of $(\F_G,\F_M,\wt\kappa_P,\F'_M,\beta_M,\F'_T, \kappa'_M)$, with
$(\F_G,\F_M,\wt\kappa_P)\in {}_{x,\infty}\BunPbw$, $(\F_M,\F'_M,\beta_M)\in {}_x\H_M$ and
$(\F'_M,\F'_T, \kappa'_M)\in {}_{x,\infty}\ol{\Bun}_{B(M)}$. As in \secref{variantM}, there are
two projections $\hl$ and $\hr$ from this stack to $_{x,\infty}\BunBP$ that ``forget'' $\F'_M$ and $\F_M$,
respectively.

By applying the projection formula, we obtain that for $\S\in\Sph_M$, $\T_1\in\Sh({}_{x,\infty}\BunPbw)$ and
$\T_2\in \Sh({}_{x,\infty}\ol{\Bun}_{B(M)})$, 
$\k_P{}_{!}\circ\CalC({}_x\Hr_{P,M}(\S,\T_1),\T_2)$ can be written as:
\begin{equation} \label{prcomptwo}
(\k_P\circ \hr)_{!}((\p'_M\circ \hl)^*(\T_1)\otimes (\qw'_P\circ \hr)^*(\T_2)\otimes \S'),
\end{equation}
where $\S'$ is the pull-back of 
$(\S\tboxtimes {\Ql}_{\Bun_M})^l\otimes (\Ql(\frac{1}{2})[1])^{\otimes-\dim(\Bun_M)}\in \Sh({}_x\H_M)$
under the map $_{x,\infty}Z_{B,P}\to {}_x\H_M$.

However, the maps $\k_P\circ \hr$ and $\k_P\circ \hl$ from $_{x,\infty}Z_{B,P}$ to $_{x,\infty}\BunBb$
coincide. Hence, \eqref{prcomptwo} can be rewritten using the projection formula as
$$\k_P{}_{!}(\p'_M{}^*(\T_1)\otimes \hl_{!}((\qw'_P\circ \hr)^*(\T_2)\otimes \S')).$$ 
Moreover, it is easy to see that for any sheaf $\T'$ on $_{x,\infty}\BunBP$
$$\T'\otimes \hl_{!}((\qw'_P\circ \hr)^*(\T_2)\otimes \S')\simeq 
\T'\otimes \qw'_P{}^*({}_x\Hl_{B(M),M}(\S,\T_2))\otimes (\Ql(\frac{1}{2})[1])^{\otimes-\dim(\Bun_M)},$$ and therefore,
\eqref{prcomptwo} can be identified with
$$\k_P{}_{!}(\p'_M{}^*(\T_1)\otimes \qw'_P{}^*({}_x\Hl_{B(M),M}(\S,\T_2)))\otimes 
(\Ql(\frac{1}{2})[1])^{\otimes-\dim(\Bun_M)},$$
which is what we had to prove.

\end{proof}

Observe now that we have the folowing natural isomorphisms of functors:
$$\k_P\circ\CalC({}_xH^{?}_{P,G}(\S,\T_1),\T_2)\simeq 
{}_xH^?_{B,G}(\S,\k_P{}_{!}\circ\CalC(\T_1,\T_2)),\,\, \S\in\Sph_G \text{ and }$$
$$\k_P\circ\CalC(\T_1,{}_xH^{?}_{B(M),T}(\S',\T_2))\simeq
{}_xH^{?}_{B,T}(\S',\k_P{}_{!}\circ\CalC(\T_1,\T_2)),\,\, \S'\in\Sph_T,$$
where $?$= either $\leftarrow$ or $\to$. It is easy to see that
the above isomorphisms of functors are compatible in the following way: 

\begin{lem}  \label{Heckecompatibility}
For $\S\in\Sph_G$ and $\T_1$ and $\T_2$ as above, the composition
\begin{align*}
&{}_x\Hl_{B,G}(\S,\k_P{}_{!}\circ\CalC(\T_1,\T_2))\to
\k_P{}_{!}\circ \CalC({}_x\Hl_{P,G}(\S,\T_1),\T_2)\overset{\text{\thmref{mainHeckenonprinceG}}}\longrightarrow \\
&\k_P{}_{!}\circ \CalC({}_x\Hr_{P,M}(\on{gRes}^G_M(\S),\T_1),\T_2)
\overset{\text{\propref{intermedHecke}}}\longrightarrow
\k_P{}_{!}\circ \CalC(\T_1,{}_x\Hl_{B(M),M}(\on{gRes}^G_M(\S),\T_2)) \\
&\overset{\text{\thmref{mainHeckenonprinceG}}}\longrightarrow 
\k_P{}_{!}\circ \CalC(\T_1,{}_x\Hr_{B(M),T}(\on{gRes}^G_T(\S),\T_2))\to 
{}_x\Hr_{B,T}(\on{gRes}^G_T(\S),\k_P{}_{!}\circ\CalC(\T_1,\T_2))
\end{align*}
equals
$${}_x\Hl_{B,G}(\S,\k_P{}_{!}\circ\CalC(\T_1,\T_2))\overset{\text{\thmref{mainHeckenonprinceG}}}\longrightarrow
{}_x\Hr_{B,T}(\on{gRes}^G_T(\S),\k_P{}_{!}\circ\CalC(\T_1,\T_2)).$$
\end{lem}

This readily implies what we need. Let $\T$ be an object of $\Sh(\Bun_T)$. Then, as in the proof of 
\thmref{commutewithHeckenonprince} given in \secref{prcomheckenonprince}, by taking $\T_1=\IC_{\BunPbw}$,
$\T_2=\IC_{\ol{\Bun}_{B(M)}}$, each of the two isomorphisms
of functors of \lemref{Heckecompatibility}
produces an isomorphism $_x\Hl_G(\S,\Eis^G_T(\T))\simeq \Eis^G_T({}_x\Hl_T(\on{gRes}^G_T(\S),\T))$:

For the latter, this will be the isomorphism of \thmref{commutewithHecke} and for the former, this
will be the composition:
\begin{align*}
&_x\Hl_G(\S,\Eis^G_T(\T))\simeq {}_x\Hl_G(\S,\Eis^G_M\circ \Eis^M_T(\T))\simeq 
\Eis^G_M({}_x\Hl_M(\on{gRes}^G_M(\S),\Eis^M_T(\T)))\simeq \\
&\Eis^G_M\circ\Eis^M_T({}_x\Hl_T(\on{gRes}^G_T(\S),\T))\simeq \Eis^G_T({}_x\Hl_T(\on{gRes}^G_T(\S),\T)).
\end{align*}
Thus, the assertion of \lemref{Heckecompatibility} implies the second part of \thmref{compos}.

\ssec{Computation in rank $1$}  \label{comprkone}

\sssec{}

From this moment on we will be concerned with the proof of \thmref{functeq}. Observe that
\thmref{compos} reduces the assertion of \thmref{functeq} to the case when $G$ is of semi-simple
rank $1$:

Indeed, by induction we may suppose that our element $\w$ is a simple reflection $s_{\i}$ corresponding
to $\i\in\I$, and according to \thmref{compos} it is enough to construct an isomorphism
$\Eis^{M_{\i}}_T(s_{\i}\cdot\S)\simeq\Eis^{M_{\i}}_T(\S)$, where $M_{\i}$ is the corresponding
sub-minimal Levi, i.e. $\I_{M_{\i}}=\{\i\}$.

Henceforth, $G$ will be a reductive group of semi--simple rank $1$. The (only) positive root, fundamental
weight and the non-trivial Weyl group elements will be denoted by $\alphach$, $\omegach$ and $\sigma$,
respectively. Note that the Weyl module $\V^{\omegach}$ is $2$-dimensional and irreducible.

Recall that the isomorphism  
$\Eis(\sigma\cdot \S)\overset{f.eq}\longrightarrow\Eis(\S)$ (for $\S\in\Sh(\Bun_T)^{reg}$) should depend on a 
specific choice of a lift of $\sigma$ to $N(\check T)$. The isomorphism that we will construct corresponds 
to the image of
$\begin{pmatrix} 0 & 1 \\ -1 & 0 
\end{pmatrix} \in SL(2)$
under the canonical map $SL(2)\to\check G$. If we multiply this element on the left by $\tau\in \check T$, 
the corresponding isomorphism must be multiplied by 
$\mu(\tau)$ for $\S$, which is concentrated on $\Bun_T^{\mu}$.

\sssec{}

Consider the fiber product $\Bun_{G,T}:=\Bun_G\underset{\Bun_{G/[G,G]}}\times\Bun_T$. We will denote
by $p_G$ and $p_T$ the projections from $\Bun_{G,T}$ to $\Bun_G$ and $\Bun_T$, respectively.
The starting point of our discussion is the following observation:

\begin{prop}
For $G$ as above, the stack $\BunBb$ classifies triples $(\F_G,\F_T,\kap)$, where $(\F_G,\F_T)\in\Bun_{G,T}$
and $\kap$ is an injective (=non-zero) map $\L^{\omegach}_{\F_T}\to \V^{\omegach}_{\F_G}$. In particular,
$\BunBb$ is smooth.
\end{prop}

\begin{proof}

Since $\omegach$ is the only fundamental weight, $\BunBb$ is a evidently a closed substack in the stack
classifying triples $(\F_G,\F_T,\kap)$ as above. To show that this closed embedding is an isomorphism,
one has to show that any map $\kap$ satisfies the Pl\"ucker relations, which is evident.

To show that $\BunBb$ is smooth we argue as follows: consider the algebraic stack $Coh_1$ which classifies
coherent sheaves on $X$ of generic rank $1$. It is known (cf. \cite{La}) that $Coh_1$ is smooth.

We have a natural map $\BunBb\to Coh_1$, which sends a triple $(\F_G,\F_T,\kap)$ as above to
the quotient $\V^{\omegach}/\on{Im}(\kap)$. It is easy to see that this map is smooth, hence $\BunBb$
is smooth too.

\end{proof}

Now, let $\BunBb^r$ be the stack that classifies triples $(\F_G,\F_T,\kappa)$, where 
$(\F_G,\F_T)\in\Bun_{G,T}$ as before, but $\kap$ is allowed to be an arbitrary map
$\L^{\omegach}_{\F_T}\to \V^{\omegach}_{\F_G}$. We will denote the natural projections from 
$\BunBb^r$ to $\Bun_G$ and $\Bun_T$ by $\p^r$ and $\q^r$, respectively.

By construction, $\BunBb$ is an open substack in $\BunBb^r$, which corresponds to non-zero $\kap$'s.
The complement, i.e. the locus $\kap=0$ projects isomorphically onto $\Bun_{G,T}$.

\medskip

For $\S\in\Sh(\Bun_T)$, we define its renormalized Eisenstein series as
$$\Eis^r(\S):=\p^r_{!}\circ \q^r{}^*(\S)\otimes (\Ql(\frac{1}{2})[1])^{\otimes ?},$$
where $?$ equals $g-1+\langle \mu,\rhoch\rangle$ over the $\mu$-th connected component of $\Bun_T$.

Note that the map $\p^r$ is non-representable, since $\Bun_T$ is a stack and not a scheme.
Therefore, $\Eis^r$ takes values not in $\Sh(\Bun_G)$, but rather in the corresponding 
``unbounded on the left'' derived category. 

\begin{prop} \label{renorm}
For $\S\in\Sh(\Bun_T)^{reg}$, there is a canonical isomorphism
$\Eis(\S)\simeq\Eis^r(\S)$.
\end{prop}

\begin{proof} 

By construction, we have a natural map $\Eis(\S)\to\Eis^r(\S)$, which corresponds to the
open embedding $\BunBb\hookrightarrow\BunBb^r$. Its cone is by definition 
$p_G{}_!\circ p_T^*(\S)$. However,the latter complex is the same as the pull-back under
$\Bun_G\to \Bun_{G/[G,G]}$ of $\f_{!}(\S)\otimes (\Ql(\frac{1}{2})[1])^{\otimes ?}$ (cf. \secref{intreq}).
This implies our assertion.

\end{proof}

Therefore, in order to construct an isomorphism $\Eis(\sigma\cdot\S)\overset{f.eq}\longrightarrow \Eis(\S)$ for an
object $\S\in\Sh(\Bun_T)^{reg}$, it suffuces to construct a functorial isomorphism
\begin{equation} \label{renequation}
\Eis^r(\sigma\cdot\S)\to \Eis^r(\S),\,\,\forall \S\in\Sh(\Bun_T).
\end{equation}

This will be done in the following general set-up:

\sssec{}  \label{developfourier}

Let $\Y$ be a stack and let $\E_0$ and $\E_1$ be two locally free coherent sheaves on $\Y$ and let 
$f:\E_0\to \E_1$ be a map between them; we will regard $\E_0\to \E_1$ as a length--$2$ complex on $\Y$; 
let $\chi(\E)$ be its Euler characteristic, i.e. $\chi(\E)=\on{rk}(\E_0)-\on{rk}(\E_1)$. We will denote
by $\E^*$ the dual complex, in other words, $\E^*$ corresponds to the adjoint map 
$\E_1^*\overset{f^*}\to \E_0^*$. 

Let $K_\E$ denote the following object of $\Sh(\Y)$: we can regard 
$\E_0$ and $\E_1$ as 
group--schemes over $\Y$, and let 
$h^0(\E)\subset \E_0$ denote the group--scheme--theoretic kernel of $f$. Then 
$K_\E$ is the !-direct image of ${\Ql}_{h^0(\E)}$ on $\Y$ tensored
by $(\Ql(\frac{1}{2})[1])^{\otimes \chi(\E)}$.

\begin{lem} \label{Fourier}
Under the above circumstances, there is a canonical isomorphism \newline
$K_\E\overset{eq_\E}\longrightarrow K_{\E*}$.
\end{lem}

\begin{proof}

Indeed, \footnote{The construction described below is due to V.~Drinfeld.}
let us consider the fiber product $\E_0\underset{\Y}\times \E_1^*$. We have the evaluatuion map
$\E_0\underset{\Y}\times \E_1^* \to \AA^1$ and let us denote by $\Psi$ the 
pull-back of the Artin-Schreier sheaf on $\AA^1$ under this map, tensored by
$(\Ql(\frac{1}{2})[1])^{\on{rk}(\E_0)}\boxtimes (\Ql(\frac{1}{2})[1])^{\on{rk}(\E^*_1)}$.

To construct the isomorphism $eq_\E$, it is enough to show that both 
$K_\E$ and $K_{\E^*}$ can be canonically
identified with the $!$--direct image of $\Psi$ under the projection $\E_0\underset{\Y}\times \E_1^*\to \Y$.

However, the !-direct image of $\Psi$ under the projection 
$\E_0\underset{\Y}\times \E_1^*\to \E_0$ is evidently isomorphic to
$(\Ql(\frac{1}{2})[1])^{\otimes \chi(\E)}_{h^0(\E)}$ and hence is further push-forward onto $\Y$
can be identified with $K_\E$. The assertion for $K_{\E^*}$ follows by symmetry.

\end{proof}

The following properties of the isomorphism $eq_\E$ follow from the construction:

\begin{lem}  \label{propfourier}

\smallskip

\noindent{\em (a)} The composition $eq_{\E^*}\circ eq_\E: K_\E\to K_\E$ is the identity map.

\smallskip

\noindent{\em (b)}
Let $f':\E'_0\to \E'_1$ be another $2$--complex consisting of locally free sheaves and let
$(\E_0\overset{f}\to \E_1)\to (\E'_0\overset{f'}\to \E'_1)$ be a quasi-isomorphism. 
Then there are canonical isomorphisms
$h^0(\E)\simeq h^0(\E')$ and $h^0(\E^*)\simeq h^0(\E'{}^*)$ and $eq_\E$ coincides with $eq_{\E'}$ up to 
$(-1)^{\on{rk}(\E_1)-\on{rk}(\E'_1)}$.

\end{lem}

\sssec{}

We will apply the above discussion to $\Y=\Bun_{G,T}$. First, let us consider the $2$-dimensional vector 
bundle $\M$ over
$X\times \Bun_{G,T}$, whose fiber over $(\F_G,\F_T)\in \Bun_{G,T}$ is 
$\underline\Hom(\L^{\omegach}_{\F_T},\V^{\omegach}_{\F_G})$. Set $\E$ to be the derived direct image of $\M$ 
with respect to the projection $X\times \Bun_{G,T}\to \Bun_{G,T}$.

It is clear that locally on $\Bun_{G,T}$, we can $\E$ represent by a $2$--complex $\E_0\to \E_1$
of vector bundles with a fixed parity of $\on{rk}(\E_1)$:
Indeed, let $y\in X$ be an arbitrary point and let $n$ be a large postive integer. Let $S_{y,n}$ be an open
substack of $\Bun_{G,T}$ which corresponds to those $(\F_G,\F_T)$ for which
$\on{Ext}^1(\L^{\omegach}_{\F_T},\V^{\omegach}_{\F_G}(n\cdot y))=0$.

Then over a point $(\F_G,\F_T)\in S_{y,n}$ we set
$$\E_0=\Hom(\L^{\omegach}_{\F_T},\V^{\omegach}_{\F_G}(n\cdot y)) \text{ and }
\E_1=\Hom(\L^{\omegach}_{\F_T},\V^{\omegach}_{\F_G}(n\cdot y)/\V^{\omegach}_{\F_G}).$$

Thus, we obtained a correctly defined object $K_\E\in\Sh(\Bun_{G,T})$. It is easy to see that $h^0(\E)\simeq \BunBb^r$
and hence $\K_E$ is the {\it kernel} of the renormalized Eisenstein series functor $\Eis^r$, i.e.
$$\Eis^r(S)\simeq p_G{}_!(p_T^*(\S)\otimes \K_E).$$

\medskip

Consider now the action of the involution $\sigma\cdot$
on $\Bun_{G,T}$. We obtain from 
\lemref{Fourier} and \lemref{propfourier}(b)
that in order to construct the isomorphism of \eqref{renequation},
it suffices to show that $\sigma\cdot\K_\E\simeq \K_{\E^*}$. 
For that, by Serre's duality, it suffices to show that over $X\times \Bun_{G,T}$
$$\sigma\cdot\M\simeq \M^*\otimes(\Omega_X\boxtimes \O_{\Bun_{G,T}}),$$ 
where $\M^*$ is the dual of $\M$.

By definition, the fiber of $\sigma\cdot\M$ over $(\F_G,\F_T)\in\Bun_{G,T}$ is
$$\Hom(\L^{\omegach-\alphach}_{\F_T}\otimes\Omega_X^{-1},\V^{\omegach}_{\F_G})\simeq
\Hom(\L_{\F_T}^{-\omegach},\L_{\F_T}^{-2\omegach+\alphach}\otimes\V_{\F_G}\otimes\Omega_X).$$
However, since $\dim(\V^{\omegach})=2$, 
$\V_{\F_G}^{\omegach}\simeq (\V_{\F_G}^{\omegach})^*\otimes\on{det}(\V^{\omegach}_{\F_G})
\simeq (\V_{\F_G}^{\omegach})^*\otimes\L_{\F_T}^{2\omegach-\alphach}$.

Therefore, the above fiber of $\sigma\cdot\M$ can be identified with 
$\Hom((\L_{\F_T}^{\omegach})^*,(\V_{\F_G}^{\omegach})^*\otimes\Omega_X)$, which is what we had to show.

\sssec{}

The above construction yields an isomorphism of functors $\Eis(\sigma\cdot\S)\overset{f.eq'}\longrightarrow\Eis(\S)$.
Our task now is to modify it in such a way that
it becomes compatible with the isomorphism of \thmref{commutewithHecke} in the sense of \secref{compfuncteq}
for the specified above choice of a lifting of $\sigma$ to $N(\check T)$.

\medskip

\noindent{\it Remark.}
There is another way to obtain the functional equation isomorphism, using the geometric analog
of the Whittaker model (cf. \cite{Ga}). By comparing the two, one can show that the isomorphism $f.eq'$
that we constructed corresponds to the element $\begin{pmatrix} 0 & 1 \\ -1 & 0 
\end{pmatrix}\cdot \epsilon \in SL(2)$, where $\epsilon=\pm 1$, depending on the genus of $X$.

\medskip

Let us first consider the case of $G=GL(2)$, in which case the Langlands dual group is $GL(2)$ as well.
For two integers $d_1,d_2$, we will denote by $(d_1,d_2)$ the corresponding coweight of $GL(2)$. 

\begin{prop}  \label{wontprove}
For $\S\in\Sh(\Bun^\mu_T)$ with $\mu=(d_1,d_2)$, there is a functorial isomorphism
$H_G^{(1,0)}(\Eis^r(\S))\simeq \Eis^r(H_T^{(1,0)}(\S))\oplus \Eis^r(H_T^{(0,1)}(\S))$,
such that:

\smallskip

\noindent{\em (a)} For $\S\in\Sh(\Bun_T)^{reg}$, it coincides with the isomorphism of \thmref{commutewithHecke}.

\smallskip

\noindent{\em (b)} The square
$$
\CD
_xH_G^{(1,0)}(\Eis^r(\sigma\cdot\S))   @>>> \Eis^r(\sigma\cdot({}_xH_T^{(1,0)}(\S)))\oplus
\Eis^r(\sigma\cdot({}_xH_T^{(0,1)}(\S)))   \\
@V{f.eq'}VV    @V{-f.eq'\cdot \epsilon(d_1,d_2) \oplus f.eq'\cdot \epsilon(d_1,d_2)}VV      \\
_xH_G^{(1,0)}(\Eis^r(\S))   @>>>   \Eis^r({}_xH_T^{(1,0)}(\S))
\oplus \Eis^r({}_xH_T^{(0,1)}(\S))
\endCD
$$
commutes, where $\epsilon(d_1,d_2)=\pm 1$. Moreover, $\epsilon(d_1+d,d_2+d)=\epsilon(d_1,d_2)$.
\end{prop}

\begin{proof}

Consider the fiber product $Z:={}_x\H_G^{(1,0)}\underset{\Bun_G}\times \BunBb^r$, where $_x\H_G^{(1,0)}\to\Bun_G$
is the map $\hr_G$. By definition, this stack classifies quintuples $(\M,\M',\L',\kap',\beta)$, where
$\M,\M'$ are rank $2$-bundles, $\beta:\M'\hookrightarrow \M$ is an embedding of coherent sheaves such that
$\M/\M'$ is a skyscraper sheaf at $x$, $\L'$ is a line bundle and $\kap':\L'\to\M'$ is an arbitrary map.

\smallskip

We have a natural map $\phi:Z\to \BunBb^r$, which
sends a quintuple $(\M,\M',\L',\kap',\beta)$ as above to $(\M,\L,\kap)$, where $\L=\L'$ and $\kap$
is the composition 
$$\L=\L'\to \M'\to \M.$$

Let $_{x,\geq (1,-1)}\BunBb^r\subset \BunBb^r$ we the closed substack of $\BunBb^r$, which corresponds to
triples $(\M,\L,\kap)$ such that $\kap:\L\to \M(-x)$. To prove point (a) of the proposition, it is enough to show
that 
$$\phi_!(\Ql_{Z})\simeq \Ql_{\BunBb^r}\oplus \Ql_{_{x,\geq (1,-1)}\BunBb^r}(-1)[-2].$$

\noindent{\it Remark.}
If we knew that $\Ql_{\BunBb^r}$ is $\IC_{\BunBb^r}$ up to cohomological shift and Tate's twist,
the last asserion would be obvious from the decomposition theorem (cf. below). Since the latter
fact is not in general true, we have to proceed differently.

\medskip

Take $y\neq x$ and let $S_{y,n}\subset \Bun_{G,T}$ be as above. Let $f:\E_0\to \E_1$ be as above as well,
and let $\delta_{\E_1}$ be the constant sheaf on the $0$-section of $\E_1$. Consider 
$Z_{\E_0}:={}_x\H_G^{(1,0)}\underset{\Bun_G}\times \E_0$; as above, we have a natural map 
$\phi_{\E_0}: Z_{\E_0}\to \E_0$. Let $_{x,\geq (1,-1)}\E_0\subset \E_0$ be the closed substack that classifies
triples $(\L,\M,\L\to \M(n\cdot y))$, for which $\on{Im}(\L)\subset \M(n\cdot y-x)$. 

It is clear that over the open subset $\E_0-{}_{x,\geq (1,-1)}\E_0$, the map $\phi_{\E_0}$ is an isomorphism
and over $_{x,\geq (1,-1)}\E_0$ this is a fibration with the typical fiber ${\mathbb P}^1$. Since $\E_0$ is smooth,
from the decomposition theorem we obtain that
$$\phi_{\E_0}{}_!(\Ql_{Z_{\E_0}})\simeq \Ql_{\E_0}\oplus \Ql_{_{x,\geq (1,-1)}\E_0}(-1)[-2].$$

By construction, $\Ql_{\BunBb^r}\simeq \Ql_{\E_0}\otimes f^*(\delta_{\E_1})$. Moreover, 
since $x\neq y$, $\E_1$ does not "distinguish" between $\M$ and $\M'$ and we obatin:
$$\phi_!(\Ql_{Z})\simeq \phi_{\E_0}{}_!(\Ql_{Z_{\E_0}}) \otimes f^*(\delta_{\E_1}).$$

Therefore,
\begin{align*}
&\phi_!(\Ql_{Z})\simeq \Ql_{\E_0}\otimes f^*(\delta_{\E_1})\oplus \Ql_{_{x,\geq (1,-1)}\E_0}(-1)[-2]\otimes 
f^*(\delta_{\E_1})\simeq \\
&\Ql_{\BunBb^r}\oplus \Ql_{_{x,\geq (1,-1)}\BunBb^r}(-1)[-2],
\end{align*}
which is what we had to prove.

The finishes the proof of point (a) of the proposition. Point (b) is a tedious but straightforward
verification, which we omit.

\end{proof}

Now, for our group $G$ let $\mu\mapsto d(\mu)$ be the map $\Lambda\to \Lambda_{G/Z(G)}\simeq\ZZ$.
Let us choose a function $\epsilon':\ZZ\to \pm 1$, so that $\epsilon'(d+1)\cdot \epsilon'(d)=\epsilon(d,0)$,
where $\epsilon$ is as in the above proposition. We set the functional equation isomorphism 
$\Eis(\sigma\cdot\S)\overset{f.eq}\longrightarrow\Eis(\S)$ to be $f.eq'\cdot \epsilon'(d(\mu))$ for
$\S\in \Sh(\Bun_T^\mu)$. Let us show that it satifies the required compatibility condition with respect
to the isomorphism of \thmref{commutewithHecke}.

\medskip 

First, let us observe that if $G_1\to G_2$ is an embedding of groups whose cokernel is a torus,
then the statement that we have to check for $G_1$ follows from that for $G_2$. It is easy to see
that any $G$ with $[G,G]\simeq SL(2)$ can be embedded into a product of $GL(2)$ and a torus. This
allows to replace our initial $G$ by $GL(2)$.

For $GL(2)$ every irreducible module $V$ over the Langlands dual group can be embedded into
a tensor power of the tautologiocal representation $V^{(1,0)}$ times a $1$-dimensional representation.
Now, since the isomorphism of \thmref{commutewithHecke} is compatible with tensor products of 
$\check G$-representaions and is evident for $1$-dimensional representations, it suffices to treat 
the case of the Hecke functor $H_{GL(2)}^{(1,0)}$. In the latter case the assertion follows from \lemref{wontprove}.

\section*{Appendix}

In this Appendix we will prove \propref{prasad}. The proof presented below is
a variant of an argument that was explained to us by G.~Prasad. 

Let $\check G$ be a connected reductive algebraic group over an
algebraically closed field of characteristic $0$ and let
$\check T$ be a Cartan subgroup of $\check G$. Denote by 
$N(\check T)$ the normalizer of $\check T$ in $\check G$.  Suppose that we are
given a closed  subgroup $\Gamma\subset \check T$ and and
element $F\in \check G$ such that $F$ normalizes $\Gamma$.
Let $\wt\Gamma$ be the minimal closed subgroup of $\check G$ which
contains $\Gamma$ and $F$. 

\medskip

\noindent{\bf Theorem A.1.}
{\it Assume that $\wt{\Gamma}$ is irreducible, i.e. there is no
proper parabolic subgroup $\check P$ in $\check G$ such that 
$\wt\Gamma\subset \check P$. Then

\smallskip

\noindent{\em (a)} For every root $\check\alpha:\check T\to{\mathbb G}_m$, the restriction of
$\alphach$ to $\Gamma$ is not identically equal to 1.\smallskip

\noindent{\em (b)} $F$ normalizes $\check T$.}

\medskip

\noindent{\bf Proposition A.2.}
{\it Let $Z$ denote the centralizer of $\Gamma$ and
let $Z^0$ denote the connected component of the identity in
$Z$. Then $Z^0$ is commutative.}

\medskip

Let us show first that the above proposition imples the Theorem A.1.

\begin{proof}

Assume that (a) does not hold,
i.e. that there exists a root $\alphach$ of $\check G$ such that
$\alphach|_{\Gamma}=1$. Let $p_{\alphach}:SL(2)\to \check G$ denote the 
corresponding map. Then it is easy to see that $\Gamma$ commutes with 
the image of $p_{\alphach}$. Hence $Z^0$ contains 
$\text{Im}(p_{\alphach})$ and therefore it is not commutative.

To prove (b) let us note that $\check T\subset Z^0$. Hence, if $Z^0$ is
commutative, it follows that $Z^0=\check T$. On the other hand, $F$
clearly normalizes $Z^0$. Hence $F\in N(\check T)$.

\end{proof}

Now let us prove Proposition A.2. The proof will use the following well-known lemma:

\medskip

\noindent{\bf Lemma A.3.}
{\it Let $H$ be an algebraic group (over an algebraically
closed field of characteristic zero) and let 
$\sigma$ be an automorphism of $H$. Then there exists a $\sigma$--stable
Borel subgroup in $H$.}

\begin{proof}(of Proposition A.2.)

Assume that $Z^0$ is not commutative. Let ${\mathfrak z}$
denote its Lie algebra. We claim that there exists a nilpotent
element $n\in {\mathfrak z}, n\neq 0$, such that
$\text{ad}_{F}(n)=c\cdot n$ for a scalar $c$.

Indeed, $\text{ad}_{F}$ defines an automorphism of $Z^0$
(it is obvious that $Z$ and hence also $Z^0$ is normalized by
$F$). Therefore, the above lemma shows that there
exists a Borel subgroup $B_{Z^0}\subset Z^0$, which
is normalized by $F$. Let ${\mathfrak b}_{Z^0}$ be its Lie algebra and let
${\mathfrak n}_{Z^0}\subset {\mathfrak b}_{Z^0}$ be the nilpotent radical 
of this Lie algebra. Then ${\mathfrak n}_{Z^0}$ is also
normalized by $F$. We choose $n\in {\mathfrak n}_{Z^0}$
to be an eigenvector of $\text{ad}_{F}$.

\medskip

Thus, let $n\in {\mathfrak z}$ be as above. It is well-known that to 
every nilpotent element $n$ in $\check{\mathfrak g}$ one
can canonically associate a parabolic subgroup $\check P_n$ in $\check G$.
The subgroup $\check P_n$ can be uniquely characterized in the following
way:

Let $\check{\mathfrak p}_n\subset \check{\mathfrak g}$
be the Lie algebra of $\check P_n$. Fix a homomorphism
$p:{\bf sl}(2)\to \check{\mathfrak g}$ in such a way that
$p
\begin{pmatrix}
0& 1\\
0& 0
\end{pmatrix}
=n.$
Let $\check{\mathfrak g}=\bigoplus\limits_{i\in\ZZ} \check{\mathfrak g}_i$
be the weight decomposition of $\check{\mathfrak g}$ with
respect to the adjoint action of $p({\bf sl}(2))$. Then 
$\check{\mathfrak p}_n=\bigoplus\limits_{i\geq 0}\check{\mathfrak g}_i$.

On the one hand, the fact that $n\neq 0$ implies that 
$\check P_n\neq \check G$. On the other hand, since $n$ is an eigenvector of 
$\wt{\Gamma}$ and
the assignment $n\mapsto \check P_n$ is canonical, 
the group $\check P_n$ is also normalized by $\wt{\Gamma}$.  But as
$\check P_n$ is a parabolic subgroup in $\check G$, it follows that
$\wt{\Gamma}\subset \check P_n$, which is a contradiction.

\end{proof}

\enddocument